\documentclass[10pt,leqno]{article}

\usepackage{amssymb}
\usepackage{amsmath}
\usepackage{amsfonts}
\usepackage{amscd}
\usepackage{amsthm}
\usepackage{graphicx}
\usepackage{enumerate}
\usepackage{empheq}
\usepackage{xcolor}

\newtheorem{theorem}{Theorem}[section]

\newtheorem{Remark}{Remark}[section]
\newtheorem{lemma}{Lemma}[section]
\newtheorem{example}{Example}[section]

\newcommand{\se}{\setcounter{equation}{0}}

\newcommand{\f}{{\bf {f}}}
\newcommand{\R}{\mathbb{R}}

\newcommand{\bu}{\textbf{u}}
\newcommand{\bU}{\textbf{U}}
\newcommand{\bW}{\textbf{W}}

\newcommand{\bV}{\textbf{V}}
\newcommand{\bS}{\textbf{S}}
\newcommand{\bR}{\textbf{R}}
\newcommand{\bP}{\textbf{P}}
\newcommand{\bH}{\textbf{H}}
\newcommand{\bL}{\textbf{L}}
\newcommand{\bv}{\mathbf{v}}
\newcommand{\bw}{\textbf{w}}
\newcommand{\de}{\textbf{e}}
\newcommand{\bn}{\textbf{n}}
\newcommand{\bz}{\textbf{z}}
\newcommand{\bq}{\textbf{q}}
\newcommand{\bJ}{\textbf{J}}
\newcommand{\brf}{\mathbf{f}}

\newcommand{\e}{{\bf e}}

\DeclareMathOperator*{\esssup}{ess\: sup}

\newcommand{\bzeta}{\mbox{\boldmath $\zeta$}}

\newcommand{\bphi}{\mbox{\boldmath $\phi$}}
\newcommand{\brho}{\mbox{\boldmath $\rho$}}

\newcommand{\bxi}{\mbox{\boldmath $\xi$}}

\newcommand{\al}{\alpha}
\newcommand{\be}{\mbox{\boldmath$\eta$}}

\newcommand{\vertiii}[1]{{\left\vert\kern-0.25ex\left\vert\kern-0.25ex\left\vert #1 
		\right\vert\kern-0.25ex\right\vert\kern-0.25ex\right\vert}}

\begin{document}
	\author{Saumya Bajpai\footnotemark[1]\thanks{School of Mathematics and Computer Science, Indian Institute of Technology Goa, Ponda, Goa-403401, India. Email: saumya@iitgoa.ac.in
},~ Deepjyoti Goswami\footnotemark[2]\thanks{Department of Mathematical Sciences, Tezpur University, Tezpur, Sonitpur, Assam-784028, India. Email: deepjyoti@tezu.ernet.in, kallol@tezu.ernet.in},~ and Kallol Ray\footnotemark[2]}
	\title{
	 A Priori Error Estimates of a Discontinuous Galerkin Finite Element Method for the Kelvin-Voigt Viscoelastic Fluid Motion Equations} 
	\date{}
	\maketitle

	\begin{abstract}
		This paper applies a  discontinuous Galerkin finite element method to the Kelvin-Voigt viscoelastic fluid motion equations when the forcing function is in $L^\infty(\bL^2)$-space.  Optimal {\it a priori} error estimates in $L^\infty(\bL^2)$-norm for the velocity and in $L^\infty(L^2)$-norm for the pressure approximations for the semi-discrete discontinuous Galerkin method are derived here. The main ingredients for establishing the error estimates are the standard elliptic duality argument and a modified version of the Sobolev-Stokes operator defined on appropriate broken Sobolev spaces. Further, under the smallness assumption on the data, it has been proved that these estimates are valid uniformly in time. Then, a first-order accurate backward Euler method is employed to discretize  the semi-discrete discontinuous Galerkin Kelvin-Voigt formulation completely.  The fully discrete optimal error estimates for the velocity and pressure are established. Finally, using the numerical experiments, theoretical results are verified. It is worth highlighting here that the error results in this article for the discontinuous Galerkin method applied to the Kelvin-Voigt model using finite element analysis are the first attempt in  this direction.  
		
	\end{abstract}

	\vspace{1em} 
	
	\section{Introduction}\label{s1}
	\se
	Consider the following system of partial differential equations arising in the Kelvin-Voigt model of
	viscoelastic fluid flow 
	\begin{eqnarray}
		\label{8.1}
		\frac {\partial \bu}{\partial t}+ \bu\cdot \nabla \bu - \kappa \Delta \bu_t 
		-\nu \Delta \bu  + \nabla p
		=\f(x,t),\,\,\, x\in \Omega ,\,\,t>0,
	\end{eqnarray}
	and incompressibility condition
	\begin{eqnarray}
		\label{8.2}
		\nabla \cdot \bu=0,\,\,\,x\in \Omega,\,t>0,
	\end{eqnarray}
	with initial and boundary conditions
	\begin{eqnarray}
		\label{8.3}
		\bu(x,0)= \bu_0 \;\;\;\mbox {in}\;\Omega,\;\;\;\;\; \bu=0,\;\; \;
		\mbox {on}\; \partial \Omega,\; t\ge 0,
	\end{eqnarray}
		where, $\Omega$ is a bounded convex polygonal or polyhedral  domain in $\mathbb {R}^d, d=2,3$  with boundary $\partial \Omega.$  
 Here, $\nu$ is the coefficient of kinematic viscosity, $\kappa$ is the retardation time or the time 
	of relaxation of deformations, $\bu=(u_1,\,u_2)$ (or $(u_1, u_2, u_3)$) is the fluid velocity, $p$ is the pressure and $\brf$ is the external force. The Kelvin-Voigt model was introduced by Oskolkov \cite{AP2} to represent the dynamics of the viscoelastic fluid motion. Later, 
 Cao {\it et al.} in \cite{CLT-06} viewed it as a smooth and inviscid regularization of the Navier-Stokes model. One may refer to \cite{bs05}, \cite{bs06}, \cite{css02} and literature therein for more detailed physical description and applications 
	of the model. Using the proof techniques of Ladyzenskaya presented in \cite{L69}, Oskolkov and his collaborators carried out a considerable amount of work in \cite{AP2}, \cite{OR2} and literature referred to therein, related to the global existence of a unique "almost"  classical solution of the Kelvin-Voigt system of equations (\ref{8.1})-(\ref{8.3}) with varying assumptions on the right-hand side force function $\f$.\\  
	\noindent
	There is significant literature devoted to the numerical approximations of the problem (\ref{8.1})-(\ref{8.3}). 
		In \cite{bnpdy}, the regularity estimates for the weak solution of (\ref{8.1})-(\ref{8.3}) with the right-hand side function $\f=0$ and the exponential decay property for the weak solution have been established. The authors defined a semidiscrete Galerkin finite element formulation for the Kelvin-Voigt model and have obtained optimal error estimates for the velocity in 
		$L^{\infty}(\bL^2)$ and $L^{\infty}(\bH_0^1)$- norms and 
		pressure in $L^{\infty}(L^2)$-norm with the initial data $\bu_0\in\bH^2\cap \bH_0^1$, which again preserve the exponential decay property.  They have further extended the analysis in \cite{ANS2} by applying a first-order accurate backward Euler method and a second order backward difference scheme to the semidiscrete model and have established fully discrete optimal error estimates for both velocity and pressure approximations. The finite element analysis of (\ref{8.1})-(\ref{8.3})  with right-hand side $\f\neq 0$ can be found in \cite{BP19} and literature referred in. 

In \cite{PDY}, Pani et al. have derived $L^{\infty}(\bL^2)$ and $L^{\infty}(\bH_0^1)$-norms optimal error estimates for the spectral 
Galerkin and modified nonlinear Galerkin methods applied to 
(\ref{8.1})-(\ref{8.3}). The results related to the spectral Galerkin method are an improvement over the Oskolkov work in \cite{AP3} as the constants in error bounds do not depend on time. Zang et al. in \cite{ZD21} have applied a pressure projection method based on the differences of two local Gauss integrations and have derived semidiscrete optimal error estimates for the lowest equal-order finite element space pairs. Then, they have employed a backward Euler method for the time discretization and have discussed the stability and convergence analysis for the fully discrete approximations.
	For more developments of numerical methods applied to the model (\ref{8.1})-(\ref{8.3}) and their finite element analysis, one may refer to  \cite{BN}, \cite{ANS6},  \cite{ZD20} and the literature mentioned therein.  

\noindent	
 As can be seen, the literature is confined to the finite element analysis for the continuous Galerkin methods applied to the problem (\ref{8.1})-(\ref{8.3}). The work dedicated to the finite element analysis for discontinuous Galerkin (DG) methods is still missing from the literature. In recent years, discontinuous Galerkin finite element methods (DGFEM) have been proven to be powerful and popular computational methods for the numerical solution of partial differential equations. The DGFEMs possess many vital properties: They are element-wise conservative and readily parallelizable compared to the continuous Galerkin method. Furthermore, they compensate the continuous Galerkin methods that fail to support nonuniform higher-order local approximation over the mesh for solutions whose smoothness exhibits variation over the computational domain and for unstructured grids. Due to these properties, the DGFEM has received significant recognition from theoretical and practical perspectives.


\noindent
DG methods, introduced in \cite{LR74,RH73}, have been applied to the Euler and  Navier-Stokes equations as early as in \cite{BO1999,LC93}.  Later on, Girault {\it et al.} \cite{GRW005} have presented the first rigorous finite element analysis for the DG method with nonoverlapping domain decomposition and approximations of order $k=1,\,2,\, 3$ for the steady state Stokes and Navier-Stokes system of equations. The authors have derived a uniform discrete inf-sup condition for the discrete discontinuous spaces and established optimal energy error estimates for the velocity and $L^2$-estimates for the pressure for the proposed DG method.  As an extension to the analysis in \cite{GRW005}, the authors in \cite{RG06} have established an improved inf-sup condition and discussed several numerical methods along with their numerical convergence rates. Kaya {\it et al.} in \cite{KR05} have extended the work to the  time-dependent Navier-Stokes equations by applying a  linear subgrid-scale eddy viscosity method combined with discontinuous Galerkin approximations. They have obtained optimal semi-discrete error estimates of the velocity and pressure with reasonable dependence on Reynolds number. They have further applied first and second-order accurate schemes in the time direction and derived fully discrete optimal velocity error estimates. In \cite{GRW05}, Girault {\it et al.} have employed a projection method to decouple the velocity and pressure with a discontinuous Galerkin method for the time-dependent incompressible Navier-Stokes equations and established the optimal error estimates for the velocity and suboptimal for pressure. A few more notable efforts that deal with DG methods for incompressible NSEs can be found in \cite{CW10,CKS09,PE10,JHYW18} and literature therein. And for the references related to the numerical methods and their convergence rates for the DG methods applied to the Navier-Stokes equations, one may refer to  \cite{LFAR19, PMB18, SFE07}.
	
\noindent	
In all the literature mentioned above, for the steady and transient Navier-Stokes equations, the error bounds for the finite element approximations have been shown to be optimal only for the energy norm for the velocity. However, the optimal convergence rate has been numerically achieved for both energy and $\bL^2$-norms. The authors in \cite{GRW005} have provided a hint about the proof of $\bL^2$-norm error estimate of velocity for the SIPG (symmetric interior penalty Galerkin) method to the steady Navier-Stokes model only. But, there is no detailed analysis available in the literature. Furthermore, no rigorous analysis exists for the optimal $\bL^2$-norm error estimate for $t>0$ to the discontinuous velocity approximations for time dependent Navier-Stokes problem. And, to the best of our knowledge, there is hardly any literature dedicated to the finite element analysis of DG methods applied to the Kelvin-Voigt equations of motion. This article can be considered as the first attempt in this direction. The article mainly focuses on deriving semidiscrete and fully discrete optimal error estimates for the symmetric interior penalty discontinuous Galerkin method applied to the problem (\ref{8.1})-(\ref{8.3}) as the  non-symmetric interior penalty discontinuous Galerkin method is known to provide suboptimal error estimates due to its  dependence on the degree of polynomial approximation (\cite{RG06}). The main ingredients in achieving the goals of the article are as follows:

\begin{itemize}
	\item[1.] Based on the analysis of Heywood and  Rannacher (\cite{HR82}) developed for the  non-conforming finite element method applied to the time-dependent Navier-Stokes equations,  we introduce an $L^2$-projection $\bP_h$ onto an appropriate DG finite element space with the help of the approximation operator  $\bR_h$ (see Section \ref{s3}) and aim at deriving the approximation properties for $\bP_h$.
	\item[2.] Next, we define a modified Sobolev-Stokes projection $\bS_h$ (\cite{bnpdy}) for broken Sobolev spaces, which plays an essential role in deriving the semidiscrete error estimates related to the DG method. Then, we concentrate on establishing the estimates for $\bS_h$ in which the approximation properties of  $\bP_h$ play an vital role. Although we apply the ideas of \cite{bnpdy}, there are differences and analytical difficulties due to the DG formulation and difference in finite element spaces. In particular, the analysis of the nonlinear term in DG formulation needs a special kind of attention. 
	\item[3.] With the help of the modified Sobolev-Stokes projection $\bS_h$ and duality arguments, we achieve
	optimal {\it a priori} error estimates for the semidiscrete discontinuous Galerkin approximations to the velocity in $L^{\infty}({\bf L}^2) $-norm
	and  pressure in $L^{\infty}( L^2)$-norm. These estimates are further shown to be uniform under the smallness assumption on the data.
	\item[4.] We then apply the backward Euler scheme to the semidiscrete discontinuous Kelvin-Voigt model and establish optimal error estimates for fully discrete velocity and pressure.
	\item[5.] Finally, we provide numerical examples and analyze the outcomes to verify the theoretical results.
	\end{itemize}
	This article is divided into the following sections: The functional spaces with notations required for the problem analysis, 
the discontinuous weak formulation, and some basic assumptions are presented in Section \ref{s2}.  The semidiscrete discontinuous Galerkin formulation, derivation of a {\it priori} bounds of the discrete solution, and some trace inequalities are dealt with in Section \ref{s3}. The modified Sobolev-Stokes operator and its properties and optimal a {\it priori} error estimates for the velocity are represented in Section \ref{s4}. The optimal a {\it priori} error 
estimates for the pressure are derived in Section \ref{s5}. The backward Euler method for the time discretization is employed, and the fully discrete error estimates are obtained in Section \ref{s6}. A few numerical examples are discussed, and the results are analyzed to verify the theoretical findings in Section \ref{s7}. Finally, the  main contributions of the article are summarized in Section \ref{s8}.
	
	\section{Preliminaries and Weak Formulations}\label{s2}
	\se
In the rest of the paper, we denote by bold face letters the $\R^2$-valued
function spaces such as $\bH_0^1 (\Omega) = (H_0^1(\Omega))^2$, $ \bL^2 (\Omega) = (L^2(\Omega))^2$ etc.
For finite element analysis, we use the standard Sobolev spaces $W^{l,r}(\Omega)$ (\cite{A75}), which is defined for any nonnegative integer $l$ and $r\geq 1$, as follows:
\[W^{l,r}(\Omega)=\{\phi\in L^r(\Omega):\, \forall\, |m|\leq l,\,\, \partial^m \phi\in L^r(\Omega)\}, \]
where $m$ is a two (or a three) dimensional multi-index, $|m|$ is the sum of its components and $\partial^m\phi$ are the partial derivatives of $\phi$ of order $|m|$. 
The norm in $W^{l,r}(\Omega)$ is denoted by $\|\cdot\|_{l,r,\Omega}$ and the seminorm by $|\cdot|_{l,r,\Omega}$. The $L^2$ inner-product is denoted by $(\cdot,\cdot)$.  For the Hilbert space $H^l(\Omega):=W^{l,2}(\Omega)$, the norm is denoted by $\|\cdot\|_{l,\Omega}$ or $\|\cdot\|_l$.
The space $\bH^1_0 (\Omega)$ is equipped with the norm
\[ \|\nabla\bv\|= \big(\sum_{i,j=1}^{2}(\partial_j v_i, \partial_j
v_i)\big)^{1/2}=\big(\sum_{i=1}^{2}(\nabla v_i, \nabla v_i)\big)^{1/2}. \]
Let $H^m (\Omega)/\R$ be the quotient space with norm $\|\phi\|_{H^m (\Omega)/\R}=\displaystyle\inf_{c\in\R}\|\phi+c\|_m.$ For $m=0$, it is denoted by $L^2 (\Omega)/\R$. 
 And for any Banach space $X$, let $L^p(0, T; X)$ denote the space of measurable $X$-valued functions $\bphi$ on  $ (0,T) $ such that
\[ \int_0^T \|\bphi (t)\|^p_X~dt <\infty~~\mbox {if}~1 \le p < \infty, ~~~\mbox{and}~~~
{\esssup_{0<t<T}}\: \|\bphi (t)\|_X <\infty~~\mbox {if}~p=\infty. \]
The dual space of $H^m(\Omega)$ is denoted by $H^{-m}(\Omega)$ with norm defined as
$$\|\phi\|_{-m}:=\sup \bigg\{\dfrac{(\phi,\psi)}{\:\|\psi\|_m}: \psi \in H^m(\Omega), \|\psi\|_m\neq 0 \bigg\}.$$
From future analysis point of view, we further define the following divergence free spaces: 
\[ \bJ_1=\{\bw\in \bH_0^1(\Omega):\,\nabla\cdot \bw=0\}\quad \text{and}\quad \bJ=\{\bw\in \bL^2(\Omega):\,\nabla\cdot \bw=0,\, \bw\cdot \bn|_{\partial\Omega}=0~~\text{holds weakly}\},\] 
where $\bn$ is the outward normal to the boundary $\partial\Omega$ and $\bw\cdot \bn|_{\partial\Omega}=0$ should be understood in the sense of trace in $\bH^{-1/2}(\partial\Omega)$. \\
Through out the paper, we make the following assumptions:\\
\noindent
({\bf A1}). The initial velocity $\bu_0$ and the external force $\brf$ satisfy, for some positive constant $M_0$ and for time $T$ with $0<T<\infty$,

$\bu_0\in \bH^2 (\Omega)\cap \bJ_1,~~\brf,\brf_t\in L^\infty(0,T;\bL^2 (\Omega))$ with $\|\bu_0\|_2\leq M_0$ and $\displaystyle\sup_{0<t<T}\{\|\brf(\cdot,t)\|,\|\brf_t(\cdot,t)\|\}\leq M_0.$ \\
({\bf A2}). For ${\bf g}\in \bL^2 (\Omega)$, let the unique solutions $\bv\in \bJ_1,~q\in L^2 (\Omega)/ \R$ of the steady state Stokes problem
\begin{align*}
	-\Delta \bv+\nabla q={\bf g},\\
	\nabla\cdot \bv=0 ~\text{in}~\Omega,~~\bv|_{\partial\Omega}=0,
\end{align*}   
satisfy 
\[\|\bv\|_2+\|q\|_{H^1 (\Omega)/ \R}\leq C\|{\bf g}\|.\]
Now the weak formulation of (\ref{8.1})-(\ref{8.3}) is as follows: Find a pair $(\bu(t),p(t))\in \bH_0^1 (\Omega)\times L^2 (\Omega)/ \R$, $t>0$, such that
\begin{align}
	(\bu_t,\bphi)+\kappa\,(\nabla\bu_t,\nabla\bphi)+&\nu(\nabla\bu,\nabla\bphi)+(\bu\cdot\nabla\bu,\bphi)-(p,\nabla\cdot\bphi)=(\brf,\bphi)\quad \forall \bphi\in \bH_0^1 (\Omega),\label{w1}\\
	&(\nabla\cdot\bu,q)=0\quad \forall q\in  L^2 (\Omega)/ \R,\label{w2}\\
	&(\bu(0),\bphi)=(\bu_0,\bphi)\quad \forall \bphi\in \bH_0^1 (\Omega).\label{w3}
\end{align}
\noindent
Now, we recall the following Sobolev inequalities which will be useful to estimate the nonlinear terms.
\begin{lemma}[\cite{T95}]\label{sobolev1}
For any open set $\Omega\subset\R^2$ and $\bv\in \bH_0^1(\Omega)$
\begin{align*}
\|\bv\|_{L^4(\Omega)}\leq 2^{1/4}\|\bv\|^{1/2}\|\nabla\bv\|^{1/2}.
\end{align*}
In addition, when $\Omega$ is bounded and $\bv\in \bH^2(\Omega)$, the following estimate holds:
\begin{align*}
\|\bv\|_{L^\infty(\Omega)}\leq C\|\bv\|^{1/2}\|\Delta\bv\|^{1/2}.
\end{align*}
\end{lemma}
\noindent
In Lemma \ref{exactpriori} below, we state the regularity estimates of weak solution pair $(\bu,\,p)$ of (\ref{w1})-(\ref{w3}), which will be used in the subsequent error analysis. 
\begin{lemma} [\cite{bnpdy}]\label{exactpriori}
	Let the assumptions ({\bf A1}) and ({\bf A2}) be satisfied. Then, there exists a positive constant $C=C(\kappa,\nu,\alpha,C_2, M_0)$, such that,  for $t>0$,  the following estimates hold true:
	\begin{align*}
		\sup_{0<t<\infty}\big\{\|\bu(t)\|_2+\|p(t)\|_{H^1 (\Omega)/R}+\|\bu_t(t)\|_2+\|p_t(t)\|_{H^1 (\Omega)/R}\big\} &\leq C, \\
		e^{-2\alpha t}\int_0^t e^{2\alpha s}\big(\|\bu(s)\|_{2}^2+\|\bu_s(s)\|_2^2+\|p(s)\|_{H^1 (\Omega)/R}^2+\|p_s(s)\|_{H^1 (\Omega)/R}^2\big)\,ds &\leq C.
	\end{align*}
\end{lemma}

\noindent
 Next, for defining the DG formulation, we consider a quasiuniform family of triangulations $\mathcal{T}_h$ of $\bar{\Omega}$,
consisting of triangles of maximum diameter $h$. Let $\Gamma_h$ denotes the set of all edges of $\mathcal{T}_h$ and $\bn_e$ represents a unit normal to each edge $e\in \Gamma_h$. 
For an edge	$e$ on the boundary $\partial\Omega$, $\bn_e$ is taken to be the unit outward vector normal to $\partial \Omega$ and for an edge $e$ shared by two elements $T_m$ and $T_n$ of $\mathcal{T}_h$, a unit normal vector $\bn_e$ is directed from $T_m$ to $T_n$.  The jump $[\phi]$ and the average $\{\phi\}$ of a function $\phi$ on an edge $e \in \Gamma_h$ are defined by
\[ [\phi]=(\phi|_{T_m})|_e-(\phi|_{T_n})|_e, \quad \{\phi\}=\frac{1}{2}(\phi|_{T_m})|_e+\frac{1}{2}(\phi|_{T_n})|_e,\] 
and on an edge $e\in \partial\Omega$ coincide with the value of $\phi$ on $e$, in both the cases.\\
We further need the following discontinuous spaces 
\begin{align*}
	\bV &=\{\bw\in  \bL^2(\Omega):\, \bw|_{T}\in \textbf{W}^{2,4/3}(T),\quad \forall T\in \mathcal{T}_h\},\\
	M &= \{q\in L^2(\Omega)/ \R:\, q|_T\in W^{1,4/3}(T),\quad\forall	T\in \mathcal{T}_h\}.
\end{align*}
equipped with the "broken" norms 
\begin{align}\label{nt1}
	&\|\bv\|_\varepsilon =(\vertiii{\nabla \bv}_0^2+J_0(\bv,\bv))^{1/2}\quad \forall \bv\in \bV,\\
	&\|q\|_{L^2/ \R} =\|q\|_{L^2(\Omega)/ \R}\quad \forall q\in M,\nonumber
\end{align}
where, 
\[\vertiii{\cdot}_l=\left(\sum_{T\in \mathcal{T}_h} \|\cdot\|^2_{l,T}\right)^{1/2}, \,\, l\ge 0.\]
The jump term $J_0$  that appeared on the right hand side of (\ref{nt1}) is defined as
\[J_0(\bv,\bw)=\sum_{e\in\Gamma_h}\frac{\sigma_e}{|e|}\int_e[\bv]\cdot [\bw]\,ds,\]
where $|e|$ denotes the measure of the edge $e$ and $\sigma_e>0$ is penalty parameter defined for each edge $e$.\\
 In the course of our analysis, we frequently use the following bounds for the  $L^p$-norm of functions in $\bV$ (\cite{GRW005}):
\begin{equation}
	 \|\bv\|_{L^p(\Omega)}\leq \gamma\|\bv\|_\varepsilon,\quad  \forall \bv\in \bV,~~2\le p<\infty, \label{Lp}
\end{equation}
where $\gamma=C(p)$ is a positive constant that depends on $p$.  \\
Below, we state  some standard trace and inverse inequalities on the discontinuous space $V$.
\begin{lemma}[\cite{GR79}]\label{trace}
	For every triangle $T$ in $\mathcal{T}_h$, the following inequalities hold true
	 \begin{align*}
		\|\bv\|_{L^2(e)} &  \leq C(h_T^{-1/2}\|\bv\|_{L^2(T)}+h_T^{1/2}\|\nabla\bv\|_{L^2(T)})\quad \forall e\in \partial T,\,\, \forall\bv\in \bV,\\
		\|\nabla\bv\|_{L^2(e)} & \leq C(h_T^{-1/2}\|\nabla\bv\|_{L^2(T)}+h_T^{1/2}\|\nabla^2 \bv\|_{L^2(T)})\quad \forall e\in \partial T,\,\,\forall \bv\in \bV,\\
		\|\bv\|_{L^4(e)} & \leq Ch_T^{-3/4}(\|\bv\|_{L^2(T)}+h_T\|\nabla \bv\|_{L^2(T)})\quad \forall e\in \partial T,\,\, \bv\in \bV,
	\end{align*}
where $h_T>0$ is the diameter of triangle $T$.
\end{lemma}

\noindent
In order to write the discontinuous Galerkin formulation of (\ref{8.1})-(\ref{8.3}), we define the bilinear forms $a:\textbf{V}\times\textbf{V}\rightarrow \mathbb{R}$ and $b:\textbf{V}\times M\rightarrow \mathbb{R}$ as follows:
\begin{align}
	a(\bw,\bv) &= \sum_{T\in \mathcal{T}_h}\int_T \nabla \bw:\nabla \bv\,dT -\sum_{e\in \Gamma_h} \int_e\{\nabla \bw\}\bn_e\cdot [\bv]\,ds +\epsilon\sum_{e\in \Gamma_h}\int_e \{\nabla \bv\}\bn_e\cdot [\bw]\,ds,\label{1.6}\\
	b(\bv,q) &= -\sum_{T\in \mathcal{T}_h}\int_T q\nabla \cdot \bv\,dT +\sum_{e\in \Gamma_h }\int_e \{q\}[\bv]\cdot \bn_e \,ds,\label{1.7}
\end{align}
Note that, here $\epsilon=\pm 1$. If $\epsilon=-1$, then the DG formulation is known as SIPG and if $\epsilon=+1$, then it is called as NIPG.  \\
\noindent
This article is dedicated to the derivation of optimal error estimates for the SIPG case. Therefore,  throughout the paper, we make the following assumption:\\
({\bf A3}).  The SIPG DG formulation is considered with $\epsilon=-1$ and the penalty parameter $\sigma_e$ is  bounded below by sufficiently large $\sigma_0>0$. \\
We further define the following trilinear form 
$c(\cdot,\cdot,\cdot)$ for the nonlinear term present in the system (\ref{8.1})-(\ref{8.3}):
\begin{align}\label{4.8}
	c^\bw(\bu,\bz,\brho)=\sum_{T\in\mathcal{T}_h} &\left(\int_T(\bu\cdot \nabla \bz)\cdot \brho\,dT+{\int_{\partial T_-}|\{\bu\}\cdot \bn_T|(\bz^{int}-\bz^{ext})\cdot \brho^{int}}\,ds\right) \nonumber \\
	& +\frac{1}{2}\sum_{T\in \mathcal{T}_h}\int_T(\nabla\cdot \bu)\bz\cdot \brho\,dT -\frac{1}{2}\sum_{e\in \Gamma_h}\int_e[\bu]\cdot \bn_e\{\bz\cdot \brho\}\,ds ,\qquad \forall \bu,\bz,\brho\in \textbf{V},
\end{align} 
where 
\[\partial T_-=\{\textbf{x}\in \partial T:\{\bw\}\cdot \bn_T<0\}.\]
The superscript $\bw$ denotes the dependence of $\partial T_-$ on $\bw$ and the superscript int ( respectively ext) refers to the trace of the function on a side of $T$ coming from the interior  (respectively exterior) of $T$ on that side. The exterior trace is considered to be zero for a side of $T$ belongs to $\partial \Omega$. The first two terms in the definition of $c(\cdot,\cdot,\cdot)$ were introduced in \cite{LR74} for solving transport problems; the last term is taken to ensure the positivity property of $c$ ((\ref{5.5})). Using the definition of $b(\cdot,\cdot)$, the trilinear form can be presented as
\[c^\bw(\bu,\bz,\brho)=\sum_{T\in\mathcal{T}_h}  \left(\int_T(\bu\cdot \nabla \bz)\cdot \brho\,dT+{\int_{\partial T_-}|\{\bu\}\cdot \bn_T|(\bz^{int}-\bz^{ext})\cdot \brho^{int}}\,ds\right)-\frac{1}{2}b(\bu,\bz\cdot \brho). \]
The above definition of trilinear form is motivated from the Lesaint-Raviart upwinding scheme (see \cite{LR74}) and is introduced in \cite{GRW005} for Navier-Stokes system of equations.
It can be easily verified that  for $\bu,\bz,\brho\in  \bH_0^1(\Omega)$, the trilinear form $c$ satisfies
\begin{equation}
	c(\bu;\bz,\brho)=\int_\Omega (\bu\cdot \nabla \bz)\cdot \brho\,dT+\frac{1}{2}\int_\Omega (\nabla \cdot \bu)\bz\cdot \brho\,dT,\,\,.\label{4.9}
\end{equation}
The superscript $\bw$ is dropped in (\ref{4.9}) since the integral on $\partial T_-$ disappears. \\
Using the integration by parts formula, the trilinear form  can be reduced to
\begin{align}
	c^{\bu}(\bu,\bz,\brho)=& -\sum_{T\in \mathcal{T}_h}\left(\int_T (\bu \cdot \nabla \brho)\cdot \bz\,dT+\frac{1}{2}\int_T(\nabla\cdot \bu)\bz\cdot\brho\,dT\right)+\frac{1}{2}\sum_{e\in \Gamma_h}\int_e [\bu]\cdot \bn_e \{\bz\cdot\brho\}\,ds\nonumber\\
	&  -\sum_{T\in \mathcal{T}_h}\int_{\partial T_-}|\{\bu\}\cdot \bn_T|\bz^{ext}\cdot (\brho^{int}-\brho^{ext})\,ds +\int_{\Gamma_+}|\bu\cdot \bn|\bz\cdot \brho\,ds,  \;\;\bu,\bz,\brho \in\textbf{V}, \label{integration}
\end{align}
where $\Gamma_+$ is the subset of $\partial\Omega$ with $\bu\cdot \bn>0$. For a proof, please refer to \cite{GRW005}.
Using a particular choice $\bz=\brho$ in (\ref{integration}), we arrive at
\begin{equation}
	c^\bu(\bu,\bz,\bz)\geq 0, \,\,\, \bu,\bz \in\textbf{V}.\label{5.5}
\end{equation}

\noindent
Now, the DG formulation of (\ref{8.1})-(\ref{8.3}) is as follows:  Find the pair $(\bu (t),p(t))\in \bV \times M,~ t>0$, such that 
\begin{align}
	(\bu _t(t),\bphi)+\kappa\,(a(\bu_t(t),\bphi)&+J_0(\bu_t(t),\bphi))+\nu (a(\bu(t),\bphi)+J_0(\bu(t),\bphi))\nonumber\\
	&+c^\bu(\bu(t),\bu(t),\bphi)+b(\bphi,p(t))= (\textbf{f}(t),\bphi)\quad \forall \bphi\in \bV,\label{8.5}\\
	b(\bu(t),q)&= 0\quad \forall q\in M,\label{8.6}\\
	(\bu(0),\bphi)&= (\bu_0,\bphi)\quad \forall \bphi\in \bV.\label{8.7}
\end{align}
The consistency proof of (\ref{8.5})-(\ref{8.7}) can be done following the similar analysis as adopted in \cite[Lemma 3.2]{KR05} for the DG formulation of NSE. 

	\section{Semidiscrete Discontinuous Galerkin Formulation}\label{s3}
\se

Consider $\bV_h$ and $M_h$ be the finite-dimensional discontinuous subspaces of $\bV$ and $M$, respectively, defined for any positive integer $k\ge 1$ as follows:
\begin{align*}
	\bV_h & = \{\textbf{v}\in  \bL^2(\Omega)\, :\,\, \forall T \in \mathcal{T}_h, \,\, \textbf{v}\in (\mathbb{P}_k(T))^2\},\\
	M_h & = \{q\in  L^2(\Omega)/ \R\, :\, \forall T\in \mathcal{T}_h,\,\, q\in \mathbb{P}_{k-1}(T)\}.
\end{align*}
Let there exists a projection operator $\textbf{R}_h\in \mathcal{L}(\bH^1(\Omega);\textbf{V}_h)$, a linear map from $\bH^1(\Omega)$ to $\textbf{V}_h$, with the properties stated below. 
\begin{lemma}[\cite{GRW005}]\label{Rhp}
	For $\bV_h$, there exists an operator $\textbf{R}_h\in \mathcal{L}(\bH^1(\Omega);\textbf{V}_h)$, such that for any $T\in \mathcal{T}_h$,
	\begin{align}
		\forall \bv\in  \bH^2(\Omega)\cap\bH_0^1(\Omega),~~ &\|\bR_h(\bv)-\bv\|_\varepsilon\leq C h |\bv|_2,\label{2.5}\\
		\forall s \in [1, 2],~ \forall \bv\in \bH^s(\Omega),~~ &\|\bv-\textbf{R}_h(\bv)\|_{L^2(T)}\leq C h_T^s|\bv|_{s,\Delta_T}, \label{2.6}
	\end{align}
	where $\Delta_T$ is a suitable macro-element containing $T$.
\end{lemma}
\noindent
Further, let there exists an operator $r_h\in \mathcal{L}(L^2(\Omega)/ \R;M_h)$ (see \cite{GRW005}), such that, for any $T\in \mathcal{T}_h$, the following properties are satisfied:
\begin{equation}
	\forall q\in L^2(\Omega)/ \R,~~\forall z_h\in \mathbb{P}_{k-1}(T),\quad  \int_T  z_h(r_h(q)-q)\,dT=0,\label{2}
\end{equation}
\begin{equation}
	\forall s\in[0,k],	\forall q\in H^s(\Omega)\cap L^2(\Omega)/ \R,\quad \|q-r_h(q)\|_{L^2(T)}\leq C h_T^s|q|_{s,T}.\label{2.1}
\end{equation}
Moreover, for the operator $\bR_h$, we have (see \cite{GRW005})
\begin{equation}
	\forall \bv\in \bH_0^1(\Omega), \quad \forall q \in M_h,\quad b(\textbf{R}_h(\bv)-\bv,q)=0, \label{3.9}
\end{equation}
and the stability property (see \cite{GRW05}): there exists a constant $C$, independent of $h$, such that 
\begin{equation}
	\forall \bv\in \bH_0^1(\Omega),\,\,	\|\bR_h \bv\|_\varepsilon\leq C  |\bv|_1.\label{Rh}
\end{equation}
We define below the semidiscrete discontinuous Galerkin approximations for the equations (\ref{w1})-(\ref{w3}) as follows:
For $t>0$, find $(\bu_h (t),p_h(t))\in \bV_h \times M_h$ such that 
\begin{align}
	(\bu_{ht}(t),\bphi_h)&+\kappa\, (a(\bu_{ht}(t),\bphi_h)+J_0(\bu_{ht}(t),\bphi_h))+\nu (a(\bu_h(t),\bphi_h)+J_0(\bu_h(t),\bphi_h))\nonumber\\
	&+c^{\bu_h}(\bu_h(t),\bu_h(t),\bphi_h)+b(\bphi_h,p_h(t))= (\textbf{f}(t),\bphi_h), \,\,\,\forall \bphi_h\in\bV_h \label{8.8}\\
	& b(\bu_h(t),q_h)= 0, \,\,\,\forall q_h\in M_h\label{8.8n}\\
	&(\bu_h(0),\bphi_h)= (\bu_0,\bphi_h) \,\,\,\forall \bphi_h\in\bV_h.\label{8.9}
\end{align}
The discrete discontinuous space $\bJ_h$, which is analogous to $ \bJ_1$, is defined as
\[\textbf{J}_h=\{\bv_h\in \textbf{V}_h:\forall q_h\in M_h,\quad b(\bv_h,q_h)=0\}.
\]
The equivalent formulation of (\ref{8.8})--(\ref{8.9}) on $\bJ_h$  is as follows:
find $\bu_h(t)\in \bJ_h$, such that for $t>0$
\begin{align}
	(\bu_{ht},\bphi_h)+\kappa\, (a(\bu_{ht}(t),&\bphi_h)+J_0(\bu_{ht}(t),\bphi_h))+\nu(a(\bu_h,\bphi_h)\nonumber\\
	&+J_0(\bu_h,\bphi_h))
	+c^{\bu_h}(\bu_h,\bu_h,\bphi_h)=({\bf f}, \bphi_h),~\forall \bphi_h\in \bJ_h. \label{ne1}
\end{align}
The next two lemmas state the boundedness properties and the coercivity of the bilinear form $(a+J_0)(\cdot,\cdot)$, respectively.
\noindent
\begin{lemma}[\cite{R08}]\label{cont}
	There exists a constant $C_1>0$, independent of $h$, such that, for all $\bv_h,\bw_h\in \bV_h$,
	\begin{equation*}
		|a(\bv_h,\bw_h) +J_0(\bv_h,\bw_h )|\leq C_1 \|\bv_h\|_\varepsilon\|\bw_h\|_\varepsilon.
	\end{equation*}
\end{lemma}
\begin{lemma}[\cite{W78}]\label{coer}
	Under the assumption of ({\bf A3}) there exists a constant $C_2>0$, independent of $h$, such that 
	\begin{equation*}
		\forall \bv_h\in \textbf{V}_h,\quad a(\bv_h,\bv_h)+J_0(\bv_h,\bv_h)\geq C_2\|\bv_h\|_\varepsilon ^2.
	\end{equation*} 
\end{lemma}
\noindent
Similar to the continuous case, we have a discrete inf-sup condition for the pair of discontinuous spaces ($\widetilde{\bV}_h,M_h$), where
\[ \widetilde{\textbf{V}}_h=\big\{\bv_h\in \textbf{V}_h: \forall e \in \Gamma_h, \quad \int_e \bq_h\cdot [\bv_h]\,ds=0,\quad \forall \bq_h\in (\mathbb{P}_{k-1}(e))^2\big\}, \]
which will be an essential tool for deriving the pressure estimates.
\begin{lemma}[\cite{GRW005}]\label{inf-sup}
	There exists a constant $\beta^*>0$, independent of $h$, such that 
	\begin{equation*}
		\inf_{p_h\in M_h}\sup_{\bv_h\in \widetilde{\textbf{V}}_h}\frac{b(\bv_h,p_h)}{\|\bv_h\|_\varepsilon \|p_h\|_0}\geq \beta ^*.
	\end{equation*}
\end{lemma}
\noindent
Analogous to the first estimate in Lemma \ref{sobolev1}, we have the following estimate for the elements in the discrete space $\bV_h$. 
\begin{lemma}[\cite{GRW05}]\label{sobolev2}
When $\Omega$ is convex, there exists a positive constant $C$, independent of $h$, such that
\begin{align*}
	\|\bv_h\|_{L^4(\Omega)}\leq C\|\bv_h\|^{1/2}\|\bv_h\|_\varepsilon^{1/2}, \quad \forall \bv_h\in \bV_h.
\end{align*}
\end{lemma}
\noindent
The following estimates of the trilinear form $c(\cdot,\cdot,\cdot)$ will be useful for our error analysis. 
\begin{lemma}[\cite{GRW05}]
	(i) Assume that $\bu\in  \bW^{1,4}(\Omega)$. There exists a positive constant $C$ independent of $h$ such that 
	\begin{align}
		|c(\bv_h,\bu,\bw_h)|\leq C\|\bv_h\| |\bu|_{\bW^{1,4}(\Omega)}\|\bw_h\|_\varepsilon, \quad \forall \bv_h,\bw_h\in\bJ_h.\label{tri1}
	\end{align}
	(ii) For any $\bv \in \bV$, $\bv_h$, $\bw_h$ and $\bz_h$ in $\bV_h$, we have the following estimate:
	\begin{align}
		|c^{\bv}(\bv_h,\bw_h,\bz_h)|\leq C\|\bv_h\|_\varepsilon\|\bw_h\|_\varepsilon\|\bz_h\|_\varepsilon.\label{tri2}
	\end{align} 
\end{lemma}

\noindent
In Lemma \ref{distrace}, we present a couple of trace inequalities along with inverse inequalities for the discrete discontinuous space $\bV_h$, similar to those stated in Lemma \ref{trace} for the space $\bV$.
\begin{lemma}[\cite{GR79}]\label{distrace}
	For every element $T$ in $\mathcal{T}_h$, the following inequalities hold
	 \begin{align}
		\|\bv_h\|_{L^2(e)}&\leq Ch_T^{-1/2}\|\bv_h\|_{L^2(T)}\quad \forall e\in \partial T, \,\,\forall \bv_h\in\bV_h,\label{T4}\\
		\|\nabla\bv_h\|_{L^2(e)}&\leq Ch_T^{-1/2}\|\nabla\bv_h\|_{L^2(T)}\quad \forall e\in \partial T, \,\,\forall \bv_h\in\bV_h,\label{T5}\\
		\|\nabla\bv_h\|_{L^2(T)}&\leq Ch_T^{-1}\|\bv_h\|_{L^2(T)}\quad \forall \bv_h\in\bV_h,\label{T6}\\
		\|\bv_h\|_{L^4(T)}&\leq Ch_T^{-1/2}\|\bv_h\|_{L^2(T)}\quad \forall \bv_h\in\bV_h,\label{T7}
	\end{align}
	where $C$ is a constant independent of $h_T$ and $\bv_h$. 
\end{lemma}

\noindent
Next, we state the regularity bounds for $\bu_h$ which will be used in deriving the existence and uniqueness of the discrete solution and fully discrete error estimates. 
\begin{lemma}\label{uhpriori}
	Let the assumptions ({\bf A1})-({\bf A3}) be satisfied and let $0<\al<\frac{\nu C_2}{2(\gamma+\kappa C_1)}$. Then, there exists a positive constant $C=C(\kappa,\nu,\alpha,C_2, M_0)$, such that, for each $t>0$, the semidiscrete discontinuous Galerkin solution $\bu_h(t)$, satisfies the following estimates :
	\begin{align*}
		\sup_{0<t<\infty}\|\bu_h(t)\|+e^{-2\al t}\int_0^t e^{2\al s}\left(\|\bu_h(s)\|_\varepsilon^2+\|\bu_{ht}(s)\|_\varepsilon^2\right)\,ds \leq C,\\
		\sup_{0<t<\infty}\left(\|\bu_{htt}(t)\|_{-1,h}+ \|\bu_{htt}(t)\|_\varepsilon\right))+	e^{-2\al t}\int_0^t e^{2\al s}\left(\|\bu_{htt}(s)\|_{-1,h}^2+ \kappa\,\|\bu_{htt}(s)\|_\varepsilon^2\right)\,ds \leq C,
	\end{align*}
	where 
	\begin{align} \|\bu_{htt}\|_{-1,h}=\sup\bigg\{\frac{\langle\bu_{htt},\bphi_h\rangle}{\|\bphi_h\|_\varepsilon},~\bphi_h\in \bV_h,\bphi_h\neq 0\bigg\}. \label{r4}
	\end{align}
	Moreover,
	\begin{align*}
		\limsup_{t\rightarrow\infty}\|\bu_h(t)\|_\varepsilon\leq \frac{ \gamma\|\brf\|_{L^\infty(\bL^2 (\Omega))}}{C_2\nu}.
	\end{align*}
\end{lemma}
\begin{proof}
	Choose $\bphi_h=\bu_h$ in (\ref{ne1}) and apply the coercivity result, positivity of $c(\cdot,\cdot,\cdot)$ (\ref{5.5}), estimate (\ref{Lp}), Cauchy-Schwarz's inequality and Young's inequality to obtain
	\begin{align}\label{pri1}
		\frac{1}{2}\frac{d}{dt}\left(\|\bu_h\|^2+ \kappa(a(\bu_h,\bu_h)+J_0(\bu_h,\bu_h))\right)+\nu C_2\|\bu_h\|_\varepsilon^2 \leq \|\brf\|\|\bu_h\|\leq \frac{\nu C_2}{2}\|\bu_h\|_\varepsilon^2+  \frac{\gamma^2}{2\nu C_2} \|\brf\|^2.
	\end{align}
	A multiplication of (\ref{pri1}) by $e^{2\al t}$, an integration from $0$ to $t$, and a use of estimate (\ref{Lp}), Lemmas \ref{cont} and \ref{coer}, lead to
	\begin{align}\label{pri2}
		e^{2\al t}(\|\bu_h(t)\|^2+C_2\kappa\|\bu_h(t)\|_{\varepsilon}^2)+ (\nu C_2-2\al(\gamma+\kappa C_1) )&\int_0^te^{2\al s}\|\bu_h(s)\|_\varepsilon^2\,ds \leq \|\bu_h(0)\|^2\nonumber\\&+ C_1\kappa \|\bu_h(0)\|_\varepsilon^2 + C\int_0^te^{2\al s}\|\brf(s)\|^2\,ds.
	\end{align}
	Again, multiply (\ref{pri2}) by $e^{-2\al t}$, use the fact that
	\[e^{-2\al t}\int_0^t e^{2\al s}\,ds=\frac{1}{2\al}(1-e^{-2\al t}) \]
	and choose  $0<\al<\frac{\nu C_2}{2(\gamma+\kappa C_1)}$ to obtain
	\begin{align}\label{pr11}
		\|\bu_h(t)\|^2+\|\bu_h(t)\|_{\varepsilon}^2+ e^{-2\alpha t}&\int_0^te^{2\al s}\|\bu_h(s)\|_\varepsilon^2\,ds \leq C.
	\end{align}
	 Again, multiply (\ref{pri1}) by $e^{2\al t}$, integrate from $0$ to $t$,  and  a use of Lemma \ref{cont} implies
	\begin{align*}
		e^{2\al t}(\|\bu_h(t)\|^2+ \kappa(a(\bu_h(t),\bu_h(t))+J_0(\bu_h(t),\bu_h(t)))+\nu C_2\int_0^t e^{2\al s}\|\bu_h(s)\|_\varepsilon^2\,ds \leq \left(\|\bu_h(0)\|^2+\kappa\, C_1\|\bu_h(0)\|_{\varepsilon}^2\right)\\
		+2\al \int_0^t e^{2\al s}(\|\bu_h(s)\|^2+ \kappa (a(\bu_h(s),\bu_h(s))+J_0(\bu_h(s),\bu_h(s))))\,ds+(e^{2\al t}-1) \frac{\gamma^2\|\brf\|^2_{L^\infty(\bL^2 (\Omega))}}{2\al\nu C_2}.
	\end{align*}
	Multiply the above inequality by $e^{-2\al t}$, take limit supremum as $t\rightarrow\infty$ and noting that,
	\[\nu C_2\limsup_{t\rightarrow\infty}e^{-2\al t}\int_0^t e^{2\al s}\|\bu_h(s)\|_\varepsilon^2\,ds=\frac{\nu C_2}{2\al }\limsup_{t\rightarrow\infty}\|\bu_h(t)\|_\varepsilon^2,\]
	we arrive at
	\begin{align}\label{r5}
		\frac{\nu C_2}{2\al }\limsup_{t\rightarrow\infty}\|\bu_h(t)\|_\varepsilon^2\leq  \frac{\gamma^2\|\brf\|^2_{L^\infty(\bL^2 (\Omega))}}{2\al\nu C_2}.
	\end{align}
	Next, differentiating (\ref{ne1}) with respect to $t$, we obtain
	\begin{align}
		(\bu_{htt},\bphi_h)&+\kappa\,(a(\bu_{htt},\bphi_h)+J_0(\bu_{htt},\bphi_h))+\nu(a(\bu_{ht},\bphi_h)+J_0(\bu_{ht},\bphi_h))\nonumber\\
		&+c^{\bu_h}(\bu_{ht},\bu_h,\bphi_h)+c^{\bu_h}(\bu_{h},\bu_{ht},\bphi_h)=(\brf_t, \bphi_h),~\forall\, \bphi_h\in \bJ_h. \label{pri4}
	\end{align}
	Substitute $\bphi_h=\bu_{ht}$ in (\ref{pri4}), apply Lemma \ref{coer} and the fact that  $c^{\bu_h}(\bu_{h},\bu_{ht},\bu_{ht})\geq 0$ from (\ref{5.5}), we obtain
	\begin{align}\label{pri5}
		\frac{d}{dt}\left(\|\bu_{ht}\|^2+ \kappa\,(a(\bu_{ht},\bu_{ht})+J_0(\bu_{ht},\bu_{ht}))\right)+\nu C_2\|\bu_{ht}\|_\varepsilon^2 \leq -2c^{\bu_h}(\bu_{ht},\bu_h,\bu_{ht})+ C \|\brf_t\|^2.
	\end{align}
	On expanding the nonlinear term in (\ref{pri5}), we find that
	\begin{align*}
		c^{\bu_h}(\bu_{ht},\bu_h,\bu_{ht})=&\sum_{T\in\mathcal{T}_h}  \int_T(\bu_{ht}\cdot \nabla \bu_h)\cdot\bu_{ht}\,dT+\sum_{T\in\mathcal{T}_h}\int_{\partial T_-}|\{\bu_{ht}\}\cdot \bn_T|(\bu_h^{int}-\bu_h^{ext})\cdot \bu_{ht}^{int}\,ds\nonumber\\
		& +\frac{1}{2}\sum_{T\in \mathcal{T}_h}\int_T(\nabla\cdot \bu_{ht})\bu_h\cdot \bu_{ht}\,dT -\frac{1}{2}\sum_{e\in \Gamma_h}\int_e[\bu_{ht}]\cdot \bn_e\{\bu_h\cdot \bu_{ht}\}\,ds\\
		=& A_1+A_2+A_3+A_4.
	\end{align*}
	Using H{\"o}lder's inequality, Young's inequality and Lemma \ref{sobolev2}, we can bound $A_1$ as follows:
	\begin{align}
		|A_1|\leq \|\bu_{ht}\|_{L^4(\Omega)} \|\nabla\bu_{h}\|\|\bu_{ht}\|_{L^4(\Omega)}\leq C\|\bu_{ht}\|\|\bu_{ht}\|_\varepsilon\|\bu_{h}\|_\varepsilon\leq \frac{\nu C_2}{8}\|\bu_{ht}\|_\varepsilon^2+C\|\bu_{ht}\|^2\|\bu_{h}\|_\varepsilon^2.\label{pr1*}
	\end{align}
	An application of H{\"o}lder's, Young's, trace inequalities and Lemma \ref{sobolev2} yield
	\begin{align}
		|A_2|\leq & \sum_{T\in\mathcal{T}_h} \|\bu_{ht}\|_{L^4(\partial T)}\|[\bu_{h}]\|_{L^2(\partial T)}\|\bu_{ht}\|_{L^4(\partial T)}\leq C\sum_{T\in\mathcal{T}_h} \|\bu_{ht}\|_{L^4(T)}\frac{1}{|e|^{1/2}}\|[\bu_{h}]\|_{L^2(e)}\|\bu_{ht}\|_{L^4(T)}\nonumber\\
		\leq & C\|\bu_{ht}\|_{L^4(\Omega)}\|\bu_{h}\|_\varepsilon\|\bu_{ht}\|_{L^4(\Omega)}\leq \frac{\nu C_2}{8}\|\bu_{ht}\|_\varepsilon^2+C\|\bu_{ht}\|^2\|\bu_{h}\|_\varepsilon^2.\label{pr2*}
	\end{align}
	The term  $A_3$ can be bounded using H{\"o}lder's, Young's inequalities and Lemma \ref{sobolev2} as
	\begin{align}
		|A_3|\leq &\frac{1}{2}\|\nabla\bu_{ht}\| \|\bu_{h}\|_{L^4(\Omega)}\|\bu_{ht}\|_{L^4(\Omega)}\leq C\|\bu_{ht}\|_\varepsilon\|\bu_{h}\|^{1/2}\|\bu_{h}\|_\varepsilon^{1/2}\|\bu_{ht}\|^{1/2}\|\bu_{ht}\|_\varepsilon^{1/2}\nonumber\\
		\leq & \frac{\nu C_2}{8}\|\bu_{ht}\|_\varepsilon^2+C\|\bu_{ht}\|^2 \|\bu_{h}\|^2\|\bu_{h}\|_\varepsilon^2.\label{pr3*}
	\end{align}
	Similar to $A_2$ and $A_3$, we can bound $A_4$ as follows
	\begin{align}
		|A_4|\leq \frac{\nu C_2}{8}\|\bu_{ht}\|_\varepsilon^2+C\|\bu_{ht}\|^2\|\bu_{h}\|^2\|\bu_{h}\|_\varepsilon^2.\label{pr4*}
	\end{align}
	Collecting the bounds form (\ref{pr1*})-(\ref{pr4*}) and applying in (\ref{pri5}), we arrive at
	\begin{align}\label{pri6}
		\frac{d}{dt}\left(\|\bu_{ht}\|^2+ \kappa\,(a(\bu_{ht},\bu_{ht})+J_0(\bu_{ht},\bu_{ht}))\right)+\nu C_2\|\bu_{ht}\|_\varepsilon^2 \leq C\|\bu_{ht}\|^2\|\bu_{h}\|_\varepsilon^2 (1+\|\bu_{h}\|^2) + C \|\brf_t\|^2.
	\end{align}
	On multiplying (\ref{pri6}) by $e^{2\al t}$, integrating from $0$ to $t$, and using Lemmas \ref{cont} and \ref{coer}, we observe that
	\begin{align*}
		e^{2\al t}(\|\bu_{ht}(t)\|^2&+\kappa C_2 \|\bu_{ht}(t)\|_\varepsilon^2)+ (\nu C_2-2\al(\gamma+\kappa C_1) )\int_0^t e^{2\al s}\|\bu_{ht}(s)\|_\varepsilon^2\,ds \leq (\|\bu_{ht}(0)\|^2+ C_1\kappa \|\bu_{ht}(0)\|_\varepsilon^2)\\
		&+C\int_0^t e^{2\al s}\|\bu_{ht}(s)\|^2\|\bu_{h}(s)\|_\varepsilon^2 (1+\|\bu_{h}(s)\|^2)\,ds
		+ C\int_0^t e^{2\al s} \|\brf_t(s)\|^2\,ds.
	\end{align*}
	Choosing $ 0<\al<\frac{\nu C_2}{2(\gamma+\kappa C_1)}$, applying the assumption ({\bf A1}), (\ref{pr11}), Gronwall's lemma and after a final multiplication by $e^{-2\al t}$, we obtain the estimate as
	\begin{align}
		\|\bu_{ht}(t)\|^2&+\|\bu_{ht}(t)\|_\varepsilon^2+e^{-2\al t}\int_0^t e^{2\al s}\|\bu_{ht}(s)\|_\varepsilon^2\,ds \leq C.\label{pr12}
	\end{align}
	Now, we substitute $\bphi_h=\bu_{htt}$ in (\ref{pri4}) and obtain
	\begin{align}
		\|\bu_{htt}\|^2&+C_2\kappa\,\|\bu_{htt}\|^2_\varepsilon =-\nu(a(\bu_{ht}, \bu_{htt})+J_0(\bu_{ht},\bu_{htt}))\nonumber\\
		&-c^{\bu_h}(\bu_{ht},\bu_h,\bu_{htt})
		-c^{\bu_h}(\bu_{h},\bu_{ht},\bu_{htt})-(\brf_t, \bu_{htt}). \label{pri5*}
	\end{align}
	Drop the first term on the left hand side of (\ref{pri5*}) as it is positive. Apply (\ref{tri2}), Lemma \ref{cont} and Young's inequality to obtain
	\begin{align}
		C_2\kappa\,\|\bu_{htt}\|^2_\varepsilon\leq  C(\|\bu_{ht}\|^2_\varepsilon+ \|\bu_{ht}\|^2_\varepsilon\|\bu_{h}\|^2_\varepsilon+\|\brf_t\|^2).\label{r2}
	\end{align}
	A use of (\ref{pr11}), (\ref{pr12}) in (\ref{r2}) yield
	\begin{align}
		\|\bu_{htt}\|_\varepsilon\leq C.\label{pr13}
	\end{align}
	Finally, we write (\ref{pri4}) as follows:
	\begin{align*}
		(\bu_{htt},\bphi_h)&=-\kappa(a(\bu_{htt},\bphi_h)+J_0(\bu_{htt},\bphi_h))-\nu(a(\bu_{ht},\bphi_h)+J_0(\bu_{ht},\bphi_h))\\
		&-c^{\bu_h}(\bu_{ht},\bu_h,\bphi_h)-c^{\bu_h}(\bu_{h},\bu_{ht},\bphi_h)-(\brf_t, \bphi_h),~\forall\quad \bphi_h\in \bJ_h. 
	\end{align*}
	An application of Lemma \ref{cont} and (\ref{tri2}) leads to
	\begin{align}
		(\bu_{htt},\bphi_h)\leq C(\kappa,\nu)(\|\bu_{htt}\|_\varepsilon+ \|\bu_{ht}\|_\varepsilon+ \|\bu_{ht}\|_\varepsilon\|\bu_{h}\|_\varepsilon+\|\brf_t\|)\|\bphi_h\|_\varepsilon. \label{r6}
	\end{align}
	Using (\ref{r4}), the estimates in (\ref{pr12}), (\ref{pr13}) and the assumption ({\bf A1}) in (\ref{r6}), we obtain
	\begin{align}\label{pri7}
		\|\bu_{htt}\|_{-1,h}\leq C(\kappa,\nu) (\|\bu_{htt}\|_\varepsilon+\|\bu_{ht}\|_\varepsilon+\|\brf_t\|)\leq C.
	\end{align}
	After taking square of (\ref{pr13}) and (\ref{pri7}), add the resulting two inequalities,  multiply by $e^{2\alpha t }$ and integrate with respect to time from $0$ to $t$. Then multipliply by $e^{-2\alpha t}$ to arrive at
	\begin{align}\label{pri8}
		e^{-2\alpha t}\displaystyle{\int_0^t}e^{2\alpha s}( \|\bu_{htt}(s)\|^2_{-1,h}+	\kappa\,\|\bu_{htt}(s)\|^2_\varepsilon) ds\leq C.
	\end{align}
	A combination of (\ref{pr11}),  (\ref{r5}), (\ref{pr12}), (\ref{pr13}), (\ref{pri7}) and (\ref{pri8}) completes the proof of Lemma \ref{uhpriori}.
\end{proof}

\noindent
Now, the existence and uniqueness of the semidiscrete discontinuous Galerkin Kelvin-Voigt model (\ref{8.8})-(\ref{8.9}) (or (\ref{ne1})) can be proved following the analysis in \cite{KR05} and using the results in (\ref{5.5}), Lemmas \ref{coer},  \ref{inf-sup}, \ref{uhpriori}.  \\
\noindent
For deriving the optimal error estimates for semidiscrete discontinuous velocity and pressure approximations, we work on the weakly divergence free spaces. Therefore, there is a need to define an approximation operator on  $\bJ_1\cap \bH^2 (\Omega)$.
\begin{lemma}\label{ih}
	For every $\bv\in  \bJ_1\cap \bH^2 (\Omega)$, there exists an approximation $i_h \bv\in \bJ_h$, such that, the following approximation property hold true
	\begin{equation*}
		\|\bv-i_h \bv\|_\varepsilon\leq Ch |\bv|_2.
	\end{equation*}
	\begin{proof}
		From Lemma \ref{Rhp}, there exists an operator $\bR_h:\bH_0^1 (\Omega)\rightarrow \bV_h$ satisfying $\|\bR_h \bv\|_\varepsilon\leq C|\bv|_1$ (see (\ref{Rh})) and 
		\[b(\bv-\bR_h \bv,q_h)=0, \quad \forall ~~ \bv\in \bH_0^1(\Omega), ~ q_h\in M_h.\]
		We observe that 
		$b(\bR_h \bv,q_h)=0$. This implies that
		$\bR_h\bv\in \bJ_h$. 
		We now restrict $\bR_h$ on $\bJ_1\cap \bH^2 (\Omega)$ and call it as $i_h$. Therefore, $i_h$ is an approximation operator from $\bJ_1\cap \bH^2 (\Omega)$ to $ \bJ_h$.  Using (\ref{2.5}), we obtain
		\[\|\bv-i_h\bv\|_\varepsilon=\|\bv-\bR_h\bv\|_\varepsilon\leq Ch|\bv|_2.\]
		This completes the rest of the proof.
	\end{proof}
\end{lemma}
\noindent We also define a projection $\bP_h:\bL^2 (\Omega)\rightarrow \bJ_h$, such that, for every $\bv\in\bL^2 (\Omega)$, 
$(\bv-\bP_h\bv,\bv_h)=0,\quad \forall \bv_h\in \bJ_h$. \\
The projection $\bP_h$ satisfies the following properties
\begin{equation}\label{P_h}
	\|\bv-\bP_h\bv\|+h\|\bv-\bP_h\bv\|_\varepsilon\leq Ch^2|\bv|_2, \quad \forall\,\, \bv\in   \bJ_1\cap\bH^2 (\Omega).
\end{equation}
The energy norm estimates in (\ref{P_h}) can be easily obtained using the definition of $\bP_h$ and the estimates in Lemma \ref{ih} and the $L^{\infty}(\bL^2)$-norm estimate can be derived using elliptic duality argument and energy norm estimate of  $\bP_h$.\\
In Theorem \ref{semierror}, we provide the main contribution of the article related to the semidiscrete error estimates.
\begin{theorem}\label{semierror}
	Let the assumptions ({\bf A1})-({\bf A3}) be satisfied and let  $0<\al<\min \left\{ \frac{\nu C_2}{2(\gamma+\kappa C_1)},\frac{\nu}{2\kappa}\right\}$. Further, let the discrete initial velocity $ \bu_{h}(0)\in \bJ_h$ with $\bu_{h}(0) =
	\bP_h\bu_0$. Then, there exists a positive constant $C=C(\kappa,\nu,\alpha,C_2, M_0)$, independent of $h$, such that
	\begin{align*}
		\|(\bu-\bu_h)(t)\|+ h\|(\bu-\bu_h)(t)\|_\varepsilon\leq Ce^{Ct}h^2,
	\end{align*}
	where $C$ is a positive constant independent of $h$.
\end{theorem}
\begin{theorem}\label{semierrorp}
	Let the assumptions ({\bf A1})-({\bf A3}) be satisfied and let  $0<\al<\min \left\{ \frac{\nu C_2}{2(\gamma+\kappa C_1)},\frac{\nu}{2\kappa}\right\}$. Further, let the discrete initial velocity $ \bu_{h}(0)\in \bJ_h$ with $\bu_{h}(0) =
	\bP_h\bu_0$. Then, there exists a positive constant $C=C(\kappa,\nu,\alpha,C_2, M_0)$, independent of $h$, such that
	\begin{align*}
		\|(p-p_h)(t)\|\leq C\,e^{Ct}h,
	\end{align*}
	where $C$ is a positive constant independent of $h$.
\end{theorem}

	\section{Error estimates for velocity}\label{s4}
	\se
This section deals with the  optimal estimates of the velocity error $\de=\bu-\bu_h$ in $\bL^2$ and energy-norms for $t>0$. We split the error into two parts, $\de=\bxi+\be$, where $\bxi=\bu-\bv_h$ represents the error inherent  due to the DG finite element approximation of (\ref{8.1}) by a linearized Kelvin-Voigt problem, and $\be=\bv_h-\bu_h$ represents the error caused by the presence of the nonlinearity in problem (\ref{8.1}). The linearized equation to be satisfied by the auxiliary function $\bv_h$ is:
	\begin{align}
		(\bv_{ht}, {\bphi}_h)+\kappa\,(a(\bv_{ht},{\bphi}_h)
		&+J_0(\bv_{ht},\bphi_h))+\nu (a (\bv_h, {\bphi}_h)+J_0(\bv_{h},\bphi_h))\nonumber\\
			=&({\bf f},{\bphi}_h)-c^{\bu}(\bu,\bu,{\bphi}_h)\quad \forall \bphi_h\in  \bJ_h.\label{9.1}
	\end{align}  
	Below we derive some estimates of $\bxi$. Subtracting (\ref{9.1}) from (\ref{8.5}), the equation in $\bxi$ can be written as 
	\begin{align}
		(\bxi_t,{\bphi}_h)+\kappa\,(a({\bxi}_t,{\bphi}_h)+J_0(\bxi_{t},\bphi_h))+\nu (a(\bxi,{\bphi}_h)+J_0(\bxi,\bphi_h))=-b(\bphi_h,p), \quad \bphi_h\in \bJ_h.\label{9.2}
	\end{align}
For deriving the optimal error estimates of $\bxi$ in $\bL^2$ and energy-norms  for $t>0$, we introduce, as in \cite{bnpdy}, the following modified Sobolev-Stokes's  projection $\bS_h\bu:[0,\infty)\rightarrow  \bJ_h$ satisfying
\begin{align}\label{E616}
	\kappa  (a(\bu_t-\bS_h \bu_t,{\bphi}_h)+J_0 (\bu_t&-\bS_h \bu_t,{\bphi}_h))+ \nu (a(\bu-\bS_h \bu,{\bphi}_h)\nonumber\\
	&+J_0 (\bu-\bS_h \bu,{\bphi}_h))=-b(\bphi_h,p)\; \; \;
	\forall {\bphi}_h\in \bJ_h,
\end{align}
where $\bS_h\bu(0)= \bP_h\bu_0.$ In other words, given $(\bu,p),$ find $ \bS_h\bu:[0,\infty)\rightarrow  \bJ_h$ satisfying (\ref{E616}). With $ \bS_h\bu$ defined as  above, we now split $\bxi$ as 
\[\bxi=\bu-\bS_h\bu+\bS_h\bu-\bv_h=:\bzeta+\brho.\] Using (\ref{E616}),  we find the equaion in $\bzeta$ to be
\begin{align}\label{E61}
	\kappa  (a(\bzeta_t,{\bphi}_h)+J_0 (\bzeta_t,{\bphi}_h))+ \nu (a( \bzeta,{\bphi}_h)+J_0 (\bzeta,{\bphi}_h))= -b(\bphi_h,p)\; \; \;
	\forall {\bphi}_h\in \bJ_h,
	\end{align}
Firstly, we will focus on deriving the estimates of $\bzeta$. Next, we will establish the estimtes of $\brho$.  Combinigng these estimates will result in the estimates of $\bxi$. 
	\begin{lemma}\label{zetaerr}
		Let the asssumptions $({\bf A1})$-$({\bf A3})$ hold true and let  $0<\al<\min \left\{ \frac{\nu C_2}{2(\gamma+\kappa C_1)},\frac{\nu}{2\kappa}\right\}$. Then, there exists a positive constant $C=C(\kappa,\nu,\alpha,C_2, M_0)$, such that,  for $t>0$, $\bzeta$ satisfies the following estimates:
		\begin{align}
			\|\bzeta(t)\|^2+h^2\|\bzeta(t)\|^2_\varepsilon+e^{-2\alpha t}	\displaystyle{\int_0^t}e^{2\alpha s}\left(\|\bzeta(s)\|^2+h^2\|\bzeta(s)\|^2_\varepsilon+h^2\|\bzeta_t(s)\|^2_\varepsilon\right)ds\leq &  Ch^{4},\label{9.4}\\
			e^{-2\alpha t}\displaystyle{\int_0^t}e^{2\alpha s}\|\bzeta_t(s)\|^2ds\leq&  Ch^{4}.\label{9.5}
		\end{align}
	\end{lemma}
	\begin{proof}
		Set $\bphi_h=\bP_h \bzeta=\bzeta-(\bu-\bP_h \bu)$ in (\ref{E61}), use the definition of space $\bJ_h$ and obtain 
		\begin{align}
			\kappa\,  (a&(\bP_h\bzeta_t,\bP_h \bzeta)+J_0 (\bP_h\bzeta_t,\bP_h \bzeta))+	\nu(a( \bP_h\bzeta, \bP_h \bzeta)+J_0(\bP_h \bzeta,\bP_h \bzeta))
			=	-\kappa  (a(\bu_t-\bP_h\bu_t,{\bP_h \bzeta})\nonumber\\&+J_0 (\bu_t-\bP_h\bu_t,\bP_h \bzeta))
			-\nu(a(\bu-\bP_h\bu,\bP_h \bzeta)+ J_0(\bu-\bP_h \bu,\bP_h \bzeta))
			-b(\bP_h\bzeta,p-r_h(p)).\label{zetaener}
		\end{align}
		Rewrite the first term on the right hand side of (\ref{zetaener})  as
		\begin{align}
			\kappa  a(\bu_t-\bP_h\bu_t,\bP_h\bzeta)= &\kappa\sum_{T\in \mathcal{T}_h}\int_T \nabla (\bu_t-\bP_h\bu_t):\nabla (\bP_h\bzeta)\,dT\nonumber \\
			&-\kappa\sum_{e\in \Gamma_h} \int_e\{\nabla (\bu_t-\bP_h\bu_t)\}\bn_e\cdot  [\bP_h\bzeta]\,ds - \kappa\sum_{e\in \Gamma_h}\int_e \{\nabla (\bP_h\bzeta)\}\bn_e\cdot [\bu_t-\bP_h\bu_t]\,ds\nonumber\\
			=& I_1+I_2+I_3, \,\,(say).\label{z1}
		\end{align}
		A use of Cauchy-Schwarz's, Young's inequalities and  (\ref{P_h}) yield
		\begin{align*}
			|I_1|\leq \kappa\sum_{T \in\mathcal{T}_h}\|\nabla(\bu_t-\bP_h\bu_t)\|_{L^2(T)}\|\nabla(\bP_h\bzeta)\|_{L^2(T)}\leq \frac{C_2\nu}{24}\|\bP_h\bzeta\|_\varepsilon^2+Ch^{2}|\bu_t|_{2}^2.
		\end{align*}
		If the edge $e$ belongs to the element $T$, using Lemma \ref{trace} and Cauchy-Schwarz's inequality, we bound $I_2$ as
		\[\bigg|\int_e\{\nabla (\bu_t-\bP_h\bu_t)\}\bn_e\cdot [\bP_h\bzeta]\,ds\bigg|\leq C(\|\nabla(\bu_t-\bP_h\bu_t)\|_{L^2(T)}+h_T\|\nabla^2(\bu_t-\bP_h\bu_t)\|_{L^2(T)})\frac{1}{|e|^{1/2}}\| [\bP_h\bzeta]\|_{L^2(e)}. \]
		By introducting the standard Lagrange interpolant $L_h(\bu)$ of degree $1$ (\cite{GRW005}) and using the triangle inequality, inverse inequality of Lemma \ref{distrace}, we obtain
		\begin{align}\label{L1}
			\|\nabla^2(\bu_t-\bP_h\bu_t)\|_{L^2(T)}&\leq \|\nabla^2(\bu_t-L_h(\bu_t))\|_{L^2(T)}+\|\nabla^2(L_h(\bu_t)-\bP_h\bu_t)\|_{L^2(T)}\nonumber\\
			&\leq \|\nabla^2(\bu_t-L_h(\bu_t))\|_{L^2(T)}+Ch_T^{-1}\|\nabla(L_h(\bu_t)-\bP_h\bu_t)\|_{L^2(T)}.
		\end{align}
		 A use of the Young inequality, (\ref{P_h}), (\ref{L1}) and the  approximation properties of $L_h$ (\cite{GRW005}) leads to
		\[|I_2|\leq C \kappa hJ_0(\bP_h\bzeta,\bP_h\bzeta)^{1/2}|\bu_t|_{2}\leq  \frac{C_2\nu}{24}\|\bP_h\bzeta\|_\varepsilon^2+Ch^{2}|\bu_t|_{2}^2.\]
		Furthermore, a use of Cauchy-Schwarz's, Young's inequalities, (\ref{P_h}) and Lemma \ref{distrace} yield
		\begin{align*}
			|I_3|\leq C\kappa\bigg(\sum_{e\in\Gamma_h}\frac{|e|}{\sigma_e}\|\{\nabla(\bP_h\bzeta)\}\|_{L^2(e)}^2\bigg)^{1/2}\bigg(\sum_{e\in\Gamma_h}\frac{\sigma_e}{|e|}\|[\bu_t-\bP_h\bu_t]\|_{L^2(e)}^2\bigg)^{1/2}\leq \frac{C_2\nu}{24}\|\bP_h\bzeta\|_\varepsilon^2+Ch^{2}|\bu_t|_{2}^2. 
		\end{align*}
	 By a virtue of the Cauchy-Schwarz, Young inequalities and (\ref{P_h}), the jump term of the first term on the right hand side of (\ref{zetaener}) is bounded as  follows:
		\begin{align}
			\kappa\,|J_0(\bu_t-\bP_h\bu_t,\bP_h\bzeta)|&\leq  \kappa J_0(\bu_t-\bP_h\bu_t,\bu_t-\bP_h\bu_t)^{1/2}J_0(\bP_h\bzeta,\bP_h\bzeta)^{1/2}\nonumber\\
			&\leq  \frac{C_2\nu}{24}\|\bP_h\bzeta\|_\varepsilon^2+Ch^2|\bu_t|_{2}^2.\label{z2}
		\end{align}
		Owing to (\ref{2}), the term involving the pressure in (\ref{zetaener}) is reduced to
		\[b(\bP_h\bzeta,p-r_hp)=\sum_{e\in \Gamma_h}\int_e\{p-r_hp\}[\bP_h\bzeta]\cdot \bn_e\,ds.\]
		Apply the Cauchy-Schwarz inequality, the result of (\ref{2.1}) and Lemma \ref{trace} to arrive at
		\begin{align}\label{z3}
			|b(\bP_h\bzeta,p-r_hp)| \leq & C\sum_{T\in \mathcal{T}_h}(\|p-r_hp\|_{L^2(T)}+h_T\|\nabla(p-r_hp)\|_{L^2(T)})J_0(\bP_h\bzeta,\bP_h\bzeta)^{1/2}\nonumber\\
			\leq & \frac{C_2\nu}{6}\|\bP_h\bzeta\|_\varepsilon^2+Ch^2|p|_1^2.
		\end{align}
		Similarly, by replacing $\bu_t-P_h\bu_t$ with $\bu-P_h\bu$, and using the analysis involved in bounding $I_1$, $I_2$ and $I_3$ and (\ref{z2}), we estimate the second term on the right hand side of (\ref{zetaener}) as
		\begin{align}\label{Rn1}
			\nu |a(\bu-\bP_h\bu,{\bP_h \bzeta})+J_0 (\bu-\bP_h\bu,\bP_h \bzeta)|\leq  \frac{C_2\nu}{6}\|\bP_h\bzeta\|_\varepsilon^2+Ch^2|\bu|_{2}^2.
		\end{align}
		Apply (\ref{z1})-(\ref{Rn1}) and the bound of Lemma \ref{coer} in (\ref{zetaener}). Then, multiply the resulting equation by $e^{2\alpha t}$, integrate from $0$ to $t$, use Lemmas \ref{cont}, \ref{coer} and observe that $P_h\bzeta(0)=0$, we obtain
		\begin{align}\label{st1}
			\kappa\,C_2\,e^{2\alpha t}\|\bP_h\bzeta \|_\varepsilon^2+	 (C_2\nu-2\alpha \kappa\,C_1)\int_0^t e^{2\alpha s} \|\bP_h\bzeta(s) \|_\varepsilon^2 ds\leq   Ch^{2}\int_{0}^t e^{2\alpha s}
			\left(|\bu_s(s)|_2^2+|\bu(s)|_2^2+|p(s)|_1^2\right)ds.
		\end{align}
		Multiply (\ref{st1}) by $e^{-2\alpha t}$ and use the regularity estimates of $\bu$ and $p$ from Lemma \ref{exactpriori} to complete the energy norm estimates of $P_h\bzeta$ as
		\begin{equation}
			\|P_h\bzeta(t)\|_\varepsilon^2+e^{-2\alpha t}\int_0^t e^{2\alpha s} \|P_h\bzeta(s) \|_\varepsilon^2 ds\leq Ch^{2}.\label{9.71}
		\end{equation}
		Since $\bzeta=\bu-P_h\bu+P_h\bzeta$, using the  triangle inequality and the bounds in (\ref{P_h}), (\ref{9.71}), we arrive at
		\begin{equation}
			\|\bzeta(t)\|_\varepsilon^2+e^{-2\alpha t}\int_0^t e^{2\alpha s} \|\bzeta(s) \|_\varepsilon^2 ds\leq Ch^{2}.\label{9.7}
		\end{equation}
		To derive the estimates of $\bzeta_t$ in energy norm, we substitute $\bphi_h=\bP_h \bzeta_t$ in (\ref{E61}). Then, apply Cauchy-Schwarz's, Young's inequalities, Lemma \ref{coer}, (\ref{P_h}) and (\ref{9.7}) to the resulting equation and arrive at
		\begin{align}\label{9.7a}
			\|P_h\bzeta_t(t) \|_\varepsilon^2\leq Ch^2.
		\end{align}
		A multiplication of (\ref{9.7a}) by $e^{2\alpha t}$, an integration from $0$ to $t$ with respect to time, and then again a multiplication by $e^{-2\alpha t}$ lead to
		\begin{equation}
			e^{-2\alpha t}\int_0^t e^{2\alpha s} \|P_h\bzeta_t(s) \|_\varepsilon^2 \,ds\leq Ch^{2}.\label{9.71*}
		\end{equation}
		A use of the triangle inequality and bounds of (\ref{P_h}), (\ref{9.7a}), (\ref{9.71*}) yields
		\begin{equation}
		 \|\bzeta_t(t)\|_\varepsilon^2+	e^{-2\alpha t}\int_0^t e^{2\alpha s} \|\bzeta_t(s) \|_\varepsilon^2 \,ds\leq Ch^{2}.\label{9.7*}
		\end{equation}
		For obtaining the $L^2$-norm estimates of $\bzeta$, we apply the Aubin-Nitsche duality argument. 
		Let $\{\bw,q\}$ be the pair of unique solution of the steady state Stokes system stated as
		\begin{align}\label{9.9}
			-\nu \Delta \bw+\nabla q & =\bzeta \quad\text{in $\Omega$},\quad\nabla\cdot \bw = 0\quad \text{in $\Omega$},\quad\bw|_{\partial \Omega} =0.
		\end{align}
		The above solution pair satisfies the following regularity results:
		\begin{equation}
			\|\bw\|_2+\|q\|_1\leq C\|\bzeta\|.\label{10.0}
		\end{equation}
		Forming $L^2$ inner product between (\ref{9.9}) and $\bzeta$, and applying the regularity estimates of $\bw$ and $q$, we arrive at
		\begin{align*}
			\|\bzeta\|^2 = &\nu\sum_{T\in\mathcal{T}_h} \int_T\nabla \bw:\nabla \bzeta\,dT- \nu\sum_{T\in\mathcal{T}_h}\int_{\partial T}(\nabla \bw \bn_T)\cdot \bzeta\,ds \nonumber\\
			&-\sum_{T\in\mathcal{T}_h}\int_T q\nabla\cdot \bzeta\,dT+\sum_{T\in\mathcal{T}_h}\int_{\partial T}q\bn_T\cdot \bzeta\,ds\nonumber\\
			= & \nu \sum_{T\in\mathcal{T}_h} \int_T\nabla \bzeta:\nabla \bw\,dT - \nu\sum_{e\in \Gamma_h}\int_e \{\nabla \bw\}\bn_e\cdot [\bzeta]\,ds+b(\bzeta,q).
		\end{align*}
		Using (\ref{E61}) with $\bphi_h$ replaced by $\bP_h\bw$, and observing that on each interior edge $[\bw]\cdot \bn_e=0$, we obtain
		\begin{align}\label{e10}
			\|\bzeta(t)\|^2=&\nu \sum_{T\in\mathcal{T}_h} \int_T\nabla \bzeta:\nabla (\bw-\bP_h\bw)\,dT-\nu\sum_{e\in \Gamma_h}\int_e \{\nabla (\bw-\bP_h\bw)\}\bn_e\cdot [\bzeta]\,ds\nonumber \\
			&+\nu\sum_{e\in \Gamma_h}\int_e \{\nabla \bzeta\}\bn_e\cdot [\bP_h\bw-\bw]\,ds +\nu J_0(\bzeta,\bw-\bP_h \bw)\nonumber\\
			&+ \kappa\,(a(\bzeta_t,\bw-\bP_h \bw)+ J_0(\bzeta_t,\bw-\bP_h \bw))+b(\bzeta,q)-b(\bP_h \bw-\bw,p-r_h(p))\nonumber\\
			&- \kappa\,(a(\bzeta_t,\bw)+J_0(\bzeta_t,\bw)).
		\end{align} 
 Again, form an $L^2$ inner product between (\ref{9.9}) and $\bzeta_t$. Then,  using integration by parts and $[\bw]\cdot \bn_e=0$ on each interior edge, we rewrite the last term on the right-hand side of (\ref{e10}) and obtain
		\begin{align}\label{e10*}
			\|\bzeta(t)\|^2=&\nu \sum_{T\in\mathcal{T}_h} \int_T\nabla \bzeta:\nabla (\bw-\bP_h\bw)\,dT-\nu\sum_{e\in \Gamma_h}\int_e \{\nabla (\bw-\bP_h\bw)\}\bn_e\cdot [\bzeta]\,ds\nonumber \\
			&+\nu\sum_{e\in \Gamma_h}\int_e \{\nabla \bzeta\}\bn_e\cdot [\bP_h\bw-\bw]\,ds +\nu J_0(\bzeta,\bw-\bP_h \bw)\nonumber\\
			&+ \kappa\,(a(\bzeta_t,\bw-\bP_h \bw)+ J_0(\bzeta_t,\bw-\bP_h \bw))+b(\bzeta,q)-b(\bP_h \bw-\bw,p-r_h(p))\nonumber\\
			&+\frac{\kappa}{\nu}b(\bzeta_t,q)-\frac{\kappa}{\nu}(\bzeta,\bzeta_t).
		\end{align} 
		   We follow the similar steps as used in bounding the first term on the right hand side of (\ref{zetaener}) to estimate the first, second, third, fourth and fifth terms in the right hand side of (\ref{e10*}). Then, we apply (\ref{P_h}) and (\ref{10.0}) to arrive at
		\begin{align}
			\bigg|\nu \sum_{T\in\mathcal{T}_h} & \int_T\nabla \bzeta:\nabla (\bw-\bP_h\bw)\,dT-\nu \sum_{e\in \Gamma_h}\int_e \{\nabla (\bw-\bP_h\bw)\}\bn_e\cdot [\bzeta]\,ds\nonumber \\
			+\nu &\sum_{e\in \Gamma_h} \int_e \{\nabla \bzeta\}\bn_e\cdot [\bP_h\bw-\bw]\,ds +\nu J_0(\bzeta,\bw-\bP_h \bw)+ \kappa\,(a(\bzeta_t,\bw-\bP_h \bw)+ J_0(\bzeta_t,\bw-\bP_h \bw))\bigg|\nonumber\\
			&\leq  Ch\|\bw\|_2\|\bzeta\|_\varepsilon+Ch^{2}|\bu|_{2}\|\bw\|_2+ Ch\|\bw\|_2\|\bzeta_t\|_\varepsilon+ Ch^{2}|\bu_t|_{2}\|\bw\|_2\nonumber\\
			&\leq  \frac{1}{8}\|\bzeta\|^2+Ch^2(\|\bzeta\|^2_\varepsilon
			+\|\bzeta_t\|^2_\varepsilon)+Ch^{4}(|\bu|^2_{2}+ |\bu_t|^2_{2}).\label{z4}
		\end{align}
		We can handle the sixth term on the right-hand side of (\ref{e10*}) as
		\begin{align}
			b(\bzeta,q)&= b(\bzeta-\bP_h\bu+\bS_h\bu,q)+b(\bP_h\bzeta,q)= b(\bu-\bP_h \bu,q)+b(\bP_h\bzeta,q-r_h(q)) \nonumber \\
			&= -\sum_{T\in \mathcal{T}_h}\int_T q\nabla\cdot (\bu-\bP_h\bu)\,dT+\sum_{e\in \Gamma_h}\int_e\{q\}[\bu-\bP_h \bu]\cdot \bn_e\,ds+b(\bP_h\bzeta,q-r_h(q)) \label{z5}.
		\end{align}
		Furthermore, using integration by parts formula to the first term on the right hand side of (\ref{z5}) and noting that $q$ is continuous, we arrive at
		\begin{align*}
			b(\bzeta,q)= \sum_{T\in\mathcal{T}_h}\int_T \nabla q\cdot (\bu-\bP_h\bu)\,dT+b(\bP_h\bzeta,q-r_h(q)). 
		\end{align*}
		From Cauchy-Schwarz's, Young's inequalities, (\ref{2.1}), (\ref{P_h}), (\ref{10.0}) and Lemma \ref{trace},  we obtain 
		\begin{align}
			|b(\bzeta,q)|  &\leq \bigg| Ch^{2}|q|_1|\bu|_{2}-\sum_{T\in\mathcal{T}_h}\int_T(\nabla \cdot \bP_h\bzeta)( q-r_h(q))\,dT+\sum_{e\in\Gamma_h}\int_e\{q-r_h(q)\}[\bP_h\bzeta]\cdot \bn_e\,ds\bigg|\nonumber\\
			& \leq Ch^{2}|\bu|_{2}\|\bzeta\|+Ch|q|_1\|\bP_h\bzeta\|_\varepsilon\leq \frac{1}{8}\|\bzeta\|^2+ Ch^2(h^2|\bu|^2_{2}+\|\bP_h\bzeta\|^2_\varepsilon).\label{z6}
		\end{align}
		Similar to (\ref{z6}), using (\ref{2.1}), (\ref{P_h}), (\ref{10.0}) and Lemma \ref{trace}, the 8th term on the right hand side of (\ref{e10*}) can be bounded as follows:
		\begin{align}
			|b(\bzeta_t,q)|  &\leq \bigg| Ch^{2}|q|_1|\bu_t|_{2}-\sum_{T\in\mathcal{T}_h}\int_T(\nabla \cdot \bP_h\bzeta_t)(q-r_h(q))\,dT+\sum_{e\in\Gamma_h}\int_e\{q-r_h(q)\}[\bP_h\bzeta_t]\cdot \bn_e\,ds\bigg|\nonumber\\
			& \leq Ch^{2}|\bu_t|_{2}\|\bzeta\|+Ch|q|_1\| \bP_h\bzeta_t\|_\varepsilon\leq \frac{1}{8}\|\bzeta\|^2+Ch^2(h^2|\bu_t|^2_{2}+\|\bP_h\bzeta_t\|^2_\varepsilon).\label{z61}
		\end{align}
	Apply Cauchy-Schwarz's inequality and (\ref{2.1}), (\ref{P_h}), (\ref{10.0}) to arrive at
		\begin{align}
			|b(\bP_h \bw-\bw,p-r_h(p))|
			\leq  Ch^{2}|p|_1\|\bw\|_2\leq\frac{1}{8}\|\bzeta\|^2+Ch^{4}|p|^2_1.\label{z7}
		\end{align}
		%
		A use of (\ref{z4}) and (\ref{z6})--(\ref{z7}) in (\ref{e10*}) leads to
		\begin{align}\label{z71}
		\frac{1}{2}	\|\bzeta(t)\|^2+\frac{\kappa}{2\nu}\frac{d}{dt}	\|\bzeta(t)\|^2\leq &  Ch^2(\|\bzeta\|_\varepsilon^2+\|\bzeta_t\|_\varepsilon^2+\|\bP_h\bzeta\|_\varepsilon^2+\| \bP_h\bzeta_t\|_\varepsilon^2)+ Ch^{4}(|\bu|^2_{2}+|\bu_t|^2_{2}+|p|^2_1).
		\end{align}
		A multiplication of (\ref{z71}) by $e^{2\alpha t}$ and an integration of the resulting equation with respect to time from $0$ to $t$ yield
		\begin{align}\label{z72}
			&e^{2\alpha t}\|\bzeta(t)\|^2+	\left(\frac{\nu-2 \kappa\alpha}{\nu}\right)\displaystyle{\int_0^t}e^{2\alpha s}	\|\bzeta(s)\|^2 ds	\leq  Ch^{4}|\bu_0|_2^2\\
				+Ch^2\displaystyle{\int_0^t} & e^{2\alpha s}\left(\|\bzeta(s)\|_\varepsilon^2+\|\bzeta_t(s)\|_\varepsilon^2
			+\|\bP_h\bzeta(s)\|_\varepsilon^2+\| \bP_h\bzeta_t(s)\|_\varepsilon^2\right)\,ds +Ch^4\displaystyle{\int_0^t}e^{2\alpha s}(|\bu(s)|^2_{2}+|\bu_t(s)|^2_{2}+|p(s)|^2_1)\,ds.\nonumber
		\end{align}
		Multiply (\ref{z72}) by $e^{-2\alpha t}$ and use (\ref{9.71}), (\ref{9.7}),  (\ref{9.71*}), (\ref{9.7*}) with Lemma \ref{exactpriori}  to arrive at
		\begin{align}\label{z72n}
			\|\bzeta(t)\|^2+e^{-2\alpha t}\displaystyle{\int_0^t}e^{2\alpha s}	\|\bzeta(s)\|^2 ds	\leq Ch^{4}.
		\end{align}
		A combination of (\ref{9.7}), (\ref{9.7*}) and (\ref{z72n}) completes the proof of (\ref{9.4}) in Lemma \ref{zetaerr}.
		
\noindent
		Following the similar steps as involved in proving the $\bL^2$ estimate of $\bzeta$ in (\ref{9.4}), we arrive at the  $\bL^2$ estimate (\ref{9.5}) involving $\bzeta_t$.  Only difference is in the dual problem, where the right hand side is changed to $\bzeta_t$. With the resulting $\bL^2$ estimate of $\bzeta_t$, we conclude the proof of Lemma \ref{zetaerr}.
		
	\end{proof}

\noindent
Below, in Lemma \ref{rhoerror}, we derive the bounds of $\brho$.
	\begin{lemma} \label{rhoerror}
	Let the assumptions ({\bf A1})-({\bf A3}) be satisfied and let  $0<\al<\min \left\{ \frac{\nu C_2}{2(\gamma+\kappa C_1)},\frac{\nu}{2\kappa}\right\}$. Then, there exists a positive constant $C=C(\kappa,\nu,\alpha,C_2, M_0)$, such that, for each $t>0$,  the following estimates hold true:
\begin{align*}
	\|\brho(t)\|^2+h^2\|\brho(t)\|^2_\varepsilon+e^{-2\alpha t}	\displaystyle{\int_0^t}e^{2\alpha s}\left(\|\brho(s)\|^2+h^2\|\brho(s)\|^2_\varepsilon\right)ds\leq &  Ch^4.
	\end{align*}
\end{lemma}

	\begin{proof}
	Subtract (\ref{E61}) from (\ref{9.2}) and write the equation in $\brho$ as
	\[(\brho_t,\bphi_h)+\kappa\,(a(\brho_t,\bphi_h)+J_0(\brho_t,\bphi_h))+\nu(a(\brho,\bphi_h)+J_0(\brho,\bphi_h))=-(\bzeta_t,\bphi_h),\quad \forall \bphi_h\in\bJ_h.\]
	Substitute $\bphi_h=\brho$  in the above equation and use Lemma \ref{coer} to obtain
	\begin{align}
		\frac{1}{2}\frac{d}{dt}\left(\|\brho\|^2+ \kappa\,(a(\brho,\brho)+J_0(\brho,\brho))\right)+\nu C_2\|\brho\|^2_\varepsilon\leq -(\bzeta_t,\brho). \label{theta1}
	\end{align}
	Multiply (\ref{theta1}) by $e^{2\alpha t}$, integrate the resulting  inequality with respect to time from $0$ to $t$, and use Lemmas \ref{coer} and \ref{cont}, $\brho(0)=0$ to arrive at
	\begin{align}
		e^{2\alpha t}\left(\|\brho\|^2+\kappa\,C_2\|\brho\|^2_\varepsilon\right)+ (\nu C_2-2\al(\gamma+\kappa C_1) )&
		\int_0^te^{2\alpha s}\|\brho(s)\|^2_\varepsilon ds \leq C\int_0^t  e^{2\alpha s}\|\bzeta_s(s)\|^2\,ds.\label{theta3}
	\end{align}
 A multiplication of (\ref{theta3}) by $e^{-2\alpha t}$ and a use of (\ref{9.5}) in the resulting inequality yield
	\begin{align}
		\|\brho(t)\|^2+\kappa\,C_2\|\brho(t)\|^2_\varepsilon+e^{-2\alpha t}
		\int_0^te^{2\alpha s}\|\brho(s)\|^2_\varepsilon ds\leq Ch^4.\label{theta4}
	\end{align}
 An application of (\ref{T6}) and (\ref{theta4}) leads to 
		\begin{align}\label{theta5}
		h^2\|\brho(t)\|^2_\varepsilon+e^{-2\alpha t}	\int_0^te^{2\alpha s}\|\brho(s)\|^2 ds\leq Ch^4.
	\end{align}
A combination of (\ref{theta4}) and (\ref{theta5}) concludes the proof of Lemma \ref{rhoerror}.
\end{proof}
\noindent
 Since $\bxi=\bzeta+( \bS_h\bu-\bv_h)=\bzeta+\brho$, we now apply Lemmas \ref{zetaerr} and \ref{rhoerror} along with the triangle inequality to obtain  the following estimates of $\bxi$.
		\begin{equation}
			\|\bxi(t)\|^2+h^2\|\bxi(t)\|^2_\varepsilon+e^{2\alpha t}\displaystyle{\int_0^t}e^{2\alpha s}	\|\bxi(s)\|^2ds\leq Ch^{4},\quad 0\leq t\leq T.\label{xierr}
		\end{equation}
	\noindent  We are now left with the estimates of $\be=\bv_h-\bu_h$.
	\begin{lemma}\label{etaerror}
	Let the assumptions ({\bf A1})-({\bf A3}) be satisfied and let $0<\al<\min \left\{ \frac{\nu C_2}{2(\gamma+\kappa C_1)},\frac{\nu}{2\kappa}\right\}$. Further, let $\bv_h(t) \in \bJ_h$ be a solution of (\ref{9.1}) corresponding to the initial value $\bv_h(0)=\bP_h\bu_0$. Then, there exists a positive constant $C=C(\kappa,\nu,\alpha,C_2, M_0)$, independent of $h$, such that, the following estimates hold true:
		\begin{equation*}
			\|\be(t)\|^2+\,\|\be(t)\|_\varepsilon^2+e^{-2\alpha t}\int_0^te^{2\al s}\|\be(s)\|_\varepsilon^2\,ds\leq Ce^{Ct}h^4.
		\end{equation*}
	\end{lemma}
		\begin{proof}
		From (\ref{ne1}) and (\ref{9.1}), we observe that
		\[ (\be_t,\bphi_h)+\kappa\,(a(\be_t,\bphi_h)+J_0(\be_t,\bphi_h))+\nu(a(\be,\bphi_h)+J_0(\be,\bphi_h))=c^{\bu_h}(\bu_h,\bu_h,\bphi_h)-c^\bu(\bu,\bu,\bphi_h) \quad \text{for $ \bphi_h\in \bJ_h$}.\]
		Substitute $\bphi_h=\be$ and use Lemma \ref{coer} to arrive at
		\begin{equation}\label{10.6}
			\frac{1}{2}\frac{d}{dt}\left(\|\be\|^2+ \kappa\,(a(\be,\be)+J_0(\be,\be))\right)+\nu C_2\|\be\|_\varepsilon^2\leq c^{\bu_h}(\bu_h,\bu_h,\be)-c^\bu(\bu,\bu,\be). 
		\end{equation}
		Since $\bu$ is continuous, we can rewrite the nonlinear term as
		\[c^\bu(\bu,\bu,\be)=c^{\bu_h}(\bu,\bu,\be).\]
		For the sake of simplicity, we drop the superscript in the nonlinear terms unless there is no confusion.
		We now rewrite the nonlinear terms  as follows:
		\begin{align}
			c(\bu_h,\bu_h,\be)-c(\bu,\bu,\be) =-c(\bu,\bxi,\be)+c(\bxi,\bxi,\be)-c(\bxi,\bu,\be)+c(\be,\bxi,\be)-c(\be,\bu,\be)-c(\bu_h,\be,\be).
		\end{align}
		Note that, the last term is non-negative due to (\ref{5.5}) and is therefore dropped.  We find that 	
		\begin{align}		
		c(\bu_h,\bu_h,\be)-c(\bu,\bu,\be)\le -c(\bu,\bxi,\be)+c(\bxi,\bxi,\be)-c(\bxi,\bu,\be)+c(\be,\bxi,\be)-c(\be,\bu,\be). \label{c0}
		\end{align}
		Apply (\ref{tri1}), Lemma \ref{sobolev1} and Young's inequality to obtain
		\begin{align}
			|c(\be,\bu,\be)|\leq C \|\be\|\|\bu\|_{2}\|\be\|_\varepsilon\leq  \frac{C_2\nu}{16}\|\be\|^2_{\varepsilon}+C\|\bu\|_2^2\|\be\|^2.\label{c4}
		\end{align}
		We bound the second and fourth nonlinear terms on the right hand side of (\ref{c0}) using the Cauchy-Schwarz, Young inequalities, (\ref{Lp}) and Lemmas \ref{trace} and \ref{distrace} as follows
		\begin{align}\label{c5}
			|c (\bxi,\bxi,\be)|= &\bigg|\sum_{T\in\mathcal{T}_h} \int_T (\bxi\cdot \nabla\bxi)\cdot \be\,dT+\sum_{T\in\mathcal{T}_h}\int_{\partial T_-}|\{\bxi\}\cdot \bn_T|(\bxi^{int}-\bxi^{ext})\cdot \be^{int}\,ds +\frac{1}{2}\sum_{T\in\mathcal{T}_h} \int_T(\nabla\cdot \bxi)\bxi\cdot \be\,dT \nonumber\\
			&-\frac{1}{2}\sum_{e\in\Gamma_h}\int_e[\bxi]\cdot \bn_e\{\bxi\cdot\be\}\,ds\bigg|\nonumber\\
			\leq  \sum_{T\in\mathcal{T}_h}\|\bxi\|_{L^4(T)}& \|\nabla\bxi\|_{L^2(T)}\|\be\|_{L^4(T)}+C\sum_{T\in\mathcal{T}_h}h_T^{-3/4}(\|\bxi\|_{L^2(T)}+h_T\|\nabla\bxi\|_{L^2(T)})|e|^{1/4}|e|^{-1/2}\|[\bxi]\|_{L^2(e)}\|\be\|_{L^4(T)}\nonumber\\
			+C\sum_{T\in\mathcal{T}_h}\|\nabla\bxi &\|_{L^2(T)}\|\bxi\|_{L^4(T)}\|\be\|_{L^4(T)}+C\sum_{T\in\mathcal{T}_h}|e|^{1/4}|e|^{-1/2}\|[\bxi]\|_{L^2(e)}h_T^{-3/4}(\|\bxi\|_{L^2(T)}+h_T\|\nabla\bxi\|_{L^2(T)})\|\be\|_{L^4(T)}\nonumber\\
			&\leq    \frac{C_2\nu}{16}\|\be\|^2_\varepsilon+C h^{-2}(\|\bxi\|^2+h^2\|\bxi\|_\varepsilon^2)\|\bxi\|_\varepsilon^2
		\end{align}
		and
		\begin{align}
			|c(\be,\bxi,\be)|= &\bigg|\sum_{T\in\mathcal{T}_h}\int_T (\be\cdot \nabla\bxi)\cdot \be\,dT+\sum_{T\in\mathcal{T}_h}\int_{\partial T_-}|\{\be\}\cdot \bn_T|(\bxi^{int}-\bxi^{ext})\cdot \be^{int}\,ds +\frac{1}{2}\sum_{T\in\mathcal{T}_h} \int_T(\nabla\cdot \be)\bxi\cdot \be\,dT \nonumber\\
			&-\frac{1}{2}\sum_{e\in\Gamma_h}\int_e[\be]\cdot \bn\{\bxi\cdot\be\}\,ds\bigg|\nonumber\\
			\leq  \sum_{T\in\mathcal{T}_h}\|\be\|_{L^4(T)}&\|\nabla \bxi\|_{L^2(T)}\|\be\|_{L^4(T)}+C \sum_{T\in\mathcal{T}_h}\|\be\|_{L^4(T)}\frac{1}{|e|^{1/2}}\|[\bxi]\|_{L^2(e)}\|\be\|_{L^4(T)}\nonumber\\
			+C\sum_{T\in\mathcal{T}_h}\|\nabla \be &\|_{L^2(T)}\|\bxi\|_{L^4(T)}\|\be\|_{L^4(T)} +C\sum_{T\in\mathcal{T}_h} h_T^{-3/4}\|\be\|_{L^2(T)} h_T^{-1/2}(\|\bxi\|_{L^2(T)}+h_T\|\nabla \bxi\|_{L^2(T)})h_T^{-1/4}\|\be\|_{L^4(T)}\nonumber\\
			& \leq  \frac{C_2\nu}{16}\|\be\|^2_\varepsilon+C h^{-4}(\|\bxi\|^2+h^2\|\bxi\|_\varepsilon^2)\|\be\|^2.\label{c3}
		\end{align}
		By	applying (\ref{integration}), the first nonlinear term on the right hand side of (\ref{c0}) can be written as
		\begin{align*}
			c(\bu;\bxi,\be)=& -\sum_{T\in\mathcal{T}_h}\int_T(\bu\cdot \nabla \be)\cdot \bxi\,dT-\frac{1}{2}\sum_{T\in\mathcal{T}_h}\int_T(\nabla\cdot \bu)\bxi\cdot \be\,dT\nonumber\\
			& +\frac{1}{2}\sum_{e\in\Gamma_h}\int_e[\bu]\cdot \bn_e\{\bxi\cdot\be\}\,ds-\sum_{T\in \mathcal{T}_h}\int_{\partial T_-}|\bu\cdot \bn_T|\bxi^{ext}\cdot (\be^{int}-\be^{ext})\,ds\nonumber\\
			& +\int_{\Gamma_+}|\bu\cdot \bn|\bxi\cdot \be \,ds\nonumber\\
			&=I_1+I_2+I_3+I_4+I_5, (\text{say}).
		\end{align*}
		Note that, $\bu$ is continuous. This will lead to $I_3=I_5=0$. An application of Cauchy-Schwarz's, Young's inequalities, Lemma \ref{sobolev1} and (\ref{Lp}) yield
		\begin{align}
			|I_1|+|I_2|&\leq C \|\bu\|_{L^\infty(\Omega)}\|\nabla \be\|\|\bxi\|+ C\|\nabla \bu\|_{L^4(\Omega)}\|\bxi\|\|\be\|_{L^4(\Omega)} \nonumber\\
			&\leq C\|\bu\|_2\|\bxi\|\|\be\|_\varepsilon\leq C\|\bu\|_2^2\|\bxi\|^2+\frac{C_2\nu}{16}\|\be\|^2_\varepsilon \label{c11}
		\end{align}
		\noindent A use of Lemmas \ref{sobolev1}, \ref{trace} and Young's inequality leads to 
		\begin{align}
			|I_4|& \leq C\|\bu\|_{L^\infty(\Omega)}\sum_{T\in \mathcal{T}_h}\|\bxi\|_{L^2(\partial T)}|e|^{1/2-1/2}\|[\be]\|_{L^2(\partial T)}\nonumber
			\\
			&\leq C\|\bu\|_2\sum_{T\in \mathcal{T}_h}|e|^{1/2}h_T^{-1/2}(\|\bxi\|_{L^2(T)}+h_T\|\nabla\bxi\|_{L^2(T)})\frac{1}{|e|^{1/2}}\|[\be]\|_{L^2(\partial T)}\nonumber\\
			& \leq C\|\bu\|_2^2(\|\bxi\|^2+h^2\|\bxi\|_\varepsilon^2)+\frac{C_2\nu}{16}\|\be\|^2_\varepsilon.\label{c13}
		\end{align}
		 Finally for the third term on the right hand side of (\ref{c0}), we first note that
		\begin{align*}
			\int_T(\nabla \cdot \bu)\bv\cdot \bw\,dT = &\int_T \nabla\cdot (\bv\otimes \bu)\cdot \bw\,dT-\int_T(\bu\cdot \nabla \bv)\cdot \bw\,dT\\
			=&-\int_T(\bu\cdot \nabla \bw)\cdot \bv\,dT+\int_{\partial T}(\bu^{int} \cdot \bn_T)\bv^{int}\cdot \bw^{int}\,ds -\int_T(\bu\cdot \nabla \bv)\cdot \bw\,dT,
		\end{align*}
		 by applying $(\nabla\cdot \bu)v_i=\nabla\cdot (\bu v_i)-\bu\cdot \nabla v_i$ and based on it, we rewrite the third term as
		\begin{align*}
			c(\bxi, \bu, \be)=& \frac{1}{2}\sum_{T\in\mathcal{T}_h}\int_T (\bxi\cdot \nabla\bu)\cdot \be\,dT- \frac{1}{2}\sum_{T\in\mathcal{T}_h}\int_T (\bxi\cdot \nabla \be)\cdot \bu\,dT\nonumber\\
			&- \sum_{T\in\mathcal{T}_h}\int_{\partial T_-} |\{\bxi\}\cdot \bn_T|(\bu^{int}-\bu^{ext})\cdot \be^{int}\,ds-\frac{1}{2}\sum_{e\in\Gamma_h}\int_e[\bxi]\cdot \bn_e\{\bu\cdot \be\}\,ds\nonumber\\
			&+ \frac{1}{2}\sum_{T\in \mathcal{T}_h}\int_{\partial T} (\bxi^{int}\cdot \bn_T)\bu^{int}\cdot \be^{int} \,ds\nonumber\\
			= & I_6+I_7+I_8+I_9+I_{10}.
		\end{align*}
		The term $I_8$ is zero due to the continuity of $\bu$. A use of (\ref{Lp}), Young's inequality and Lemma \ref{sobolev1}, bounds $I_{6}$ and $I_{7}$ as 
		\begin{align}
			|I_{6}| & \leq \|\bxi\|\|\|\nabla\bu\|_{L^4(\Omega)}\|\be\|_{L^4(\Omega)}  \leq  \frac{C_2\nu}{16}\|\be\|^2_\varepsilon+C\|\bu\|_2^2\|\bxi\|^2, \label{c21} \\
			|I_{7}|&\leq \|\bxi\|\|\|\nabla\be\|\|\bu\|_{L^\infty(\Omega)} \leq  \frac{C_2\nu}{16}\|\be\|^2_\varepsilon+C\|\bu\|_2^2\|\bxi\|^2. \label{c22}
		\end{align}
		Next, sum the integral in $I_{10}$ over all $T$ and consider the contribution of this sum to one interior edge $e$. Let us assume that $e$ is shared by two triangles $T_r$ and $T_s$, with exterior normal $\bn_r$ and $\bn_s$. Then, we arrive at
		\[\int_e (\bxi|_{T_r}\cdot \bn_r)\bu|_{T_r}\cdot \be|_{T_r}\,ds+\int_e(\bxi|_{T_s}\cdot \bn_s)\bu|_{T_s}\cdot \be|_{T_s}\,ds=\int_e[(\bxi\cdot \bn_e)\bu\cdot \be]\,ds. \]
		An application of Lemma \ref{trace} leads to
		\begin{align}
			I_{9}+I_{10}&=\frac{1}{2}\sum_{e\in\Gamma_h}\int_e\{\bxi\}\cdot \bn_e[\bu\cdot \be]\,ds \leq C \|\bu\|_{L^\infty(\Omega)}\sum_{e\in \Gamma_h}\frac{1}{|e|^{1/2}}\|[\be]\|_{L^2(e)}|e|^{1/2}\|\bxi\|_{L^2(e)}\nonumber\\
			& \leq C\|\bu\|_2 \sum_{T\in \mathcal{T}_h}\|\be\|_\varepsilon(\|\bxi\|_{L^2(T)}+h_T\|\nabla\bxi\|_{L^2(T)}) \leq  \frac{C_2\nu}{16}\|\be\|^2_\varepsilon+C \|\bu\|_2^2(\|\bxi\|^2+h^2\|\bxi\|_\varepsilon^2).\label{c23}
		\end{align}
		\noindent
		Substitute  (\ref{c4})-(\ref{c23}) in (\ref{c0}), and thereby in (\ref{10.6}), and multiply the resulting inequality  by $e^{2\al t}$ to obtain
		\begin{align}
		\frac{1}{2}	\frac{d}{dt} &\left(e^{2\al t}(\|\be\|^2+ \kappa\,(a(\be,\be)+J_0(\be,\be)))\right)+ \left(\frac{\nu C_2}{2}-2\al(\gamma+\kappa C_1)\right ) e^{2\al t}\|\be\|_\varepsilon^2\nonumber\\
			&\leq  C(\|\bu\|_2^2+h^{-4}\|\bxi\|^2+h^{-2}\|\bxi\|_\varepsilon^2)e^{2\al t}(\|\be\|^2+\kappa\,\|\be\|_\varepsilon^2)
			+Ce^{2\al t}(\|\bu\|_2^2+h^{-2}\|\bxi\|_\varepsilon^2)(\|\bxi\|^2+h^2\|\bxi\|_\varepsilon^2).\label{c6}
		\end{align}
		An integration of (\ref{c6}) with respect to time from $0$ to $t$,  a use of 
		$\be(0)=0$, Lemmas \ref{coer} and \ref{cont}, Gronwall's inequality, (\ref{xierr}) and Lemma \ref{exactpriori} leads to the following estimates of $\be$
		\begin{align*}
			e^{2\al t}(\|\be(t)\|^2+\kappa\,C_2\|\be(t)\|_\varepsilon^2)+\nu C_2\int_0^te^{2\al s}\|\be\|_\varepsilon^2\,ds\leq Ce^{Ct}h^4e^{2\alpha t}.
		\end{align*}
		Multiply the resulting inequality by $e^{-2\al t}$ to arrive at the desired estimates. 
	\end{proof}

	\noindent
{\it Proof of Theorem \ref{semierror}.}  A use of $\e=\bxi+\be$, triangle's inequality, the estimates in (\ref{xierr}) and Lemma \ref{etaerror} completes the proof. \hfill{$\Box$}
	
\begin{Remark}
	Under the following smallness assumption on the data
	\begin{equation}
		N=\sup_{\bv_h,\bw_h\in\bV_h}\frac{c(\bw_h,\bv_h,\bw_h)}{\|\bw_h\|_\varepsilon^2\|\bv_h\|_\varepsilon}\quad \text{and}\quad \frac{N\gamma}{C_2^2\nu^2}\|\brf\|<1,\label{uniqueness}
	\end{equation}
	the bounds of Theorem \ref{semierror} can be shown to be uniform in time, that is,
	\begin{align*}
		\|(\bu-\bu_h)(t)\|+h\|(\bu-\bu_h)(t)\|_\varepsilon\leq Ch^2,
	\end{align*}
	where the constant $C$ is independent of time $t$ and $h$.
\end{Remark}
\noindent
Firstly, rewrite the nonlinear terms as follows:
\begin{align}
	c^{\bu_h}(\bu_h,\bu_h,\be)-c^{\bu_h}(\bu,\bu,\be)=- c^{\bu_h}(\bu,\bxi,\be)+c^{\bu_h}(\bxi,\bxi,\be)-c^{\bu_h}(\bxi,\bu,\be)-c^{\bu_h}(\be,\bu_h,\be)-c^{\bu_h}(\bv_h,\be,\be).\label{nl1}
\end{align}
The last nonlinear term on the right hand side of (\ref{nl1}) is positive due to the positivity property (\ref{5.5}). Using the proof of the Lemma \ref{etaerror}, the first three nonlinear terms on the right hand side of (\ref{nl1}) can be bounded as
\[|c^{\bu_h}(\bu,\bxi,\be)+c^{\bu_h}(\bxi,\bu,\be)+c^{\bu_h}(\bxi,\bxi,\be)|\leq C(\|\bu\|_2+ch^{-1}\|\bxi\|_\varepsilon)(\|\bxi\|+h\|\bxi\|_\varepsilon)\|\be\|_\varepsilon,\]
For the fourth nonlinear term on the right hand side of (\ref{nl1}), we use the smallness condition and find that
\[c^{\bu_h}(\be,\bu_h,\be)\leq N\|\be\|_\varepsilon^2\|\bu_h\|_\varepsilon.\]
Now,  in order to achieve the uniform in time velocity error estimates, we modify the  proof of Lemma \ref{etaerror} as follows:
From (\ref{10.6}) and using (\ref{xierr}), Lemma \ref{exactpriori}, we arrive at 
\begin{align}
	\frac{d}{dt}(\|\be\|^2+\kappa\,(a(\be,\be)+J_0(\be,\be)))+2(\nu C_2-N\|\bu_h\|_\varepsilon)\|\be\|_\varepsilon^2\leq  Ch^2\|\be\|_\varepsilon. \label{uni1}
\end{align}
After multiplying (\ref{uni1}) by $e^{2\al t}$, integrate with respect to time from $0$ to $t$. Then, multiply the resulting equation by $e^{-2\al t}$ and use Lemma \ref{cont} to find that
\begin{align*}
	\|\be\|^2&+ \kappa\,(a(\be,\be)+J_0(\be,\be))+2e^{-2\al t}\int_0^te^{2\al s}(\nu C_2-N\|\bu_h\|_\varepsilon)\|\be(s)\|_\varepsilon^2\,ds\leq e^{-2\al t}(\|\be(0)\|^2+ C_1\kappa \|\be(0)\|^2_\varepsilon) \nonumber\\
	&+2\al e^{-2\al t}\int_0^te^{2\al s}(\|\be(s)\|^2+ \kappa\,(a(\be(s),\be(s))+J_0(\be(s),\be(s))))\,ds+Ch^2e^{-2\al t}\int_0^te^{2\al s}\|\be(s)\|_\varepsilon\, ds.
\end{align*}
Letting $t\rightarrow\infty$, apply the L'Hospital rule and use Lemma \ref{uhpriori} to arrive at
\[\frac{1}{\al}\left(\nu C_2-\frac{N \gamma}{C_2\nu}\|\brf\|_{L^\infty(0,\infty;L^2(\Omega))}\right)\limsup_{t\rightarrow \infty}\|\be(t)\|_\varepsilon^2\leq \frac{Ch^2}{\al}\limsup_{t\rightarrow \infty}\|\be(t)\|_\varepsilon.\]
Owing to the smallness condition (\ref{uniqueness}), we arrive at
\[\limsup_{t\rightarrow \infty}\|\be(t)\|_\varepsilon\leq Ch^2.\]
Hence,
\begin{align}\label{eta1}
	\limsup_{t\rightarrow \infty}\|\be(t)\|\leq Ch^2.
\end{align}
Combining the estimates of (\ref{xierr}) and (\ref{eta1}), we find that
\[\limsup_{t\rightarrow \infty}\|\bu(t)-\bu_h(t)\|\leq Ch^2.\]

		\section{Error estimates for pressure}\label{s5}
	\se
	This section presents the derivation of error estimates for the semidiscrete discontinuous Galerkin approximation of the pressure.  We begin by proving a lemma which is crucial for establishing these error estimates.

\begin{lemma}\label{etLinfL2}
	Let the assumptions ({\bf A1})-({\bf A3}) be satisfied and let $0<\al<\min \left\{ \frac{\nu C_2}{2(\gamma+\kappa C_1)},\frac{\nu}{2\kappa}\right\}$. Then, the error $\de=\bu-\bu_h$ in approximating the velocity satisfies for $t>0$
	 \begin{equation*}
		\|\de_t(t)\|+\kappa\,\, \|\de_t(t)\|_\varepsilon\leq Ce^{Ct}h. 
	\end{equation*}
\end{lemma}
\begin{proof}%
	Subtract (\ref{ne1}) from (\ref{8.5}) to write the equation in error $\e=\bu-\bu_h$ as
	\begin{align}
		(\de_t,\bphi_h)&+\kappa\,(a(\de_t,\bphi_h)+J_0(\de_t,\bphi_h))+\nu(a(\de,\bphi_h)\nonumber\\
		&+J_0(\de,\bphi_h))+c^{\bu}(\bu,\bu,\bphi_h)-c^{\bu_h}(\bu_h,\bu_h,\bphi_h)+b(\bphi_h,p)=0\quad \forall \bphi_h\in  \bJ_h.\label{et4}
	\end{align}
	  Choose
	$\bphi_h=\bP_h\de_t=\de_t-(\bu_t-\bP_h\bu_t)$ in (\ref{et4}) and and use Lemma \ref{coer} to observe that
	\begin{align}
		\|\de_t\|^2+\kappa\,C_2 \| \bP_h\de_t\|_\varepsilon^2 \leq & (\de_t,\bu_t-P_h\bu_t)+\kappa\,(a(\bP_h\bu_t-\bu_t,\bP_h\de_t)+J_0(\bP_h\bu_t-\bu_t,\bP_h\de_t))\nonumber\\
		&-\nu(a(\e,\bP_h\de_t)+J_0(\e,\bP_h\de_t)) +c^{\bu_h}(\bu_{h},\bu_h,\bP_h\de_t)-c^{\bu}(\bu,\bu,\bP_h\de_t) -b(\bP_h\de_t,p).\label{et5}
	\end{align}
 Drop the superscripts of the nonlinear terms in (\ref{et5}) as $\bu$ is continuous. Now, rewrite the nonlinear terms as
	\begin{align}\label{nw1}
		c(\bu,\bu,\bP_h\de_t)-c(\bu_h,\bu_h,\bP_h\de_t)=-c(\de,\de,\bP_h\de_t)+c(\de,\bu,\bP_h\de_t)+c(\bu,\de,\bP_h\de_t).
	\end{align}
	Similar to Lemma \ref{etaerror}, the nonlinear term $c(\de,\de,\bP_h\de_t)$ can be bounded by  applying (\ref{Lp}), Lemmas \ref{trace}, \ref{distrace}, and Theorem \ref{semierror} as follows
	\begin{align}
		|c(\de,\de,\bP_h\de_t)|\leq C\|\de\|_\varepsilon \|\bP_h\de_t\|_\varepsilon\leq   \frac{\kappa C_2}{12}\|\bP_h\de_t\|_\varepsilon^2+C\|\de\|_\varepsilon^2.\label{pr1}
	\end{align}
	Using the continuity of $\bu$, Lemmas \ref{exactpriori}, \ref{trace} and \ref{sobolev1}, we arrive at 
	\begin{align}
		|c(\de,\bu,\bP_h\de_t)|=&\bigg|\sum_{T\in\mathcal{T}_h}\int_T(\de\cdot \nabla \bu)\cdot \bP_h\de_t\,dT+\frac{1}{2}\sum_{T\in\mathcal{T}_h}\int_T(\nabla\cdot \de)\bu\cdot \bP_h\de_t\,dT-\frac{1}{2}\sum_{e\in\Gamma_h}\int_e[\de]\cdot n_e\{\bu\cdot\bP_h\de_t\}\,ds\bigg|\nonumber\\
		\leq &\sum_{T\in\mathcal{T}_h}\|\de\|_{L^4(T)}\|\nabla\bu\|_{L^4(T)}\|\bP_h\de_t\|_{L^2(T)}+C\|\bu\|_{L^\infty(\Omega)}\sum_{T\in\mathcal{T}_h}\|\nabla\de\|_{L^2(T)}\|\bP_h\de_t\|_{L^2(T)}\nonumber\\
		&+C\|\bu\|_{L^\infty(\Omega)}\sum_{T\in\mathcal{T}_h}\frac{1}{|e|^{1/2}}\|[\de]\|_{L^2(e)}\|\bP_h\de_t\|_{L^2(T)}\leq  \frac{\kappa C_2}{12}\|\bP_h\de_t\|_\varepsilon^2+C\|\de\|_\varepsilon^2.\label{pr2}
	\end{align}
	Similarly, using the steps involved in obtaining estimates in (\ref{pr2}), we bound
	\begin{align}
		|c(\bu,\de,\bP_h\de_t)|\leq  \frac{\kappa C_2}{12}\|\bP_h\de_t\|_\varepsilon^2+C\|\de\|_\varepsilon^2.\label{pr3}
	\end{align}
	Using (\ref{pr1})-(\ref{pr3}) in (\ref{nw1}), we arrive at
	\begin{align}\label{nw2}
		|(c(\bu_{h},\bu_h,\bP_h\de_t)-c(\bu,\bu,\bP_h\de_t))|\leq   \frac{\kappa C_2}{4}\|\bP_h\de_t\|_\varepsilon^2+ C\|\e\|^2_\varepsilon.
	\end{align}
	\noindent
	 Applying the Cauchy-Schwarz and Young inequalities, definition of space $\bJ_h$, (\ref{2.1}) and (\ref{P_h}), we bound the second, third and sixth terms on the right hand side of (\ref{et5}) as
	\begin{align}
		\kappa\,|a(\bP_h\bu_t-\bu_t,\bP_h\de_t)+J_0(\bP_h\bu_t-\bu_t,\bP_h\de_t)| &\leq \frac{C_2\kappa}{12}\|\bP_h\de_t\|_\varepsilon^2 +Ch^{2}|\bu_t|_2^2,\\
		 \nu|a(\e,\bP_h\de_t)+J_0(\e,\bP_h\de_t)| & \leq \frac{C_2\kappa}{12}\|\bP_h\de_t\|_\varepsilon^2+C\|\e\|_\varepsilon^2+ Ch^2|\bu|_2^2,\\ 
		|b(\bP_h\de_t,p))|=|b(\bP_h\de_t,p-r_h(p)))| & \leq \frac{C_2\kappa}{12}\|\bP_h\de_t\|_\varepsilon^2+Ch^{2}|p|_1^2.\label{pr6}
	\end{align}
	An application of the bounds from (\ref{nw2})-(\ref{pr6}) in (\ref{et5}) leads to
	 \begin{align}
		\frac{1}{2}\left(\|\de_t\|^2+C_2\,\kappa\,\|\bP_h\de_t\|_\varepsilon^2\right)&\leq C\left( \|\bu_t-\bP_h\bu_t\|^2 +h^{2}(|\bu_t|_2^2+|\bu|_2^2+|p|_1^2)+\|\e\|_{\varepsilon}^2\right).\label{et6}
	\end{align}
	Finally, a use of triangle inequality, (\ref{P_h}), (\ref{et6}), Lemma \ref{exactpriori} and Theorem \ref{semierror} yield 
	\begin{align*}
		\|\de_t\|^2+K\,\kappa\,\|\de_t\|_\varepsilon^2\leq Ch^2e^{CT}.
	\end{align*}
	This completes the rest of the proof.
	
\end{proof}
\noindent
{\it Proof of Theorem \ref{semierrorp}: } A use of (\ref{8.5}), (\ref{8.6}), (\ref{8.8}) and (\ref{8.8n}) yield
\begin{align}
	(\bu_{ht}-\bu_t,\bv_h)+&\kappa\,(a(\bu_{ht}-\bu_t,\bv_h)+J_0(\bu_{ht}-\bu_t,\bv_h))
	+\nu(a(\bu_h-\bu,\bv_h)+J_0(\bu_h-\bu,\bv_h))\nonumber\\&+c^{\bu_h}(\bu_h,\bu_h,\bv_h)-c^{\bu}(\bu,\bu,\bv_h)-b(\bv_h,p-r_h(p))=-b(\bv_h,p_h-r_h(p)),\quad \forall \bv_h\in\bV_h.\label{30.16}
\end{align}
Using the inf-sup condition stated in Lemma \ref{inf-sup}, there exists $\bv_h\in \bV_h$ such that
\begin{equation}
	b(\bv_h,p_h-r_h(p))=-\|p_h-r_h(p)\|^2,\quad \|v_h\|_\varepsilon\leq \frac{1}{\beta_0}\|p_h-r_h(p)\|.\label{30.17}
\end{equation}
 A combination of (\ref{30.16}) and (\ref{30.17}) leads to
\begin{align}
	\|p_h-r_h(p)\|^2=&(\bu_{ht}-\bu_t,\bv_h)+\kappa\,(a(\bu_{ht}-\bu_t,\bv_h)+J_0(\bu_{ht}-\bu_t,\bv_h))\nonumber\\
	&+\nu(a(\bu_h-\bu,\bv_h)+J_0(\bu_h-\bu,\bv_h))\nonumber\\
	&+c(\bu_h,\bu_h,\bv_h)-c(\bu,\bu,\bv_h)-b(\bv_h,p-r_h(p)).\label{30.18}
\end{align}
 Since $\bu$ is continuous, the superscripts in nonlinear terms of (\ref{30.18}) are dropped. Now, following the analysis used in Lemmas \ref{etaerror} and \ref{etLinfL2},  we bound terms on the right hand side of (\ref{30.18}) as 
\begin{align*}
	\|p_h-r_h(p)\|^2\leq C\left(\|\bu_{ht}-\bu_t\|^2+\|\bu_{ht}-\bu_t\|_\varepsilon^2+\|\bu_h-\bu\|_\varepsilon^2+\|\bu_h-\bu\|^2+h^2(|\bu|_2^2+|p|_1^2+|\bu_t|_2^2)\right).
\end{align*}
Using the triangle inequality and (\ref{2.1}), we arrive at 
\begin{align}
	\|p-p_h\|^2\leq \left(\|\bu_{ht}-\bu_t\|^2+\|\bu_{ht}-\bu_t\|_\varepsilon^2+\|\bu_h-\bu\|_\varepsilon^2+\|\bu_h-\bu\|^2+h^2(|\bu|_2^2+|p|_1^2+|\bu_t|_2^2)\right).\label{n10}
\end{align}
An application of Lemmas \ref{exactpriori}, \ref{etLinfL2} and Theorem \ref{semierror} in (\ref{n10})  leads to the desired pressure error estimate. This completes the proof. \hfill{$\Box$}
\section{Fully discrete approximation and error estimates}\label{s6}
\se
In this section, we employ the backward Euler method for temporal discretization of the semidiscrete discontinuous Galerkin Kelvin-Voigt system represented by (\ref{8.8})-(\ref{8.9}).  Let $\Delta t$, $0<\Delta t<1$, denotes the time step and $0=t_0<t_1<\cdots<t_M=T$ be a subdivision of the interval $(0,T)$ with  $t_n=n\Delta t$, $n\geq 0$. Further, let the function $\psi(t)$ evaluated at time $t=t_n$ is denoted by $\psi^n$ with $ \bar\partial_t\psi^n$ is defined as
\[\bar\partial_t\psi^n=\frac{1}{\Delta t}(\psi^n-\psi^{n-1}), \,\,n\geq 0.\]
 Now, the backward Euler approximations for (\ref{8.8})-(\ref{8.9}) is defined as follows:
Given $\bU^0$, find $\{\bU^n\}_{n\geq 1}\in \bV_h$ and $\{P^n\}_{n\geq 1}\in M_h$, such that
\begin{align}
	(\bar \partial_t \bU^n,\bv_h&)+\kappa\,(a(\bar\partial_t\bU^n,\bv_h)+J_0(\bar\partial_t\bU^n,\bv_h))+\nu(a(\bU^n,\bv_h)+J_0(\bU^n,\bv_h))\nonumber\\
	&+c^{\bU^{n-1}}(\bU^{n-1},\bU^n,\bv_h) +b(\bv_h,P^n)=(\brf^n,\bv_h),\quad \forall \bv_h\in \bV_h,\label{fullyxh1}\\
	b(\bU^n,q_h)&=0,\quad \forall q_h\in M_h, \label{fullyxh2}
\end{align}
where $\bU^0=\bu_h(0)=\bP_h\bu_0$. \\
 An equivalent formulation of (\ref{fullyxh1})-(\ref{fullyxh2}) is defined as follows: For $\bv_h\in \bJ_h$, we seek $\{\bU^n\}_{n\geq 1}\in \bJ_h$, such that,
\begin{align}
	(\bar\partial_t \bU^n,\bv_h)+\kappa\,(a(\bar\partial_t\bU^n,\bv_h)&+J_0(\bar\partial_t\bU^n,\bv_h))+\nu(a(\bU^n,\bv_h)+J_0(\bU^n,\bv_h ))\nonumber\\
	&+c^{\bU^{n-1}}(\bU^{n-1},\bU^n,\bv_h) =(\brf^n,\bv_h), \label{fullysol}
\end{align}
where $\bU^0=\bu_h(0)=\bP_h\bu_0$. \\
Next in Lemma \ref{fullypriorisol}, we present {\it a priori} estimates of backward Euler solution $\bU^n$ of (\ref{fullysol}).
\begin{lemma}\label{fullypriorisol}
	Let the assumptions ({\bf A1})-({\bf A3}) be satisfied and let $0<\al<\frac{\nu C_2}{2(\gamma+\kappa C_1)}$. Then, the solution $\{\bU^n\}_{n\geq 1}$ of (\ref{fullysol}) satisfies the following estimates:
	\begin{align*}
		\|\bU^n\|^2+\|\bU^n\|^2_\varepsilon+e^{-2\al t_M}\Delta t\sum_{n=1}^{M}e^{2\al t_n}\|\bU^n\|^2_\varepsilon &\leq C,\quad n=0,1,...,M,
	\end{align*}
	where $C$ depends on the given data.
\end{lemma}
\begin{proof}
 Substitute $\bv_h=\bU^n$ in (\ref{fullysol}). Then, using
\begin{align*}
(\bar\partial _t\bU^n,\bU^n)&=\frac{1}{\Delta t}(\bU^n-\bU^{n-1},\bU^n)\geq \frac{1}{2\Delta t}(\|\bU^n\|^2-\|\bU^{n-1}\|^2)=\frac{1}{2}\bar\partial_t \|\bU^n\|^2,\\
a(\bar\partial_t\bU^n,\bU^n)&= \frac{1}{2}\bigg(\frac{1}{\Delta t}a(\bU^n,\bU^n)-\frac{1}{\Delta t}a(\bU^{n-1},\bU^{n-1})+\Delta t a(\bar\partial_t\bU^n,\bar\partial_t\bU^n)\bigg)\geq\frac{1}{2}\bar\partial_t a(\bU^n,\bU^n), 	
\end{align*}
 (\ref{5.5}) and Lemma \ref{coer} in the resulting equation and then the Cauchy-Schwarz inequality, we arrive at
\begin{align}
\bar\partial_t \|\bU^n\|^2+\kappa\,\bar\partial_t \left(a(\bU^n,\bU^n)+J_0(\bU^n,\bU^n)\right)+2\nu C_2\|\bU^n\|_\varepsilon^2\leq 2\|\brf^n\|\|\bU^n\|.\label{pri01}
\end{align}
 Note that,
\begin{align}\label{2r2}
 \sum_{n=1}^m \Delta t e^{2\al t_n}\bar\partial_t\|\bU^n\|^2=e^{2\al t_m}\|\bU^m\|^2-\sum_{n=1}^{m-1}  e^{2\al t_n}( e^{2\al \Delta t}-1)\|\bU^n\|^2-e^{2\al \Delta t}\|\bU^0\|^2.
 \end{align}
  Multiply (\ref{pri01}) by $\Delta t e^{2\al t_n}$, sum over $n=1$ to $m$, and use (\ref{Lp}), (\ref{2r2}) and Lemmas \ref{cont}, \ref{coer} to obtain
 \begin{align}\label{pri02}
 e^{2\al t_m}\|\bU^m\|^2+C_2\kappa\,e^{2\al t_m}\|\bU^m\|_\varepsilon^2+\left(\nu C_2-(\gamma+\kappa\,C_1)\frac{(e^{2\al \Delta t}-1)}{\Delta t}\right)\Delta t\sum_{n=1}^m e^{2\al t_n}\|\bU^n\|_\varepsilon^2\leq e^{2\al \Delta t}\|\bU^0\|^2\nonumber\\
 +\kappa\,C_1 e^{2\al \Delta t}\|\bU^0\|_\varepsilon^2+\frac{2\gamma^2}{\nu C_2}\Delta t\sum_{n=1}^m e^{2\al t_n}\|\brf\|^2.
 \end{align}
 With $\al$ as $0<\al<\frac{\nu C_2}{2(\gamma+\kappa C_1)}$ we have
 \[1+\frac{\nu C_2\Delta t}{\gamma+\kappa\,C_1}\geq e^{2\al \Delta t}.\]
\noindent
  On multiplying (\ref{pri02}) by $e^{-2\al t_m}$ and applying assumption ({\bf A1}), we establish our desired estimates.
 \end{proof}
\noindent
Below, we focus on the derivation of error estimates for the backward Euler method. \\
We denote $\de_n=\bU^n-\bu_h(t_n)=\bU^n-\bu_h^n$, $n\in\mathbb{N},\, 1<n\leq M$.
Now, consider the semidiscrete formulation (\ref{ne1}) at $t=t_n$ and subtract it from (\ref{fullysol}) to arrive at
	\begin{align}\label{11.0}
 (\bar\partial_t\de_n,\bphi_h)+\kappa(  a(\bar\partial_t\de_n,\bphi_h)+&J_0(\bar\partial_t\de_n,\bphi_h))+\nu(  a(\de_n,\bphi_h)+J_0(\de_n,\bphi_h))=(\bu_{ht}^n,\bphi_h)-( \bar\partial_t\bu_h^n,\bphi_h) \nonumber \\+\kappa\,(a(\bu_{ht}^n,\bphi_h)+&J_0\,(\bu_{ht}^n,\bphi_h))-\kappa\,\,(a (\bar\partial_t\bu_h^n,\bphi_h)+J_0 (\bar\partial_t\bu_h^n,\bphi_h))+c^{\bu_h^n}(\bu_h^n,\bu_h^n,\bphi_h) \nonumber \\&-c^{\bU^{n-1}}(\bU^{n-1},\bU^n,\bphi_h).
\end{align}
\begin{lemma}\label{fully}
	Under the assumptions of Theorem \ref{semierror}, there exists a positive constant $C=C(\kappa,\nu,\alpha,C_2, M_0)$, independent of $h$ and $\Delta t$, such that, the following estimates hold true:
	\begin{align*}
		(\|\de_n\|+C_2 \kappa\,\|\de_n\|_\varepsilon)+\left(\nu C_2 e^{-2\al t_M}\Delta t\sum_{n=1}^{M}e^{2\al t_n}\|\de_n\|^2_\varepsilon\right)^{1/2}\leq Ce^{CT}\Delta t.
	\end{align*}
\end{lemma}
\begin{proof}
	Choose $\bphi_h=\de_n$ in (\ref{11.0})  and use the facts 
	\[(\bar\partial _t\de_n,\de_n)=\frac{1}{\Delta t}(\de_n-\de_{n-1},\de_n)\geq \frac{1}{2\Delta t}(\|\de_n\|^2-\|\de_{n-1}\|^2)=\frac{1}{2}\bar\partial_t \|\de_n\|^2,\]
	\begin{align*}
			a(\bar\partial_t\de_n,\de_n)= \frac{1}{2}\bigg(\frac{1}{\Delta t}a(\de_n,\de_n)-\frac{1}{\Delta t}a(\de_{n-1},\de_{n-1})+\Delta t a(\bar\partial_t\de_n,\bar\partial_t\de_n)\bigg)\geq\frac{1}{2}\bar\partial_t a(\de_n,\de_n) 
	\end{align*}
	and Lemma \ref{coer} to arrive at
	\begin{align}
		\bar\partial_t(&\|\de_n\|^2+\kappa(  a(\de_n,\de_n)+J_0(\de_n,\de_n)))+2C_2\nu\|\de_n\|^2_\varepsilon\leq 2(\bu^n_{ht},\de_n)-2(\bar\partial_t\bu_h^n,\de_n)\nonumber\\
		&+2\kappa\,(a(\bu^n_{ht},\de_n)+J_0\,(\bu_{ht}^n,\de_n))-2\kappa\,(a(\bar\partial_t\bu_h^n,\de_n)+J_0 (\bar\partial_t\bu_h^n,\de_n))\nonumber\\ &+2c^{\bu_h^n}(\bu_h^n,\bu_h^n,\de_n)-2c^{\bU^{n-1}}(\bU^{n-1},\bU^n,\de_n).\label{11.2}
	\end{align}
	For the nonlinear terms in (\ref{11.2}), we note that
	\begin{align}
		c^{\bu_h^n}(\bu_h^n,\bu_h^n,\de_n)-c^{\bU^{n-1}}(\bU^{n-1},\bU^n,\de_n)=&c^{\bu_h^n}(\bu_h^n,\bu_h^n,\de_n)-c^{\bu_h^n}(\bU^{n-1},\bU^n,\de_n)\nonumber\\
		&+c^{\bu_h^n}(\bU^{n-1},\bU^n,\de_n)-c^{\bU^{n-1}}(\bU^{n-1},\bU^n,\de_n).\label{en1}
	\end{align}
	After dropping the superscripts for the first two nonlinear terms on the right hand side of (\ref{en1}), rewrite them as
	\begin{align}
		c(\bu_h^n, \bu_h^n,\de_n)-c(\bU^{n-1}  ,\bU^n,\de_n)= &-c(\bU^{n-1},\de_n,\de_n) +c(\de_{n-1},\bu^n-\bu_h^n,\de_n)-c(\de_{n-1},\bu^n,\de_n)\nonumber\\
		&+c(\bu_h^n-\bu_h^{n-1},\bu_h^n-\bu^n,\de_n)+c(\bu_h^n-\bu_h^{n-1},\bu^n,\de_n).\label{en2}
	\end{align}
	 The first term on the right hand side of (\ref{en2}) is positive due to (\ref{5.5}). The second term on the right hand side of  (\ref{en2}) can be bounded following the similar steps as involved in estimating $c(\be,\bxi,\be)$ in Lemma \ref{etaerror}.  A use of Cauchy-Schwarz's, Young's  inequalities and Theorem \ref{semierror} yield
	\begin{align}
		|c(\de_{n-1},\bu^n-\bu_h^n,\de_n)|\leq C\|\de_{n-1}\|\|\de_n\|_\varepsilon\leq \frac{C_2\nu}{8}\|\de_n\|^2_\varepsilon+C\|\de_{n-1}\|^2.\label{en3}
	\end{align}
	 Apply (\ref{tri1}), the Cauchy-Schwarz and Young inequalities and Lemmas \ref{sobolev1}, \ref{exactpriori} to bound the third term as 
	\begin{align}
		|c(\de_{n-1},\bu^n,\de_n)|\leq C\|\de_{n-1}\|\|\bu^n\|_2\|\de_n\|_\varepsilon\leq  \frac{C_2\nu}{8}\|\de_n\|^2_\varepsilon+C\|\de_{n-1}\|^2. \label{en4}
	\end{align}
	A use of $\bu_h^n-\bu_h^{n-1}=\displaystyle\int_{t_{n-1}}^{t_n}\bu_{ht}(t)\,dt$ and H{\"o}lder's, Young's inequalities, (\ref{Lp}), Lemmas \ref{trace}, \ref{distrace}, Theorem \ref{semierror} in the fourth term on the right hand side of (\ref{en2}) leads to
	\begin{align}
		|c(\bu_h^n-\bu_h^{n-1},\bu_h^n-\bu^n,\de_n)|&\leq  \int_{t_{n-1}}^{t_n}\sum_{T\in \mathcal{T}_h}\|\bu_{ht}\|_{L^4(T)}\|\nabla(\bu_h^n-\bu^n)\|_{L^2(T)}\|\de_n\|_{L^4(T)}\,dt\nonumber\\
		+\int_{t_{n-1}}^{t_n}\sum_{e\in\Gamma_h}\|\bu_{ht}\|_{L^4(e)}\|[\bu_h^n &-\bu^n]\|_{L^2(e)}\|\de_n\|_{L^4(e)}\,dt +C\int_{t_{n-1}}^{t_n}\sum_{T\in \mathcal{T}_h}\|\nabla\bu_{ht}\|_{L^2(T)}\|\bu_h^n-\bu^n\|_{L^4(T)}\|\de_n\|_{L^4(T)}\,dt\nonumber\\
		+C\int_{t_{n-1}}^{t_n} \sum_{e\in \Gamma_h}\|[\bu_{ht}]\|_{L^2(e)}&\|\bu_h^n-\bu^n\|_{L^4(e)}\|\de_n\|_{L^4(e)}\,dt \nonumber\\
		\leq   \frac{C_2\nu}{8}\|\de_n\|^2_\varepsilon +C\Delta t &\|\bu_{ht}\|^2_{L^2(t_{n-1},t_n;\varepsilon)}.
	\end{align}
	Using  the form of $c(\cdot,\cdot,\cdot)$ observed in (\ref{4.8}), (\ref{Lp}) and Lemmas \ref{sobolev1}, \ref{distrace} along with the fact that $\bu^n$ has zero jump, the last nonlinear term $c(\bu_h^n -\bu_h^{n-1},\bu^n,\de_n)$ on the right hand side of (\ref{en2}) can be bounded as
	\begin{align}
		|c (\bu_h^n &-\bu_h^{n-1} ,\bu^n,\de_n)|\nonumber\\
		\leq &\Delta t^{1/2}\|\bu_{ht}\|_{L^2(t_{n-1},t_n;L^2(\Omega))}\|\nabla\bu^n\|_{L^4(\Omega)}\|\de_n\|_{L^4(\Omega)}+C\Delta t^{1/2} \|\bu_{ht}\|_{L^2(t_{n-1},t_n;\varepsilon)}\|\bu^n\|_{L^\infty([0,T]\times\Omega)}\|\de_n\|\nonumber\\
		& +C\int_{t_{n-1}}^{t_n}\sum_{T\in \mathcal{T}_h}\|\bu^n\|_{L^\infty([0,T]\times\Omega)}|e|^{-1/2}\|[\bu_{ht}]\|_{L^2(e)}|e|^{1/2}h^{-1/2}\|\de_n\|_{L^2(T)}\,dt\nonumber\\
		\leq & \frac{C_2\nu}{8}\|\de_n\|^2_\varepsilon+C\Delta t\|\bu_{ht}\|^2_{L^2(t_{n-1},t_n;\varepsilon)}.
	\end{align}
	Using  a result in \cite[Proposition 4.10]{GR09},  we observe that 
	\begin{align*}
		|c^{\bu_h^n}(\bU^{n-1},\bU^n,\de_n)-c^{\bU^{n-1}}(\bU^{n-1},\bU^n,\de_n)|\leq C\|\bu_h^n-\bU^{n-1}\|_{L^4(\Omega)}\|\bU^n\|_\varepsilon\|\de_n\|_\varepsilon.
	\end{align*}
	An application of the triangle inequality, Lemmas \ref{sobolev2}, \ref{uhpriori}, \ref{fullypriorisol},  and Cauchy-Schwarz's inequality leads to
	\begingroup
	\allowdisplaybreaks
	\begin{align}
		| c^{\bu_h^n}&(\bU^{n-1},\bU^n,\de_n)-c^{\bU^{n-1}}(\bU^{n-1},\bU^n,\de_n)|\nonumber\\
		\leq &C\big(\|\bu_h^n-\bu_h^{n-1}\|_{L^4(\Omega)}\|\de_n\|_\varepsilon\|\de_n\|_{L^4(\Omega)} +\|\de_{n-1}\|_{L^4(\Omega)}\|\de_n\|_\varepsilon\|\de_n\|_{L^4(\Omega)}\nonumber\\
		&+ \|\bu_h^n-\bu_h^{n-1}\|_{L^4(\Omega)}\|\bu_h^n\|_\varepsilon\|\de_n\|_{L^4(\Omega)} +\|\de_{n-1}\|_{L^4(\Omega)}\|\bu_h^n\|_\varepsilon\|\de_n\|_{L^4(\Omega)}\big)\nonumber\\
		\leq &C\big(\|\bu_h^n-\bu_h^{n-1}\|^{1/2}\|\bu_h^n-\bu_h^{n-1}\|_\varepsilon^{1/2}\|\de_n\|_\varepsilon^{3/2}\|\de_n\|^{1/2} +\|\de_{n-1}\|^{1/2}\|\de_{n-1}\|_\varepsilon^{1/2}\|\de_n\|_\varepsilon^{3/2}\|\de_n\|^{1/2}\nonumber\\
		&+ \|\bu_h^n-\bu_h^{n-1}\|_\varepsilon\|\de_n\|_\varepsilon +\|\de_{n-1}\|^{1/2}\|\de_{n-1}\|_\varepsilon^{1/2}\|\de_n\|_\varepsilon\big)\nonumber\\
		\leq &   \frac{C_2\nu}{8}\|\de_n\|^2_\varepsilon+ \frac{C_2\nu}{8}\|\de_{n-1}\|^2_\varepsilon+C\|\de_{n-1}\|^2+C\|\de_{n-1}\|^2\|\de_{n-1}\|^2_\varepsilon+C\Delta t\|\bu_{ht}\|^2_{L^2(t_{n-1},t_n;\varepsilon)}.
	\end{align}
	Now to bound the first term on the right hand side of (\ref{11.2}), we use the Taylor series expansion and observe that 
	\begin{equation}\label{11.1}
		(\bu^n_{ht},\bphi_h)-(\bar\partial_t\bu_h^n,\bphi_h)=\frac{1}{\Delta t}\int_{t_{n-1}}^{t_n} (t-t_{n-1})(\bu_{htt},\bphi_h)\, dt.
	\end{equation}
	From (\ref{11.1}) and using Cauchy-Schwarz's inequality, we obtain 
	\begin{align}
		2(\bu^n_{ht},\de_n)-2(\bar\partial_t\bu_h^n,\de_n)& \leq C\Delta t^{1/2}\left(\int_{t_{n-1}}^{t_n}\|\bu_{htt}(t)\|^2_{-1,h}\,dt\right)^{1/2}\|\de_n\|_\varepsilon\nonumber\\
		&\leq   \frac{C_2\nu}{8}\|\de_n\|^2_\varepsilon+C\Delta t\int_{t_{n-1}}^{t_n}\|\bu_{htt}(t)\|^2_{-1,h}\,dt.\label{en91}
	\end{align}
	Similarly, following the steps involved in bounding (\ref{en91}), we arrive at
	\begin{align}
		2\kappa\,(a(\bu^n_{ht},\de_n)+J_0\,(\bu_{ht}^n,\de_n)-a(\bar\partial_t\bu_h^n,\de_n) -J_0 (\bar\partial_t\bu_h^n,\de_n))& \leq C\Delta t^{1/2}\left(\int_{t_{n-1}}^{t_n}\|\bu_{htt}(t)\|^2_\varepsilon\,dt\right)^{1/2}\|\de_n\|_\varepsilon\nonumber\\
		&\leq  \frac{C_2\nu}{8}\|\de_n\|^2_\varepsilon+C\Delta t\int_{t_{n-1}}^{t_n}\|\bu_{htt}(t)\|^2_\varepsilon\,dt.\label{en9}
	\end{align}
	\endgroup
	Apply (\ref{en1})-(\ref{en9}) in (\ref{11.2}), multiply the resulting inequality by $\Delta t e^{2 \al n \Delta t}$ and sum over $1\leq n\leq m\leq M$, where $T=M\Delta t$. Then, using  
	\begin{align*}
		\sum_{n=1}^m \Delta t e^{2\al n \Delta t}\bar\partial_t\|\de_n\|^2 &=e^{2\al m \Delta t}\|\de_m\|^2-\sum_{n=1}^{m-1}  e^{2\al n \Delta t}( e^{2\al \Delta t}-1)\|\de_n\|^2,\\
			\sum_{n=1}^m   \Delta t e^{2\al n \Delta t}\bar\partial_t(a+J_0)(\de_n,\de_n) &=e^{2\al m \Delta t}(a+J_0)(\de_m,\de_m)-\sum_{n=1}^{m-1}  e^{2\al n \Delta t}(e^{2\al \Delta t}-1)(a+J_0)(\de_n,\de_n),
	\end{align*}
	and Lemmas \ref{cont}, \ref{coer}, we arrive at
	\begin{align}
		e^{2\al m \Delta t}&(\|\de_m\|^2+ C_2 \kappa\,\|\de_m\|^2_\varepsilon)+ C_2\nu\Delta t\sum_{n=1}^me^{2\al n\Delta t}\|\de_n\|_\varepsilon^2\leq \sum_{n=1}^{m-1}  e^{2\al n \Delta t}( e^{2\al \Delta t}-1)(\|\de_n\|^2 +C_1\kappa\|\de_n\|_\varepsilon^2)\nonumber\\
		&+\Delta t\sum_{n=1}^m e^{2\al n \Delta t}(C+C\|\de_{n-1}\|_\varepsilon^2)\|\de_{n-1}\|^2+ C{\Delta t}^2\sum_{n=1}^m e^{2\al n \Delta t}\int_{t_{n-1}}^{t_n}\|\bu_{ht}(t)\|_\varepsilon^2\,dt\nonumber\\
		&+C{\Delta t}^2\sum_{n=1}^m e^{2\al n \Delta t}\int_{t_{n-1}}^{t_n}(\|\bu_{htt}(t)\|_{-1,h}^2+\|\bu_{htt}(t)\|^2_\varepsilon)\,dt.\label{11.5}
	\end{align}
	Using $e^{2\al \Delta t}-1\leq C\Delta t$, the first term on right hand side of (\ref{11.5}) can be merged with the second term on right hand side of (\ref{11.5}). Apply Lemma \ref{uhpriori} to observe that 
	\begin{align}
		\sum_{n=1}^m e^{2\al n \Delta t}\int_{t_{n-1}}^{t_n}\|\bu_{ht}(t)\|_\varepsilon^2\,dt=\sum_{n=1}^m\int_{t_{n-1}}^{t_n}e^{2\al (t_n-t)}e^{2\al t}\|\bu_{ht}(t)\|_\varepsilon^2\,dt&\leq e^{2\al \Delta t}\int_0^{t_m}e^{2\al t}\|\bu_{ht}(t)\|_\varepsilon^2\,dt\nonumber\\
		&\leq Ce^{2\al (m+1)\Delta t}.\label{nw11}
	\end{align}
	The last term on the right hand side of (\ref{11.5}) can be bounded in a similar manner as in (\ref{nw11}) by $C{\Delta t}^2e^{2\al (m+1)\Delta t}$.  An application of (\ref{nw11}) in (\ref{11.5}) leads to
	\begin{align*}
		e^{2\al m \Delta t}(\|\de_m\|^2+ C_2 \kappa\,\|\de_m\|^2_\varepsilon)+C_2\nu\Delta t\sum_{n=1}^me^{2\al n\Delta t}\|\de_n\|_\varepsilon^2\leq & C\Delta t\sum_{n=1}^m e^{2\al (n-1) \Delta t}(1+\|\de_{n-1}\|_\varepsilon^2)(\|\de_{n-1}\|^2+\|\de_{n-1}\|_\varepsilon^2)\\
		& +C{\Delta t}^2e^{2\al (m+1)\Delta t}.
	\end{align*}
	An application of discrete Gronwall's lemma completes the rest of the proof.
\end{proof}
\begin{theorem}\label{velofinal}
	Under the assumptions of the Theorem \ref{semierror} and Lemma \ref{fully}, the following holds true:
	\begin{align*}
		\|\bu(t_n)-\bU^n\|\leq Ce^{CT}(h^2+\Delta t),\\
		\|\bu(t_n)-\bU^n\|_\varepsilon \leq Ce^{CT}(h+\Delta t).
		\end{align*}
\end{theorem}
\begin{proof}
	A combination of Theorem \ref{semierror} and  Lemma \ref{fully} leads to the desired result.
\end{proof}
\begin{lemma}\label{neget}
Under the assumptions of Theorem \ref{semierror} and  Lemma \ref{fully}, there exists a positive constant $C=C(\kappa,\nu,\alpha,C_2, M_0)$, independent of $h$ and $\Delta t$, such that, the error $\de_n=\bU^n-\bu_h^n$, satisfies
	\begin{align*}
		\|\bar\partial_t\de_n\|+C_2\kappa \|\bar\partial_t\e_n\|_\varepsilon\leq C\Delta t.
	\end{align*}
\end{lemma}
\begin{proof}
	Rewrite the non-linear terms in (\ref{11.0})  as
	\begin{align}\label{neget2}
		c^{\bu_h^n}(\bu_h^n,\bu_h^n,\bphi_h)-c^{\bU^{n-1}}(\bU^{n-1},\bU^n,\bphi_h)=&-c^{\bu_h^n}(\bu_h^{n-1},\de_n,\bphi_h) +c^{\bu_h^n}(\bu_h^n-\bu_h^{n-1},\bu_h^n-\bu^n,\bphi_h)\nonumber\\
		&+c^{\bu_h^n}(\bu_h^n-\bu_h^{n-1},\bu^n,\bphi_h) -c^{\bu_h^n}(\de_{n-1},\de_n,\bphi_h)-c^{\bu_h^n}(\de_{n-1},\bu_h^n,\bphi_h)\nonumber\\
		& +c^{\bu_h^n}(\bU^{n-1},\bU^n,\bphi_h)-c^{\bU^{n-1}}(\bU^{n-1},\bU^n,\bphi_h).
	\end{align}
	 A use of (\ref{tri2}) yield
	\begin{align}\label{neget3}
&		|c^{\bu_h^n}(\bu_h^{n-1},\de_n,\bphi_h)+c^{\bu_h^n}(\de_{n-1},\de_n,\bphi_h)+c^{\bu_h^n}(\de_{n-1},\bu_h^n,\bphi_h)| \nonumber \\
		\leq ~&C\|\bu_h^{n-1}\|_\varepsilon\|\de_n\|_\varepsilon\|\bphi_h\|_\varepsilon+C\|\de_{n-1}\|_\varepsilon\|\de_n\|_\varepsilon\|\bphi_h\|_\varepsilon
		+C\|\de_{n-1}\|_\varepsilon\|\bu_h^n\|_\varepsilon\|\bphi_h\|_\varepsilon.
	\end{align}
	 Similar to Lemma \ref{fully}, using (\ref{tri1}), Theorem \ref{semierror} and Lemma \ref{exactpriori},  we bound the nonlinear terms as follows:
	\begin{align}\label{neget4}
		|c^{\bu_h^n}(\bu_h^n-\bu_h^{n-1},\bu_h^n-\bu^n,\bphi_h)+c^{\bu_h^n}(\bu_h^n-\bu_h^{n-1},\bu^n,\bphi_h)|\leq C\|\bu_h^n-\bu_h^{n-1}\|\|\bphi_h\|_\varepsilon.
	\end{align}
	 Further, using the similar steps as in Lemma \ref{fully} and the estimates of Theorem \ref{semierror}, we arrive at
	\begin{align}\label{neget5}
		|c^{\bu_h^n}(\bU^{n-1} &,\bU^n,\bphi_h)-c^{\bU^{n-1}}(\bU^{n-1},\bU^n,\bphi_h)|\nonumber\\
		\leq &\|\bu_h^n-\bU^{n-1}\|_{L^4(\Omega)}\|\bU^n-\bu^n\|_\varepsilon\|\bphi_h\|_{L^4(\Omega)}\nonumber\\
		\leq & \|\bu_h^n-\bU^{n-1}\|_{L^4(\Omega)}\|\de_n\|_\varepsilon\|\bphi_h\|_\varepsilon+ \|\bu_h^n-\bU^{n-1}\|_{L^4(\Omega)}\|\bu^n-\bu_h^n\|_\varepsilon\|\bphi_h\|_\varepsilon\nonumber\\
		\leq & \|\bu_h^n-\bu_h^{n-1}\|_\varepsilon\|\de_n\|_\varepsilon\|\bphi_h\|_\varepsilon+\|\de_{n-1}\|_\varepsilon\|\de_n\|_\varepsilon\|\bphi_h\|_\varepsilon+C\|\bu_h^n-\bu_h^{n-1}\|\|\bphi_h\|_\varepsilon+C\|\de_{n-1}\|\|\bphi_h\|_\varepsilon.
	\end{align}
	A use of Lemma \ref{cont} leads to
	\begin{align}\label{neget6}
		|a(\de_n,\bphi_h)+J_0(\de_n,\bphi_h)|\leq  C_1\|\de_n\|_\varepsilon\|\bphi_h\|_\varepsilon.
	\end{align}
	Apply (\ref{11.1}), Cauchy-Schwarz's and Young's inequalities  to arrive at
	\begin{align}
		(\bu^n_{ht},\bphi_h)-(\bar\partial_t\bu_h^n,\bphi_h)& \leq C\Delta t^{1/2}\left(\int_{t_{n-1}}^{t_n}\|\bu_{htt}(s)\|^2_{-1,h}\,ds\right)^{1/2}\|\bphi_h\|_\varepsilon\nonumber\\
		&\leq C 	\left(\sup_{0<t<\infty}\|\bu_{htt}(t)\|^2_{-1,h}\right)\left(\int_{t_{n-1}}^{t_n} 1  ds\right)^{1/2} \Delta t^{1/2} \|\bphi_h\|_\varepsilon.\label{neget7}
	\end{align}
	Following the similar analysis as in (\ref{neget7}), we obtain
	\begin{align}\label{neget71}
		\kappa\,(a(\bu^n_{ht},\bphi_h)+J_0\,(\bu_{ht}^n,\bphi_h)-a(\bar\partial_t\bu_h^n,\bphi_h) -J_0 (\bar\partial_t\bu_h^n, \bphi_h))&\leq C 	\left(\sup_{0<t<\infty}\|\bu_{htt}(t)\|^2_\varepsilon\right)\left(\int_{t_{n-1}}^{t_n} 1  ds\right)^{1/2}\times\nonumber\\
		& \Delta t^{1/2} \|\bphi_h\|_\varepsilon.
	\end{align}
	Substitute $\bphi_h=\bar\partial_t\de_n$ in (\ref{11.0}) and  use (\ref{neget2})--(\ref{neget71}) with $\bphi_h$ replaced by $\bar\partial_t\de_n$. Then, apply Young's inequality and estimates from Lemmas \ref{coer}, \ref{uhpriori}, \ref{fully} to arrive at the desired result. This completes the rest of the proof.
\end{proof}
\begin{lemma}\label{fullypressure}
	Under the assumptions of Theorem \ref{semierror} and  Lemma \ref{fully}, there exists a positive constant $C=C(\kappa,\nu,\alpha,C_2, M_0)$, independent of $h$ and $\Delta t$, such that for $n=1,2,\cdots,M$
	\begin{align*}
		\|P^n-p_h^n\|\leq C\Delta t.
	\end{align*}
\end{lemma}
\begin{proof}
	Consider (\ref{8.8}) at $t=t_n$ and subtract it from (\ref{fullyxh1}) to arrive at
	\begin{align}
		b(\bphi_h,P^n-p_h^n)= &(\bu_{ht}^n,\bphi_h) -(\bar\partial_t\bu_h^n,\bphi_h)+\kappa\,(a(\bu_{ht}^n,\bphi_h)+J_0(\bu_{ht}^n,\bphi_h)) -\kappa\,(a(\bar\partial_t\bu_h^n,\bphi_h)+J_0(\bar\partial_t\bu_h^n,\bphi_h))\nonumber\\
		&-(\bar\partial_t\de_n,\bphi_h)-\kappa\,( a(\bar\partial_t\de_n,\bphi_h)+J_0(\bar\partial_t\de_n,\bphi_h))
		-\nu(  a(\de_n,\bphi_h)+J_0(\de_n,\bphi_h))\nonumber\\
		&+c^{\bu_h^n}(\bu_h^n,\bu_h^n,\bphi_h)-c^{\bU^{n-1}}(\bU^{n-1},\bU^n,\bphi_h),\quad \forall \bphi_h\in \bV_h.\label{fullyerrpn1}
	\end{align}
	 Following the steps involved in the proof of Lemma \ref{neget} and applying Lemma \ref{inf-sup}, we obtain  
	\begin{align}
		\|P^n-p_h^n\|\leq C (\|\partial_t\de_n\|+\kappa\|\partial_t\de_n\|_{\varepsilon}+\|\de_n\|_\varepsilon+\Delta t).\label{fullyerrpn2}
	\end{align}
	Finally, an application of the  Lemmas \ref{fully} and \ref{neget} concludes the proof.
\end{proof}
\noindent
A combination of Lemma \ref{fullypressure} and Theorem \ref{semierrorp} lead to the following fully discrete pressure error estimates.
\begin{theorem}\label{pressfinal}
	Under the assumptions of Theorem \ref{semierrorp} and Lemma \ref{fullypressure}, the following hold true:
	\begin{align*}
		\|p(t_n)-P^n\|\leq C(h+\Delta t).
	\end{align*}
\end{theorem}
\section{Numerical experiments}\label{s7}
\begin{figure}[h!]
	\centering
	\includegraphics[scale=0.5]{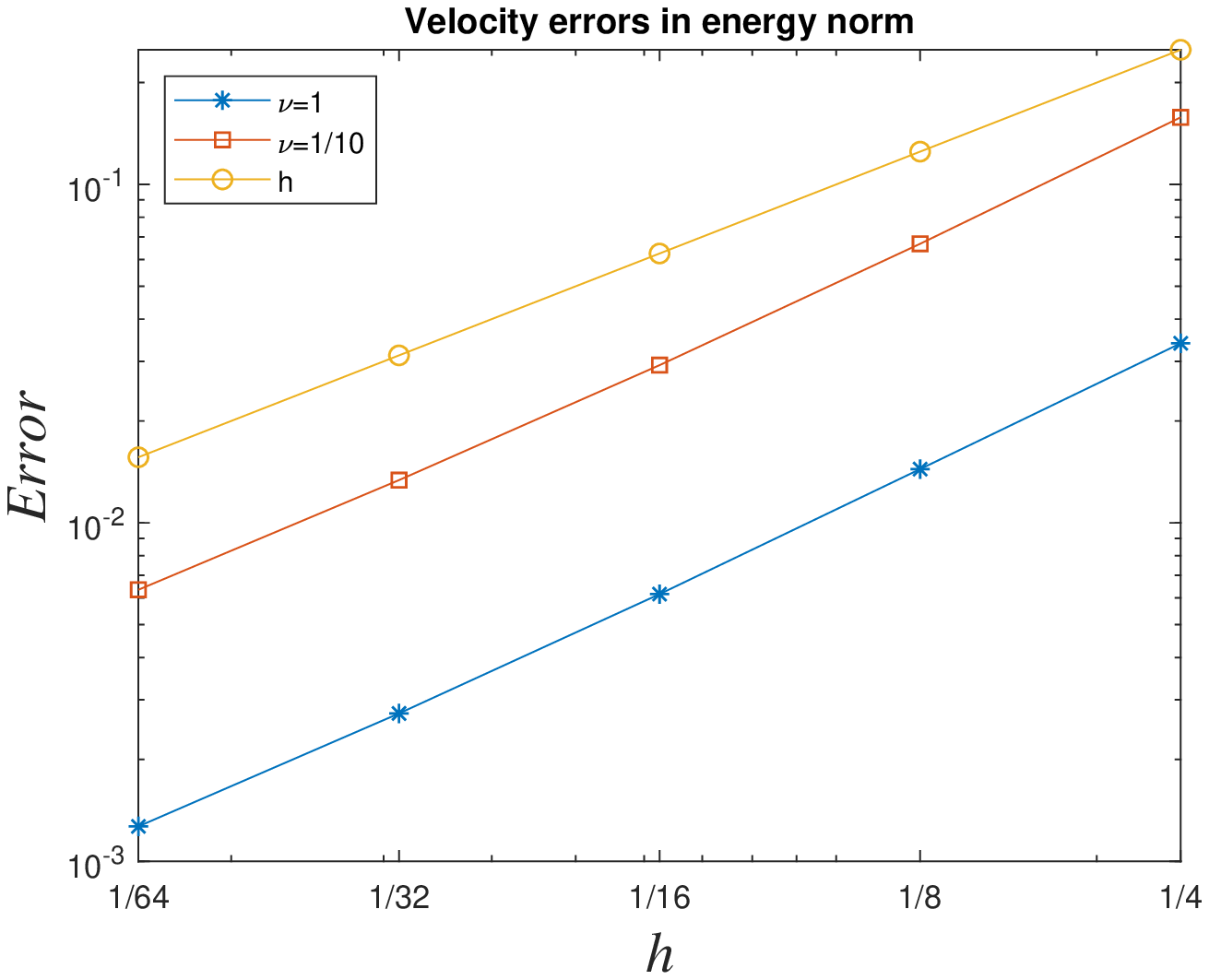}
	\includegraphics[scale=0.5]{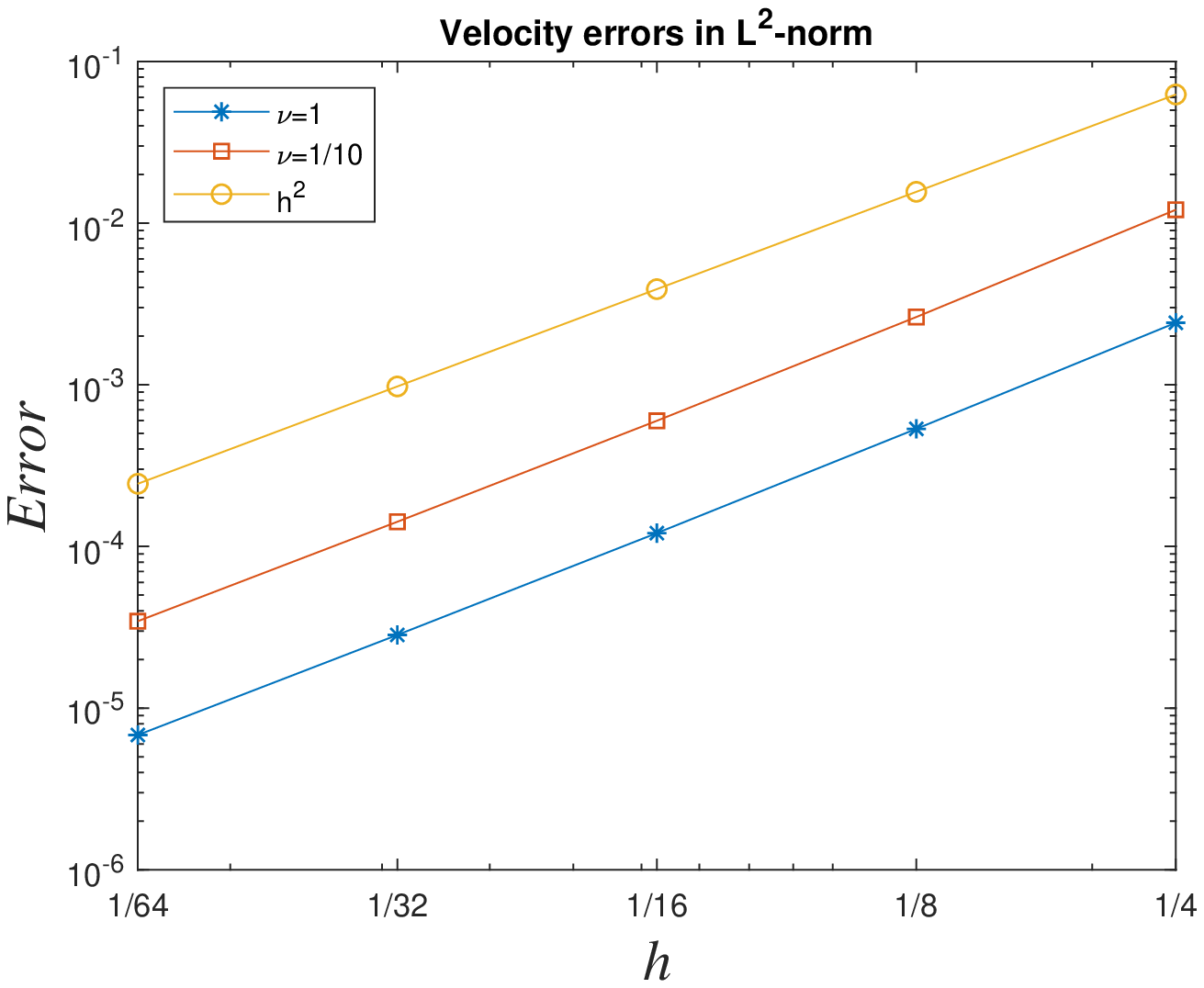}
	\includegraphics[scale=0.5]{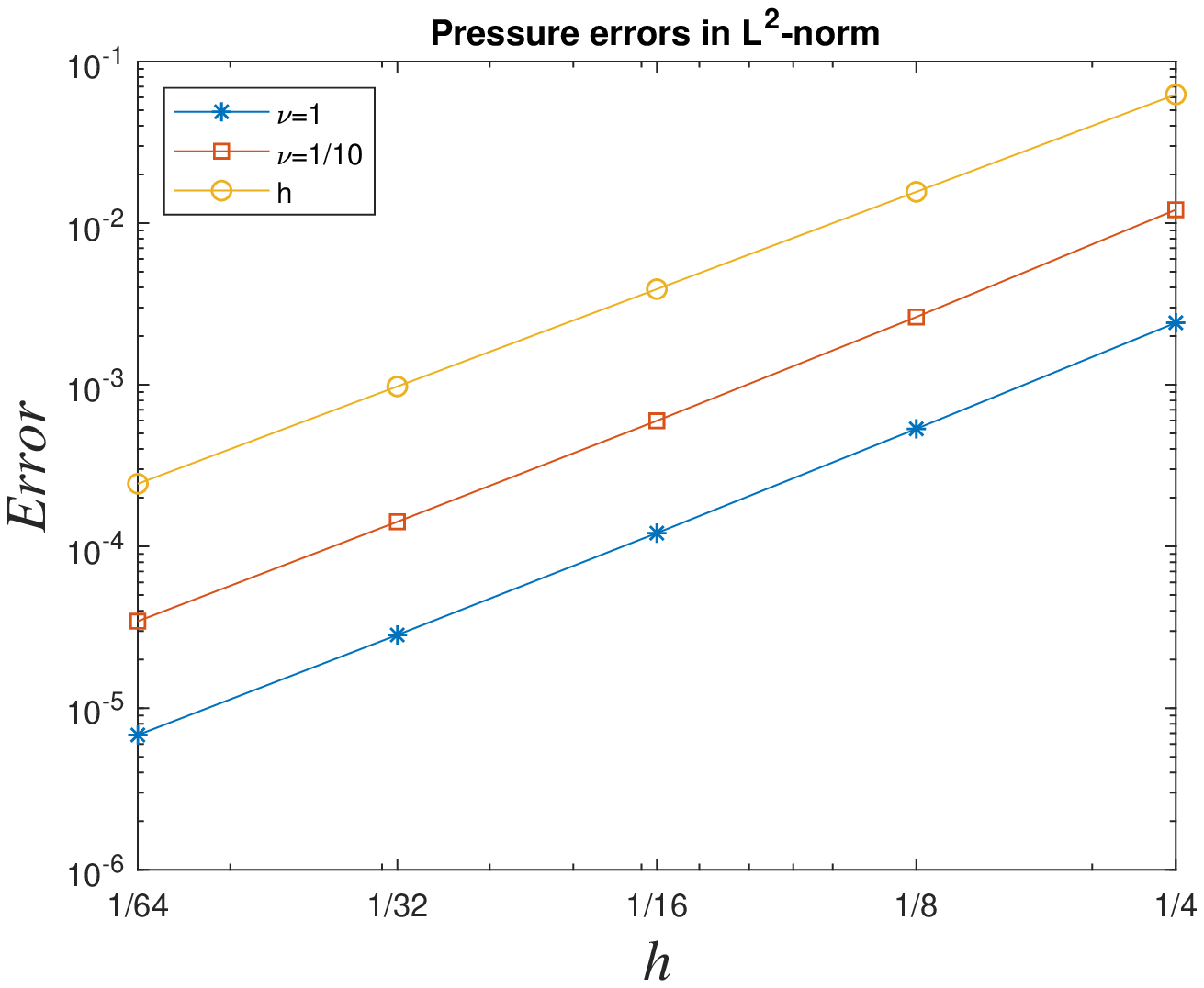}
	\caption{Velocity and pressure errors for $P_1-P_0$  with $\nu=1$ and $1/10$ for Example \ref{ex1}.}
	\label{fig1}
\end{figure}

\se
In this section, we present a couple of numerical experiments to verify the theoretical results stated in the  Theorems \ref{velofinal} and \ref{pressfinal}. The domain $\Omega$ is considered as $[0, 1]\times[0, 1]$.
We use the mixed finite element spaces $P_1-P_0$ and $P_1-P_1$ for the space discretization and a first order accurate backward Euler method for the time discretization. The time interval is chosen as $[0,1]$ with final time $T=1$. 
\begin{example}\label{ex1}
	In our first example, the right hand side function $\brf$ is chosen in such a way that the exact solution is
	\[\bu = (2 x^2 (x-1)^2 y (y-1) (2y-1)\cos t,-2x(x-1)(2x-1)y^2(y-1)^2 \cos t),\quad p=2(x-y)\cos t.\]
\end{example}
Tables 1 and 2 show the errors and the convergence rates for the mixed finite element space $P_1-P_0$  with  kinematic viscosities $\nu = 1$ and $1/10$, respectively. And Figure \ref{fig1} represents the velocity and pressure errors for $P_1-P_0$ element with $\nu=1$ and $1/10$.
In Table 3, we choose $P_1-P_1$ mixed finite element space with $\nu = 1$.  Note that, we choose a constant penalty parameter $\sigma_e = 10$ and retardation $\kappa=10^{-2}$ for Tables 1-3.  In Table 4, we have represented the errors and the convergence rates for the backward Euler method applied to continuous Galerkin finite element method with $\nu=1$ and $\kappa=10^{-2}$.  It can be observed that the numerical results represented in Tables 1-3  and Figure \ref{fig1} validate our theoretical findings in Theorems \ref{velofinal} and \ref{pressfinal}.
 Further, the results in Tables 3 and 4 represent that the discontinuous Galerkin finite element method works well for equal order element $P_1-P_1$, whereas continuous Galerkin finite element method fails to approximate the exact solution \cite{GR79}.
{\small
	\renewcommand{\arraystretch}{1.3}
	\begin{table}[ht!]
		\centering
		\caption{Numerical errors and convergence rates, for $P_1$--$P_0$ ($\nu$=1, $\kappa=10^{-2}$, $\sigma_e=10$) for Example \ref{ex1} with $\Delta t=\mathcal{O}(h^2)$.}
		\vspace{.4cm}
		\begin{tabular}{|c|c c| c c| c c|}
			
			\hline
			
			$h$ & $\|\bu(T)-\bU^M\|_\varepsilon$ & Rate & $\|\bu(T)-\bU^M\|_{L^2(\Omega)}$ & Rate & $\|p(T)-P^M\|_{L^2(\Omega)}$ & Rate\\
			
			\hline  
			1/4  & $3.3927\times 10^{-2}$  &  		  & $2.4193\times 10^{-3}$ & 		 & $1.0667\times 10^{-2}$ & 		\\
			1/8  & $1.4413\times 10^{-2}$  &  1.2350  & $5.3378\times 10^{-4}$ & 2.1802 & $7.0680\times 10^{-3}$ & 0.5937\\
			1/16 & $6.1566\times 10^{-3}$  &  1.2271  & $1.2089\times 10^{-4}$ & 2.1424 & $4.1982\times 10^{-3}$ & 0.7515 \\
			1/32 & $2.7335\times 10^{-3}$  &  1.1713  & $2.8341\times 10^{-5}$ & 2.0928 & $2.3113\times 10^{-3}$ & 0.8610 \\
			1/64 & $1.2676\times 10^{-3}$  &  1.1086  & $6.8138\times 10^{-6}$ & 2.0563 & $1.2137\times 10^{-3}$ & 0.9292 \\
			\hline
			
		\end{tabular}
		\vspace{.1cm}
		
	\end{table}
}

{\small
	\renewcommand{\arraystretch}{1.3}
	\begin{table}[ht!]
		\centering
		\caption{Numerical errors and convergence rates, for $P_1$--$P_0$ ($\nu=1/10$, $\kappa=10^{-2}$, $\sigma_e$=10) for Example \ref{ex1} with $\Delta t=\mathcal{O}(h^2)$.}
		\vspace{.4cm}
		\begin{tabular}{|c|c c| c c| c c|}
			
			\hline
			
			$h$ & $\|\bu(T)-\bU^M\|_\varepsilon$ & Rate & $\|\bu(T)-\bU^M\|_{L^2(\Omega)}$ & Rate & $\|p(T)-P^M\|_{L^2(\Omega)}$ & Rate\\
			
			\hline 
			1/4  & $1.5785\times 10^{-1}$  &  	     & $1.2098\times 10^{-2}$ & 	   & $1.0790\times 10^{-2}$ &	 \\
			1/8  & $6.6764\times 10^{-2}$  & 1.2414  & $2.6259\times 10^{-3}$ & 2.2038 & $7.4719\times 10^{-3}$ & 0.5302\\
			1/16 & $2.9235\times 10^{-2}$  & 1.1913  & $5.9873\times 10^{-4}$ & 2.1328 & $4.3042\times 10^{-3}$ & 0.7957\\
			1/32 & $1.3373\times 10^{-2}$  & 1.1283  & $1.4209\times 10^{-4}$ & 2.0750 & $2.3030\times 10^{-3}$ & 0.9022 \\
			1/64 & $6.3422\times 10^{-3}$  & 1.0762  & $3.4480\times 10^{-5}$ & 2.0430 & $1.1906\times 10^{-3}$ & 0.9518 \\
			\hline
			
		\end{tabular}
		\vspace{.1cm}
		
	\end{table}
}

{\small
	\renewcommand{\arraystretch}{1.3}
	\begin{table}[ht!]
		\centering
		\caption{Numerical errors and convergence rates, for $P_1$--$P_1$ ($\nu=1$, $\kappa=10^{-2}$, $\sigma_e$=10) for Example \ref{ex1} with $\Delta t=\mathcal{O}(h^2)$.}
		\vspace{.4cm}
		\begin{tabular}{|c|c c| c c| c c|}
			
			\hline
			
			$h$ & $\|\bu(T)-\bU^M\|_\varepsilon$ & Rate & $\|\bu(T)-\bU^M\|_{L^2(\Omega)}$ & Rate & $\|p(T)-P^M\|_{L^2(\Omega)}$ & Rate\\
			
			\hline   
			1/4  & $1.0869\times 10^{-2}$ &  		 & $7.9834\times 10^{-4}$ & 	   & $2.6782\times 10^{-2}$ &		\\
			1/8  & $6.1457\times 10^{-3}$ &  0.8226  & $4.4474\times 10^{-4}$ & 0.8440 & $2.1336\times 10^{-2}$ & 0.3279\\
			1/16 & $2.6823\times 10^{-3}$ &  1.1961  & $1.5730\times 10^{-4}$ & 1.4993 & $1.2302\times 10^{-2}$ & 0.7943\\
			1/32 & $1.1742\times 10^{-3}$ &  1.1916  & $4.4715\times 10^{-5}$ & 1.8147 & $6.5318\times 10^{-3}$ & 0.9133\\
			1/64 & $5.5135\times 10^{-4}$ &  1.0907  & $1.1724\times 10^{-5}$ & 1.9312 & $3.3487\times 10^{-3}$ & 0.9638\\
			\hline
			
		\end{tabular}
		\vspace{.1cm}
		
	\end{table}
}

{\small
	\renewcommand{\arraystretch}{1.3}
	\begin{table}[ht!]
		\centering
		\caption{Numerical errors and convergence rates for continuous finite element method using $P_1$--$P_1$ ($\nu=1$,  $\kappa=10^{-2}$) for Example \ref{ex1} with $\Delta t=\mathcal{O}(h^2)$.}
		\vspace{.4cm}
		\begin{tabular}{|c|c c| c c| c c|}
			
			\hline
			
			$h$ & $\|\bu(T)-\bU^M\|_{H^1(\Omega)}$ & Rate & $\|\bu(T)-\bU^M\|_{L^2(\Omega)}$ & Rate & $\|p(T)-P^M\|_{L^2(\Omega)}$ & Rate\\
			
			\hline   
			1/4  & $1.8439\times 10^{-2}$  & 		& $1.9654\times 10^{-3}$ & 		  & $0.20708$ &	\\
			1/8  & $8.4605\times 10^{-3}$  & 1.1240 & $6.1972\times 10^{-4}$ & 1.6651 & $0.11661$ & 0.8284\\
			1/16 & $3.8693\times 10^{-3}$  & 1.1286 & $1.5828\times 10^{-4}$ & 1.9690 & $0.07076$ & 0.7205\\
			1/32 & $1.8594\times 10^{-3}$  & 1.0572 & $3.9421\times 10^{-5}$ & 2.0055 & $0.05129$ & 0.4642 \\
			1/64 & $9.1443\times 10^{-4}$  & 1.0238 & $9.8055\times 10^{-6}$ & 2.0072 & $0.04443$ & 0.2072\\
			\hline
			
		\end{tabular}
		\vspace{.1cm}
		
	\end{table}
}

\begin{figure}[h!]
	\centering
	\includegraphics[scale=0.5]{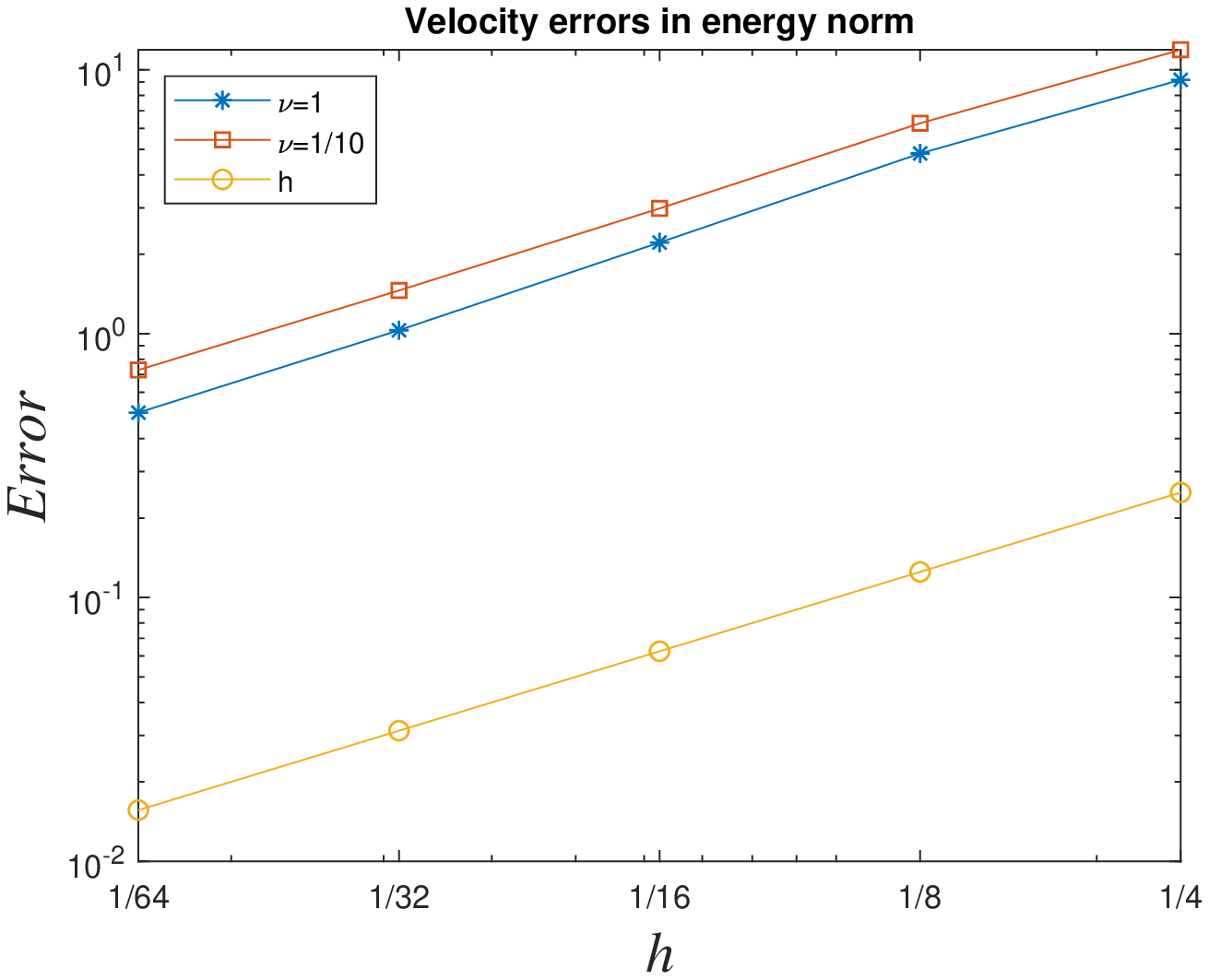}
	\includegraphics[scale=0.5]{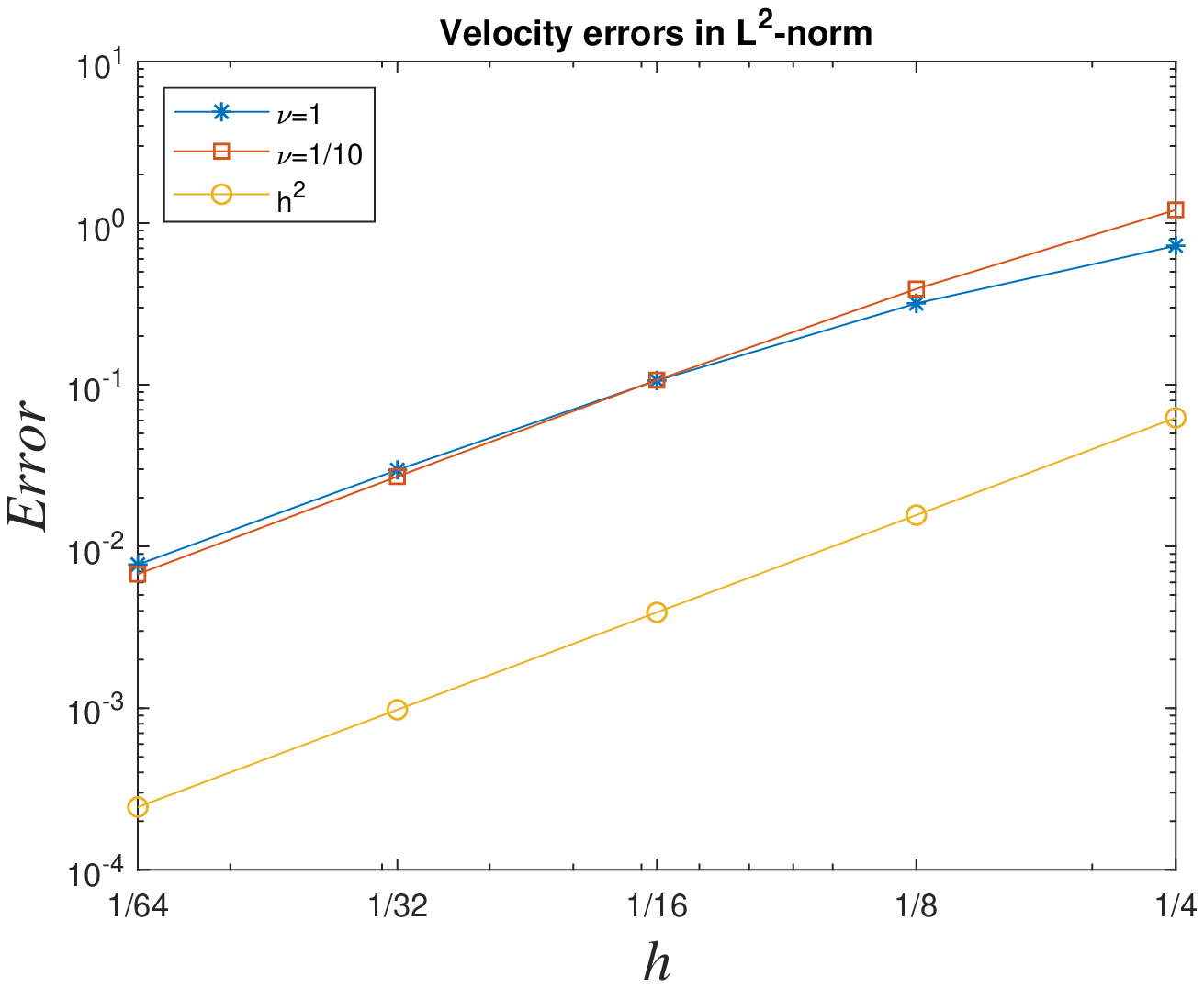}
	\includegraphics[scale=0.5]{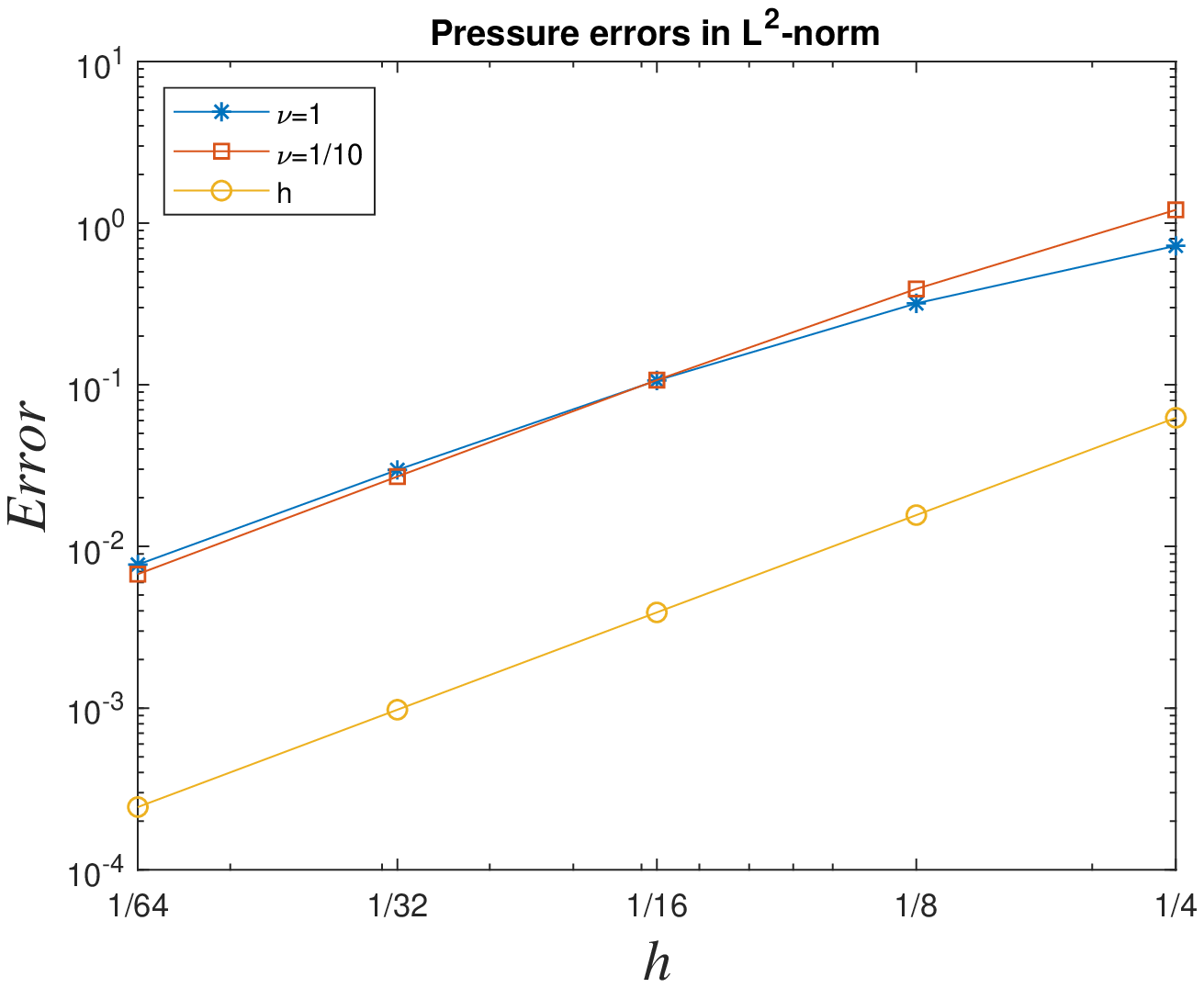}
	\caption{Velocity and pressure errors for $P_1-P_0$  with $\nu=1$ and $1/10$ for Example \ref{ex2}.}
	\label{fig2}
\end{figure}

\begin{example}\label{ex2}
	In this example, we take the forcing term  $\brf$ such that the solution of the problem to be
	\[\bu = (e^t(-cos(2\pi x)sin(2\pi y) + sin(2\pi y)),e^t(sin(2\pi x)cos(2\pi y) - sin(2\pi x)),\quad p=e^t (2\pi (cos(2\pi y)- cos(2\pi x))).\]
\end{example}
\noindent
In Tables 5 and 6, we have shown the errors and convergence rates for the mixed finite element space $P_1-P_0$  with  viscosities $\nu = 1$ and $ 1/10$, respectively.  And Figure \ref{fig2} represents the velocity and pressure errors for $P_1-P_0$ element with $\nu=1$ and $1/10$. 
For the cases $\nu = 1$ and $\nu = 1/10$, we have choosen  $\kappa=10^{-2}$ and  $\kappa=10^{-3}$, respectively. Table 7 depicts the results for the mixed space $P_1-P_1$ with $\nu = 1$ and $\kappa=10^{-2}$. We choose a constant penalty parameter $\sigma_e = 10$  for Tables 5-7.
 For continuous finite element method, the errors and the convergence rates for the case $P_1-P_1$ with $\nu=1$ and $\kappa=10^{-2}$  have shown in Table 8 for Example \ref{ex2}.  Here again, the numerical outcomes verify the derived theoretical results. Also, it can be concluded that the discontinuous Galerkin finite element method works well for equal order elements in comparison to the continuous Galerkin method. 

{\small
	\renewcommand{\arraystretch}{1.3}
	\begin{table}[ht!]
		\centering
		\caption{Numerical errors and convergence rates, for $P_1$--$P_0$ ($\nu=1$, $\kappa=10^{-2}$, $\sigma_e=10$) for Example \ref{ex2} with $\Delta t=\mathcal{O}(h^2)$.}
		\vspace{.4cm}
		\begin{tabular}{|c|c c| c c| c c|}
			
			\hline

			$h$ & $\|\bu(T)-\bU^M\|_\varepsilon$ & Rate & $\|\bu(T)-\bU^M\|_{L^2(\Omega)}$ & Rate & $\|p(T)-P^M\|_{L^2(\Omega)}$ & Rate\\
			
			\hline   
			1/4  & $9.17093$  &  	    & $0.72430$ & 		 & $9.37006$ &	     \\
			1/8  & $4.82057$  &  0.9278 & $0.31886$ & 1.1836 & $6.19123$ & 0.5978\\
			1/16 & $2.21524$  &  1.1217 & $0.10630$ & 1.5846 & $4.15326$ & 0.5759\\
			1/32 & $1.03043$  &  1.1042 & $0.02973$ & 1.8381 & $2.38552$ & 0.7999 \\
			1/64 & $0.50189$  &  1.0377 & $0.00771$ & 1.9454 & $1.25309$ & 0.9288\\
			\hline
		\end{tabular}
		\vspace{.1cm}
		
	\end{table}
}

{\small
	\renewcommand{\arraystretch}{1.3}
	\begin{table}[ht!]
		\centering
		\caption{Numerical errors and convergence rates, for $P_1$--$P_0$ ($\nu=1/10$, $\kappa=10^{-3}$, $\sigma_e=10$) for Example \ref{ex2} with $\Delta t=\mathcal{O}(h^2)$.}
		\vspace{.4cm}
		\begin{tabular}{|c|c c| c c| c c|}
			
			\hline
			
			$h$ & $\|\bu(T)-\bU^M\|_\varepsilon$ & Rate & $\|\bu(T)-\bU^M\|_{L^2(\Omega)}$ & Rate & $\|p(T)-P^M\|_{L^2(\Omega)}$ & Rate\\
			
			\hline   
			1/4  & $11.91970$ &        & $1.20830$ &        & $5.89632$ &        \\
			1/8  & $6.27597$  & 0.9254 & $0.39185$ & 1.6246 & $2.34282$ & 1.3315\\
			1/16 & $2.98342$  & 1.0728 & $0.10691$ & 1.8738 & $0.83518$ & 1.4880\\
			1/32 & $1.45801$  & 1.0329 & $0.02704$ & 1.9831 & $0.33008$ & 1.3392\\
			1/64 & $0.72815$  & 1.0016 & $0.00673$ & 2.0045 & $0.14996$ & 1.1381\\
			\hline
		\end{tabular}
		\vspace{.1cm}
		
	\end{table}
}

{\small
	\renewcommand{\arraystretch}{1.3}
	\begin{table}[ht!]
		\centering
		\caption{Numerical errors and convergence rates, for $P_1$--$P_1$ ($\nu=1$, $\kappa=10^{-2}$, $\sigma_e=10$) for Example \ref{ex2} with $\Delta t=\mathcal{O}(h^2)$.}
		\vspace{.4cm}
		\begin{tabular}{|c|c c| c c| c c|}
			
			\hline
			
			$h$ & $\|\bu(T)-\bU^M\|_\varepsilon$ & Rate & $\|\bu(T)-\bU^M\|_{L^2(\Omega)}$ & Rate & $\|p(T)-P^M\|_{L^2(\Omega)}$ & Rate\\
			
			\hline   
			1/4  & $9.68528$  &  	    & $0.84240$ & 		 & $20.72070$ &   	\\
			1/8  & $5.42842$  &  0.8352 & $0.45434$ & 0.8907 & $18.31680$ & 0.1779\\
			1/16 & $2.26422$  &  1.2615 & $0.15667$ & 1.5360 & $10.42180$ & 0.8135\\
			1/32 & $0.94989$  &  1.2531 & $0.04369$ & 1.8421 & $5.44668$  & 0.9361 \\
			1/64 & $0.43939$  &  1.1122 & $0.01130$ & 1.9500 & $2.76205$  & 0.9796\\
			\hline
		\end{tabular}
		\vspace{.1cm}
		
	\end{table}
}

{\small
	\renewcommand{\arraystretch}{1.3}
	\begin{table}[ht!]
		\centering
		\caption{Numerical errors and convergence rates, for continuous finite element method using $P_1$--$P_1$ ($\nu=1$, $\kappa=10^{-2}$) for Example \ref{ex2} with $\Delta t=\mathcal{O}(h^2)$.}
		\vspace{.4cm}
		\begin{tabular}{|c|c c| c c| c c|}
			
			\hline
			
			$h$ & $\|\bu(T)-\bU^M\|_{H^1(\Omega)}$ & Rate & $\|\bu(T)-\bU^M\|_{L^2(\Omega)}$ & Rate & $\|p(T)-P^M\|_{L^2(\Omega)}$ & Rate\\
			
			\hline   
			1/4  & $14.94269$ &  		& $1.60828$ &  		 & $53.92851$ &		\\
			1/8  & $4.97675$  &  1.5861 & $0.46035$ & 1.8047 & $28.08640$ & 0.9411\\
			1/16 & $2.01532$  &  1.3041 & $0.12058$ & 1.9327 & $23.25936$ & 0.2720\\
			1/32 & $0.92299$  &  1.1266 & $0.03039$ & 1.9879 & $21.65743$ & 0.1029 \\
			1/64 & $0.44765$  &  1.0439 & $0.00759$ & 2.0000 & $21.00652$ & 0.0440\\
			\hline
		\end{tabular}
		\vspace{.1cm}
	\end{table}
}

\pagebreak

\begin{figure}[h!]
	\centering
	\includegraphics[scale=0.5]{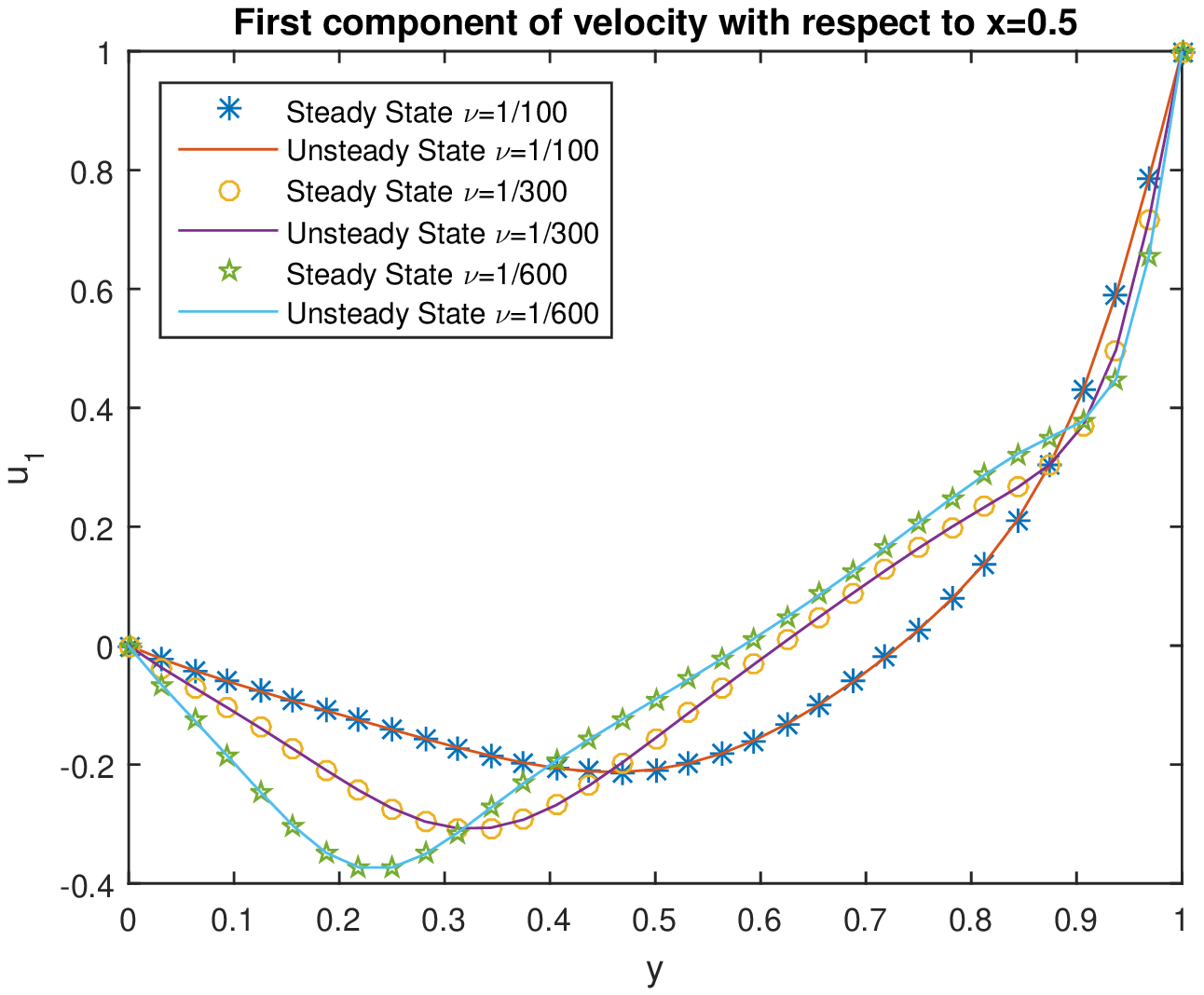}
	\includegraphics[scale=0.5]{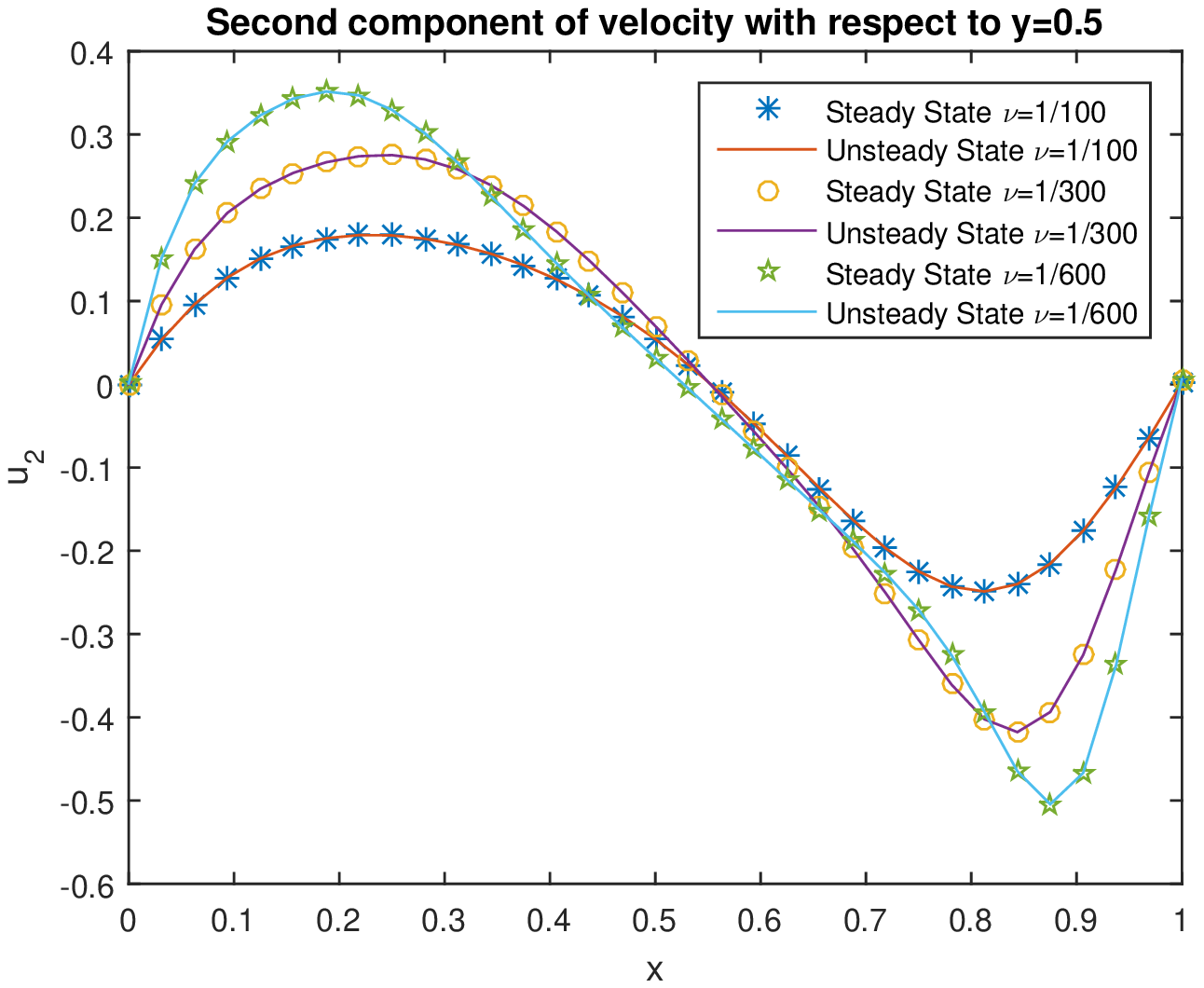}
	\caption{Velocity components for Example \ref{ex3}.}
	\label{fig3}
\end{figure}

\begin{figure}[h!]
	\centering
	\includegraphics[scale=0.5]{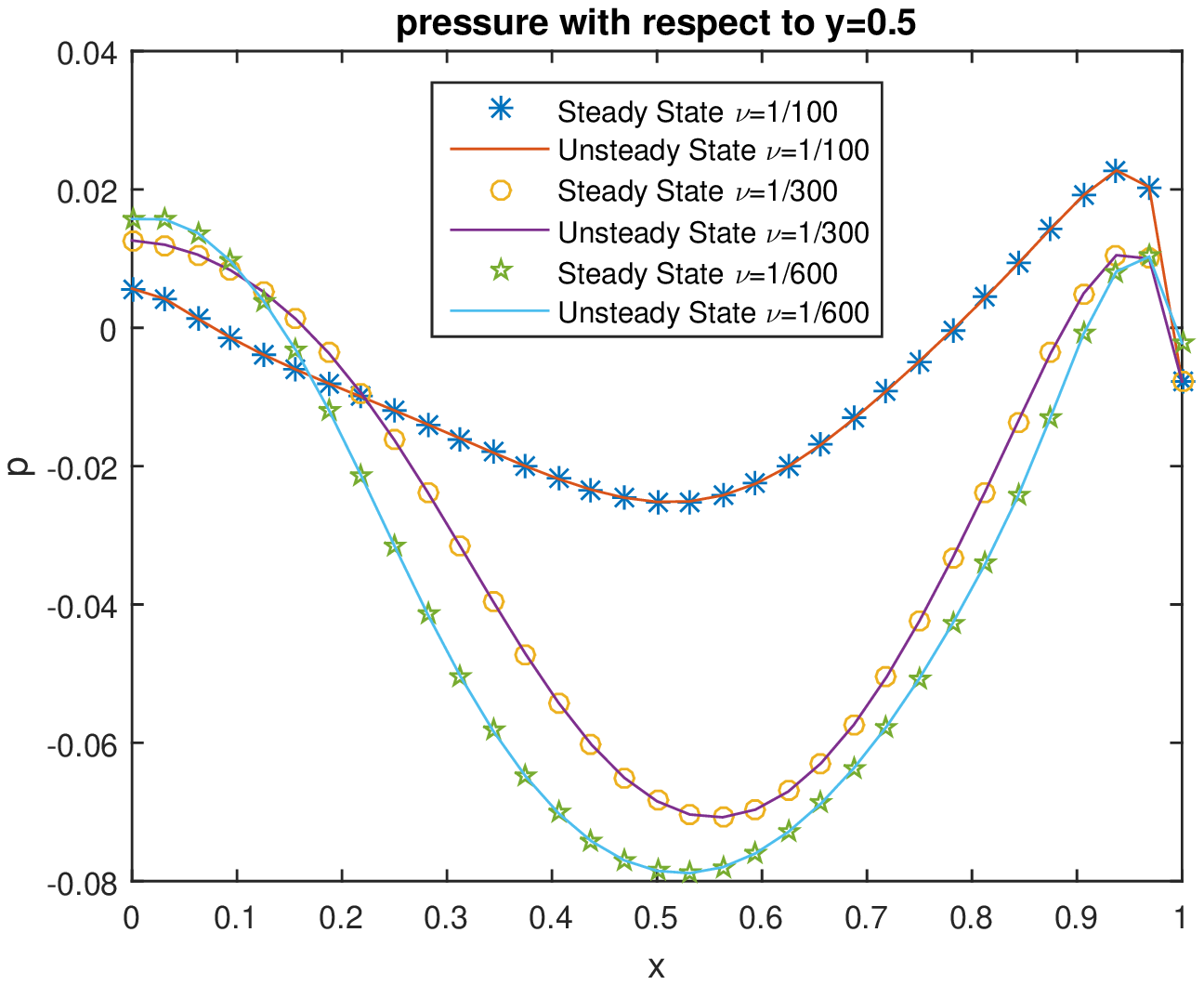}
	\includegraphics[scale=0.5]{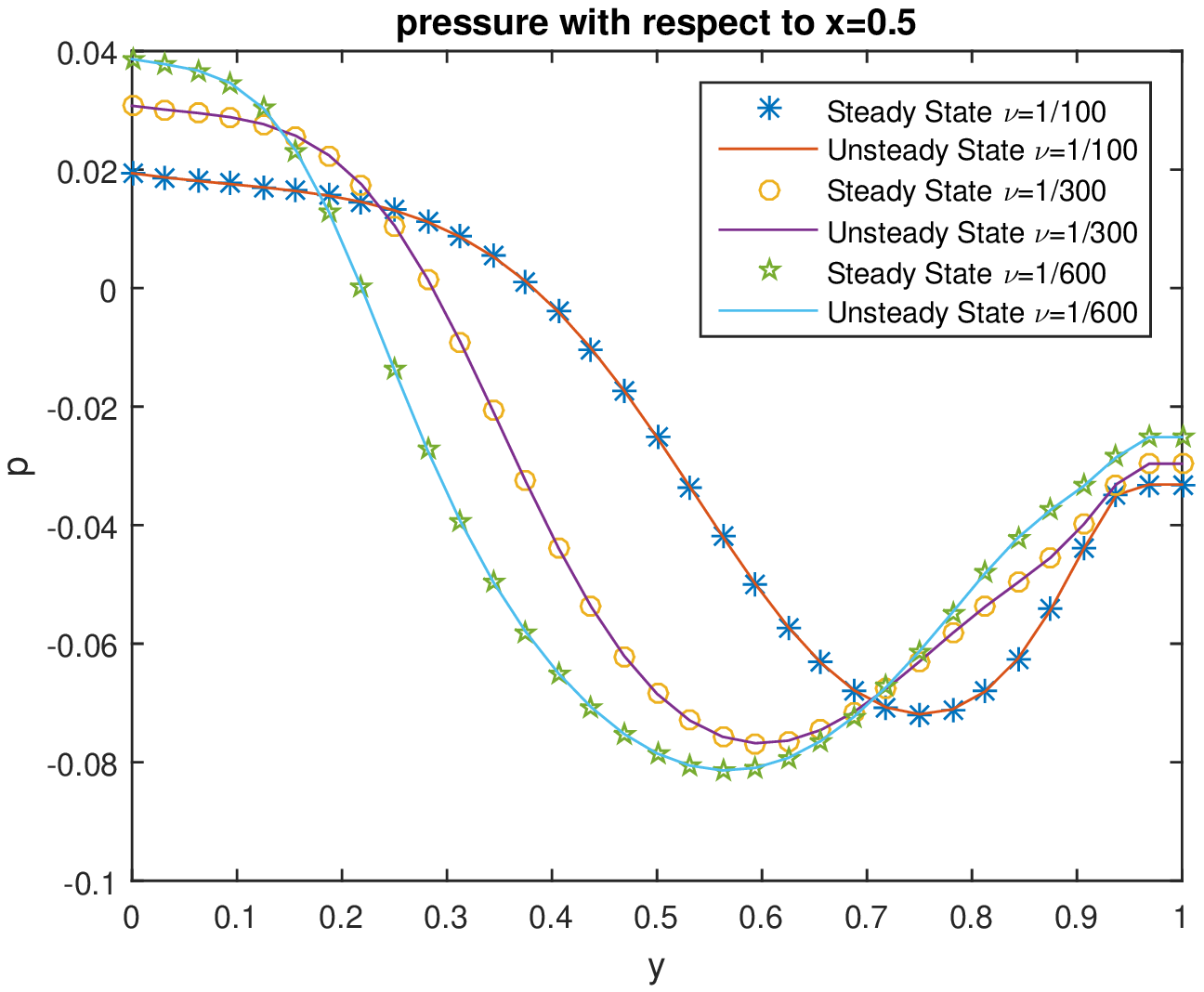}
	\caption{Pressure for Example \ref{ex3}.}
	\label{fig4}
\end{figure}

\begin{example}\label{ex3}
	(2D Lid Driven Cavity Flow Benchmark Problem). In this example, we consider a benchmark problem related to a lid driven cavity flow on a unit square with zero body forces. Further, no slip boundary conditions are considered everywhere, except non zero velocity $(u_1,u_2)=(1,0)$ on the upper part of boundary, that is, the lid of the cavity is moving horizontally with a prescribed velocity. For numerical experiments, we have chosen lines $(0.5,y$) and $(x,0.5)$. We choose $P_1-P_0$ mixed finite element space for the space discretization. In Figure \ref{fig3}, we have
	presented the comparison between unstedy  backward Euler and steady state velocities, whereas Figure \ref{fig4} depicts the comparison of backward Euler and steady state pressure for different values of viscosity $\nu = \{1/100,1/300,1/600\}$, final time $T = 75$, $h = 1/32$ and $\Delta t = \mathcal{O}(h^2)$ with $\sigma_e=40$ and $\kappa=0.1\times\nu$.  From the graphs, it is observed that the Kelvin-Voigt solutions converge to its steady state solutions for large time.
\end{example}
\section{Summary}\label{s8}
\se

In this paper, we have applied the symmetric interior penalty discontinuous Galerkin method to the Kelvin-Voigt equations of motion represented by (\ref{8.1})-(\ref{8.3}). We have defined the semidiscrete discontinuous Galerkin formulation to (\ref{8.1})-(\ref{8.3}) and have derived a priori bounds to the velocity approximation. In order to establish error estimates, we have introduced  the $L^2$-projection $\bP_h$ and a  modified Sobolev-Stokes projection $\bS_h$ on appropriate DG spaces and proved the approximation properties.  Note that, the bounds of $\bP_h$ play a crucial role in deriving the bounds for $\bS_h$. Then, by using regularity estimates for the weak solution and duality arguments along with the approximation properties of $\bP_h$ and $\bS_h$, we have obtained optimal error estimates for the velocity in $L^{\infty}(\bL^2)$ and pressure in $L^{\infty}(L^2)$-norms. Moreover, under the smallness assumption on the data, we have shown that the semidiscrete error estimates are uniform in time. Further, we have employed a backward Euler method for full discretization and have achieved optimal convergence rates for the approximate solution. Finally, we have conducted the numerical experiments and have shown that the outcomes verify the theoretical results.

	\section*{References}
	\begingroup
	\renewcommand{\section}[2]{}%
	
\end{document}